\documentclass[12pt,reqno] {amsart}

 \usepackage{amsmath, wasysym}
 \usepackage{bbm, dsfont}
 \usepackage{amssymb}
\usepackage{graphicx}%
\usepackage{hyperref}
\usepackage{cleveref}
\usepackage{comment}
 \usepackage{mathrsfs}
 \usepackage{float}
 \usepackage{tikz, tikz-cd}
  \usetikzlibrary{arrows.meta, knots, calc}
 \usepackage{bm}
 \usepackage{mathtools} 
\usepackage{extarrows} 
\usepackage{bbm}

 \usepackage{enumerate}
  \usepackage{mathdots}
\usetikzlibrary{decorations.markings}
\usetikzlibrary{decorations.pathreplacing}

\theoremstyle{plain}
\newtheorem{theorem}{Theorem}[section]

\newtheorem{corollary}[theorem]{Corollary}
\newtheorem{lemma}[theorem]{Lemma}
\newtheorem{proposition}[theorem]{Proposition}
\newtheorem*{ac}{Acknowledgement}

\newtheorem{remark}[theorem]{Remark}

\newtheorem{definition}[theorem]{Definition}

\newtheorem{notation}{Notation}[section]

\newcommand{\cM}{\mathcal{M}}
\newcommand{\cN}{\mathcal{N}}
\newcommand{\cL}{\mathcal{L}}
\newcommand{\cJ}{\mathcal{J}}

\newcommand{\cR}{\mathcal{R}}
\newcommand{\cA}{\mathcal{A}}

\newcommand{\bC}{\mathbb{C}}

\newcommand{\bE}{\mathbb{E}}
\newcommand{\bN}{\mathbb{N}}
\newcommand{\bR}{\mathbb{R}}
\newcommand{\bfx}{\mathbf{x}}
\newcommand{\bfy}{\mathbf{y}}
\newcommand{\id}{\mathrm{Id}}
\newcommand{\Ker}{\mathrm{Ker}}
\newcommand{\K}{\mathbf{K}}
\newcommand{\oM}{\widetilde{\mathcal{M}_1^{\oplus}}}

\newcommand{\fF}{\mathfrak{F}}
\newcommand{\Div}{\mathrm{Div}}
\newcommand{\Ran}{\mathrm{Ran }}
\newcommand{\Dom}{\mathrm{Dom }}

\newcommand{\diag}{\mathrm{Diag }}
\newcommand{\Span}{\mathrm{Span }}
\newcommand{\Tr}{\mathrm{Tr}}
\newcommand{\grad}{\mathrm{grad}}
\newcommand{\fix}{\mathscr{M}}
\newcommand{\nfix}{\mathscr{N}}

\begin{document}

\title{Bimodule KMS Symmetric Quantum Markov Semigroups and Gradient Flows}

\author{Chunlan Jiang}
\address{Chunlan Jiang, Hebei Normal University}
\email{cljiang@hebtu.edu.cn }

\author{Jincheng Wan}
\address{Jincheng Wan, Tsinghua University, Beijing}
\email{wanjc23@mails.tsinghua.edu.cn }

\author{Jinsong Wu}
\address{Jinsong Wu, Beijing Institute of Mathematical Sciences and Applications, Beijing, 101408, China}
\email{wjs@bimsa.cn}

\date{}

\begin{abstract}
The bimodule KMS symmetry of a bimodule quantum Markov semigroup extends the classical KMS symmetry of a quantum Markov semigroup. 
Compared with (bimodule) GNS symmetry, the (bimodule) KMS symmetry retains significantly more of the underlying noncommutativity. 
In this paper, we study bimodule KMS symmetric quantum Markov semigroups and introduce directional matrices for such semigroups, which reduce to diagonal matrices in the GNS symmetric setting. 
Using these directional matrices, we establish a corresponding gradient-flow structure. 
As a consequence, we obtain both a modified logarithmic Sobolev inequality and a Talagrand inequality for bimodule KMS symmetric quantum Markov semigroups.
\end{abstract}

\maketitle

\section{Introduction}
The quantum Markov semigroup \cite{Lin76, Fri78, FriVer82} is a fundamental tool in quantum statistical mechanics for modeling open quantum systems. 
In this framework, an open system interacts with a heat bath in thermal equilibrium, which is mathematically represented by an equilibrium state. 
Owing to the noncommutative nature of the setting, the symmetries of the heat flow with respect to the equilibrium state are considerably more intricate than in the classical case. 
Two fundamental examples of such symmetries are the Gelfand–Naimark–Segal (GNS) symmetry \cite{Wirth2022a, wirth2022b} and the Kubo–Martin–Schwinger (KMS) symmetry \cite{FagReb15, FagUma08, FagUma10, KFGV77, AmorimCarlen21, VerWir23}.
In \cite{VerWir23}, Vernooij and Wirth show that the generator of uniformly continuous KMS-symmetric quantum Markov semigroups is the square of a derivation with values in a Hilbert bimodule.

Recently, Zhao and the third author \cite{WuZha25} systematically investigated bimodule quantum Markov semigroups associated with finite quantum symmetries \cite{Jon83, Jon85}, which provides a novel perspective on the study of quantum Markov semigroups. 
They introduced the notions of bimodule GNS symmetry, and bimodule KMS symmetry for bimodule quantum channels and bimodule quantum semigroups. 
It is worth noting that the structures of bimodule GNS symmetry are significantly richer than their classical counterparts. 
Moreover, the corresponding gradient flow and its consequences were also established in \cite{WuZha25}. 
The main difficulty in formulating the gradient flow for bimodule KMS symmetry lies in the noncommutative nature of the underlying duality.

In this paper, we introduce directional matrices for the bimodule KMS symmetric quantum Markov semigroups.
This matrix characterizes the angle between symmetric Laplacian and the bimodule modular operator.
By using this directional matrix, we write the Laplacian for bimodule KMS symmetric quantum Markov semigroups in terms of divergence and gradient.
Consequently, we set up the gradient flow.
Under the intertwining property of the bimodule KMS symmetric quantum Markov semigroups, we obtain the bimodule logarithmic Sobolev inequality and bimodule Talagrand inequality for bimodule KMS symmetric quantum semigroups.

The rest of the paper is organized as follows. 
In Section 2, we recall the Jones tower for $\lambda$-extensions of finite von Neumann algebras and bimodule quantum channels.
In Section 3 we recall the bimodule KMS symmetric bimodule quantum Markov semigroups.
In Section 4 we introduce the directional matrix and related divergence ad gradient. 
In Section 5 we compute the Laplacian of the bimodule KMS symmetric quantum semigroups with respect to the divergence and gradient. 
In Section 6, we setup the gradient flow for the bimodule KMS symmetric quantum semigroups. 
In Section 7, we present one example of bimodule KMS symmetric quantum Markov semigroups which is not KMS symmetric.
In Section 8, we derive the modified logarithmic Sobolev inequality and Talagrand inequality with respect to the hidden density. 
In Section 9, we explore the KMS semigroup for Clifford algebra and obtain a new modified logarithmic Sobolev inequality.

\begin{ac}
J. Wu was supported by grants from Beijing Institute of Mathematical Sciences and Applications.
J.~W. was supported by NSFC (Grant no. 12371124) and partially supported by NSFC (Grant no. 12031004). 
C. ~Jiang was supported by Hebei Natural Science Foundation (Grant No. A2023205045) and NSFC (Grant no. 12471120).
\end{ac}

\section{Preliminaries}
We briefly review the Jones tower of $\lambda$-extensions.
Let $\cN\subset \cM$ be an unital inclusion of finite von Neumann algebras and $\tau$ be a normal faithful tracial state on $\cM$. 
The Hilbert space $L^2(\cM, \tau)$ is the Gelfand-Naimark-Segal Hilbert space of $\tau$, with cyclic separating vector $\Omega$ and modular conjugation $J$ given by $Jx\Omega = x^*\Omega$, \ $x\in\cM$. 
We denote by $e_1$ the Jones projection and $\bE_{\cN}$ the $\tau$-preserving conditional expectation of $\cM$ onto $\cN$.  
The basic construction $\cM_1=\langle \cM, e_1\rangle$ is the von Neumann algebra generated by $\cM$ and $e_1$. 
The inclusion $\cN\subset\cM$ is called finite if $\cM_1$ is a finite von Neumann algebra, and irreducible if $\cN'\cap \cM=\mathbb{C}$.  
For a finite inclusion, we have $\cM_1 = J\cN'J$, where $\cN'$ is the commutant of $\cN$ on $L^2(\cM,\tau)$. 

Suppose $\tau_1$ is a faithful normal trace on $\cM_1$ extending $\tau$, and let $\mathbb{E}_{\cM}$ be the trace-preserving conditional expectation onto $\cM$. 
The pair $(\cM_1, \tau_1)$ is called a $\lambda$-extension of $\cN\subset\cM$ if $\tau_1|_{\cM}=\tau$ and $\bE_{\cM}(e_1)=\lambda$ for some positive constant $\lambda$. 
The index of the extension is defined as $[\cM:\cN] = \lambda^{-1}$. 
We denote by $\Omega_1$ the cyclic and separating vector in $L^2(\cM_1, \tau_1)$, and by $e_2$ the Jones projection onto $L^2(\cM,\tau_1)$. 
The $\lambda$-extension is called extremal if for all $x\in\cN'\cap \cM$, $\tau_1(x) = \tau_1(Jx^*J)$. 
In this paper, $\lambda$-extensions are always assumed to be extremal. 
Let $\mathbb{E}_{\cN'}$ be the $\tau_1$-preserving conditional expectation from $\cM_1$ onto $\cN'\cap \cM_1$. 
Let $\cM_2 = \langle \cM_1,e_2\rangle$ be the basic construction of the inclusion $\cM\subset \cM_1$, with a normal faithful trace $\tau_2$ extending $\tau_1$. 
We assume $\cM_1\subset \cM_2$ is a $\lambda$-extension of $\cM\subset\cM_1$, i.e. $\mathbb{E}_{\cM_1}(e_2) = \lambda$, where $\mathbb{E}_{\cM_1}$ is the $\tau_2$-preserving conditional expectation onto $\cM_1$. 
Denote by $\mathbb{E}_{\cM'}$ the $\tau_2$-preserving conditional expectation from $\cM_2$ onto $\cM'\cap \cM_2$.

For a finite inclusion, there exists a finite set $\{\eta_j\}_{j=1}^m$ of operators in $\cM$ called Pimsner-Popa basis for $\cN\subset \cM$, which satisfies $\displaystyle x=\sum_{j=1}^m  \bE_{\mathcal{N}}(x\eta^*_j)\eta_j$, for all $x\in\cM$\cite{PimPop86}. 
In terms of Jones projection, this condition is expressed as $\displaystyle \sum_{j=1}^m\eta_j^*e_1\eta_j =1$. 
This implies that any operator in $\cM_1$ is a finite sum of operators of the form $ae_1b$ with $a,b\in \cM$. 
As a consequence, for any $y\in \cM_1$, there is a unique $x\in\cM$ such that $ye_1 = xe_1$.
This indicates that $x=\lambda^{-1}\bE_{\cM}(ye_1)$.
That $\cN\subset\cM$ is a $\lambda$-extension implies $\displaystyle \sum_{j=1}^m\eta_j^*\eta_j =\lambda^{-1}$. 
We shall assume that the basis is orthogonal, that is $\bE_{\cN}(\eta_k\eta_j^*)=0$ for $k\neq j$. 
The conditional expectation $\bE^{\cN'}_{\cM'}$ from $\cN' = J\cM_1 J$ onto $\cM' = J\cM J$ can be written as 
\begin{align*}
    \bE^{\cN'}_{\cM'} (x) = \lambda\sum^m_{j=1} \eta^*_j x \eta_j,\quad x\in \cM_2.
\end{align*}
Note that this implies that $\mathbb{E}_{\cM'}(e_1) = \lambda$. 
We also have $\mathbb{E}_{\cM'}(yx) = \mathbb{E}_{\cM'}(xy)$ for all $y\in \cM$ and $x\in \cM_2$. 

The Pimsner-Popa inequality \cite{PimPop86} for the inclusion states that $\bE_{\cN}(x)\geq \lambda_{\cN\subset\cM} x$ for any $0\leq x\in \cM$, where $\lambda_{\cN\subset\cM}$ is the Pimsner-Popa constant. 

The basic construction from a $\lambda$-extension is assumed to be iterated to produce the Jones tower 
\begin{align*}
\cN\subset\cM\subset\cM_1\subset\cM_2\subset\cdots.
\end{align*}
The sequence of higher relative commutants consists of the standard invariant of the inclusion initial. 
The standard invariants are axiomatized by planar algebras in \cite{Jones2021}. 
The Fourier transform $\mathfrak{F}: \cN'\cap \cM_1 \to \cM'\cap \cM_2$ is
\begin{align}\label{eqn:: Fourier transform}
    \mathfrak{F}(x)=&\lambda^{-3/2}\bE_{\cM'}(xe_2e_1), \quad x\in \cN'\cap \cM_1.
\end{align}
Recall that the contragredient $\overline{x}$ is $J x^*J$ for $x\in \cN'\cap \cM_1$.
We have $\fF^2(x)=\overline{x}$ and $(\fF(x))^*=\fF(\overline{x^*})$ for any $x\in \cN'\cap \cM_1$.
The convolution between $x,y\in \cM'\cap \cM_2$ is
\begin{align*}
    x*y =& \mathfrak{F}^{-1}(\mathfrak{F}(y)\mathfrak{F}(x))\\
    =& \lambda^{-9/2}\mathbb{E}_{\cM'}(e_1e_2\mathbb{E}_{\cM_1}(e_2e_1y)\mathbb{E}_{\cM_1}(e_2e_1x)). 
\end{align*}
The shift $\gamma_{1, +}: \cM_1'\cap \cM_3 \to \cN'\cap \cM_1$ is an isomorphism given by
\begin{align*}
\gamma_{1, +}(x)e_3=\lambda^{-2}e_3e_2e_1 x e_1e_2 e_3, \quad x\in \cM_1'\cap \cM_3.
\end{align*}
The inverse $\gamma_{1, +}^{-1}: \cN'\cap \cM_1\to \cM_1'\cap \cM_3$ is given by
\begin{align*}
\gamma_{1, +}^{-1}(x)e_1=\lambda^{-2}e_1e_2 e_3 x e_3e_2e_1, \quad x\in \cN'\cap \cM_1.
\end{align*}
For more details on this string Fourier transform, we refer to \cite{Liu16,JLW16,JLW2019,JJLRW20} etc.

A linear map $\Phi:\cM\rightarrow \cM$ is called positive if it preserves the positive cone $\cM_+$. 
The map $\Phi$ is called completely positive if $\Phi\otimes id_n$ is positive on $\cM\otimes M_n(\mathbb{C})$ for all $n\geq 1$, and completely bounded if $\displaystyle \sup_{n \geq 1}\Vert \Phi\otimes id_n\Vert$ is finite. 
The map $\Phi$ is unital if $\Phi(1)=1$.
A quantum channel is a normal unital completely positive map on $\cM$. 

For a finite inclusion $\cN\subset \cM$ of finite von Neumann algebras, a linear map $\Phi$ is said to be $\cN$-bimodule, if 
\begin{align*}
\Phi(y_1 x y_2) = y_1\Phi(x)y_2,
\end{align*}
for all $y_1, y_2\in\cN$ and $x\in\cM$. 
The Fourier multiplier $\widehat{\Phi}$ is the unique element in $\cM'\cap\cM_2$ such that 
\begin{equation}\label{eqn:: bilinear form induced by Fourier multiplier}
    \langle \widehat{\Phi}(x_1e_1y_1\Omega_1), x_2e_1y_2\Omega_1\rangle = \lambda^{3/2}\tau(y^*_2\Phi(x^*_2x_1)y_1),\quad x_1, x_2,y_1, y_2\in \cM.
\end{equation}
The bimodule map $\Phi$ can be written in terms of the Fourier multiplier $\widehat{\Phi}$ as follows:
\begin{align}\label{eq:multiplierphi}
    \Phi(x)=\lambda^{-5/2} \bE_{\mathcal{M}}(e_2e_1\widehat{\Phi} x e_1e_2).
\end{align}
We shall denote $\lambda^{-5/2} \bE_{\mathcal{M}}(e_2e_1\widehat{\Phi} x e_1e_2)$ by $x*\widehat{\Phi}$ for simplicity.
The limit $\displaystyle \mathbb{E}_{\Phi} = \lim_{\ell\rightarrow \infty}\frac{1}{\ell}\sum^\ell_{k=1}\Phi^{k}$ exists as a $\cN$-bimodule quantum channel, with the property that $\mathbb{E}^2_{\Phi} = \mathbb{E}_{\Phi}$. 
Moreover, the image of $\mathbb{E}_{\Phi}$ is $\fix(\Phi)$. 
Taking Fourier multiplier gives: 
\begin{align*}
    \widehat{\mathbb{E}}_{\Phi} = \lim_{\ell \rightarrow \infty}\frac{1}{\ell}\sum^\ell _{k=1} \widehat{\Phi}^{(*k)},
\end{align*}
with $\widehat{\mathbb{E}}_{\Phi}*\widehat{\mathbb{E}}_{\Phi} =  \widehat{\mathbb{E}}_{\Phi}$. 
Therefore, $\mathbb{E}_{\Phi}$ is an idempotent with positive Fourier multiplier. 
For more details on bimodule quantum channels, we refer to \cite{HJLW23,HJLW24}.

\section{Bimodule KMS Symmetry Semigroups}
In this section, we provide a brief review of the KMS symmetry for quantum channels and the bimodule KMS symmetry for bimodule quantum Markov semigroups.

Suppose that $\rho$ is a faithful normal state on a finite von Neumann algebra $\cM$ with a faithful normal tracial state $\tau$ and $\Phi$ is a quantum channel.
The quantum channel $\Phi$ is KMS symmetry with respect to $\rho$ if for any $x, y\in \cM$,
\begin{align*}
    \langle J\Phi(x)^*J\Delta_\rho^{1/2}\Omega, y\Delta_\rho^{1/2}\Omega\rangle = \langle Jx^*J\Delta_\rho^{1/2}\Omega, \Phi(y)\Delta_\rho^{1/2}\Omega\rangle,
\end{align*}
where $\Delta_\rho$ is the relative modular operator with respect to $\tau$.
By taking $y=1$, we see that $\rho(\Phi(x))=\rho(x)$ for all $x\in \cM$, i.e. $\Phi$ is equilibrium with respect to $\rho$.
Suppose that $e_1\in \Dom(\sigma_{-i/2})$, where $\sigma_t$ is the modular automorphism of $\cM$ arising from $\rho$.
Let 
\begin{align*}
\widehat{\Delta}_{\rho,1/2}
=\lambda^{-1/2}\mathfrak{F}( \bE_{\cN'} (\Delta_{\rho}^{1/2} e_1 \Delta_{\rho}^{-1/2} ) ). 
\end{align*}
By Proposition 4.27 in \cite{WuZha25}, we have that 
\begin{enumerate}[(1)]
    \item $\widehat{\Delta}_{\rho,1/2}\in \cM'\cap \cM_2$ is positive and $\widehat{\Delta}_{\rho,1/2} e_2=e_2$.
    \item $\widehat{\Delta}_{\rho,1/2}$ is invertible.
\end{enumerate}

Suppose $\cN\subset \cM$ is a finite inclusion of finite von Neumann algebras.
A continuous family $\{\Phi_t:\cM\to \cM\}_{t\geq 0}$ of quantum channels is a quantum Markov semigroup if
\begin{enumerate}[(1)]
    \item $\Phi_t\Phi_s=\Phi_{t+s}$ for all $t, s\geq 0$.
    \item $\Phi_0=\id$.
    \item $\Phi_t$ is normal for all $t\geq 0$.
\end{enumerate}
We say $\{\Phi_t\}_{t\geq 0}$ is a bimodule quantum Markov semigroup if $\Phi_t$ is a bimodule quantum channel for $t\geq 0$ with respect to $\cN$.
For more details on bimodule quantum Markov semigroups, we refer to \cite{WuZha25}.

Let $\cL$ be the generator of a bimodule quantum Markov semigroup $\{\Phi_t\}_{t\geq 0}$, i.e. $e^{-t\cL}=\Phi_t$, which is also called Lindbladian.
The Fourier multiplier $\widehat{\cL}$ is related to the Fourier multipliers $\{\widehat{\Phi}_t\}_{t\geq 0}$ as $\displaystyle \widehat{\cL}=\lim_{t\to 0} \frac{\lambda^{-1/2} e_2- \widehat{\Phi}_t}{t}$. 
We are interested in the following two components of the Fourier multiplier $\widehat{\cL}$:
\begin{align*}
    \widehat{\cL}_0= -(1-e_2)\widehat{\cL}(1-e_2)\geq 0. \quad \text{ and } \quad 
    \widehat{\cL}_1 =e_2\widehat{\cL}(1-e_2).
\end{align*}
The generator $\cL$ can be written in terms of $\widehat{\cL}_0$ and $\widehat{\cL}_1$ as follows:
\begin{align}
 \cL(x)=\lambda^{-1/2}\bE_{\cM}(e_2 \widehat{\cL} e_2)x+\lambda^{-3/2} \bE_{\cM}(e_2e_1\widehat{\cL}_1^*) x+ \lambda^{-3/2} x \bE_{\cM}(\widehat{\cL}_1 e_1 e_2)-  x*\widehat{\cL}_0.   
\end{align}
The Laplacian $\cL_a$ of $\{\Phi_t\}_{t\geq 0}$ is 
\begin{align}
    \cL_a(x)= \frac{1}{2}(1*\widehat{\cL}_0) x
        + \frac{1}{2} x  (1*\widehat{\cL}_0)-x* \widehat{\cL}_0, \quad x\in \cM.
\end{align}
Define 
    \begin{align*}
\cL_w (x)
=&  i[x, \Im\bE_{\cM}(\mathfrak{F}^{-1}(\widehat{\cL}_1))], \\
\mathbf{y}_1=\Im\bE_{\cM}(\mathfrak{F}^{-1}(\widehat{\cL}_1))
= & \frac{i}{2}\left(\bE_{\cM}(\mathfrak{F}^{-1}(\widehat{\cL}_1))^*- \bE_{\cM}(\mathfrak{F}^{-1}(\widehat{\cL}_1))\right).
    \end{align*}
    Then $\mathcal{L}$ is decomposed as 
    \begin{align*}
        \cL=\cL_a+ \cL_w.
    \end{align*}  
The gradient form $\Gamma$ of the semigroup is
\begin{align*}
    \Gamma(x, y)=\frac{1}{2} (y^*\cL(x)+\cL(y)^*x-\cL(y^*x)), \quad x, y\in \cM.
\end{align*}
Hence $\Gamma(x, y)=\displaystyle \frac{1}{2} (y^*\cL_a(x)+\cL_a(y)^*x-\cL_a(y^*x))$.
We shall focus on the Laplacian part of a bimodule quantum Markov semigroup.

The fixed point subspace $\fix(\{\Phi_t\}_{t\geq 0})$ is defined to be
\begin{align*}
    \fix(\{\Phi_t\}_{t\geq 0})=\{x\in \cM: \Phi_t(x)=x, t\geq 0\}.
\end{align*}
The multiplicative domain $\nfix(\{\Phi_t\}_{t\geq 0})$ of $\{\Phi_t\}_{t\geq 0}$ is defined to be
\begin{align*}
    \nfix(\{\Phi_t\}_{t\geq 0})=\{x\in \cM: \Phi_t(x^*)\Phi_t(x)=\Phi_t(x^*x) , \Phi_t(x)\Phi_t(x^*)=\Phi_t(xx^*) , t\geq 0\}.
\end{align*}

\begin{definition}[Bimodule KMS Symmetric Semigroups]
Suppose $\{\Phi_t\}_{t\geq 0}$ is a bimodule quantum Markov semigroup and $\widehat{\Delta} \in \cM'\cap \cM_2$ is strictly positive such that $\widehat{\Delta} e_2= e_2$ and $\cR(\widehat{\cL}) \overline{\widehat{\Delta}}= \cR(\widehat{\cL}) \widehat{\Delta}^{-1}$, where $\cR(\widehat{\cL})$ is the range projection of $\widehat{\cL}$.
We say $\{\Phi_t\}_{t\geq 0}$ is bimodule KMS symmetric with respect to $\widehat{\Delta}$ if  $\overline{\widehat{\cL}}
= \overline{\widehat{\Delta}}\widehat{\cL} \overline{\widehat{\Delta}}$ 
and $\Phi_t$ admits a faithful normal equilibrium state for all $t \ge 0$.
\end{definition}

\begin{remark}
The element $\widehat{\Delta}$ generalizes the element $\widehat{\Delta}_{\rho, 1/2}$ when $e_1\in \Dom(\sigma_{-i/2})$.
The requirement that $\Phi_t$ admits a  faithful normal equilibrium state for all $t \ge 0$ ensures that the limit $\displaystyle \lim_{t\to\infty}\Phi_t$ exists. (See Theorem \ref{thm:kmslimit} for details.)
\end{remark}

\begin{remark}
The condition $\cR(\widehat{\cL}) \overline{\widehat{\Delta}}= \cR(\widehat{\cL}) \widehat{\Delta}^{-1}$ is obtained from Theorem 4.25 in \cite{WuZha25}.
It is automatically true for matrix case.
Suppose that $\bC\subset M_n(\bC)$ is the inclusion and $\rho$ is a faithful state on $M_n(\bC)$.
Then $\widehat{\Delta}_{\rho, 1/2}$ always satisfies $\widehat{\Delta}_{\rho, 1/2} e_2=e_2$ and $\overline{\widehat{\Delta}_{\rho, 1/2}}=\widehat{\Delta}_{\rho, 1/2}^{-1}$.
When the semigroup is bimodule GNS symmetric, the condition $\cR(\widehat{\cL}) \overline{\widehat{\Delta}}= \cR(\widehat{\cL}) \widehat{\Delta}^{-1}$ is implied from $\overline{\widehat{\cL}}
= \overline{\widehat{\Delta}}^2\widehat{\cL} $.
\end{remark}

\begin{proposition}
Suppose that $\{\Phi_t\}_{t\geq 0}$ is a bimodule quantum Markov semigroup which is KMS symmetric with respect to a normal faithful state $\rho$ with $e_1\in \Dom(\sigma_{-i/2})$.
Then $\{\Phi_t\}_{t\geq 0}$ is bimodule KMS symmetric with respect to $\widehat{\Delta}_{\rho, 1/2}$.
\end{proposition}
\begin{proof}
By Theorem 4.26 in \cite{WuZha25}, we have that for all $t\geq 0$, 
\begin{align*}
    \overline{\widehat{\Delta}_{\rho, 1/2} }\widehat{\Phi}_t\overline{\widehat{\Delta}_{\rho, 1/2}}
    =\overline{\widehat{\Phi}_t}, \quad \widehat{\Phi}_t\overline{\widehat{\Delta}_{\rho, 1/2}}= \widehat{\Phi}_t \widehat{\Delta}_{\rho, 1/2}^{-1}.
\end{align*}
Note that $\widehat{\Delta}_{\rho, 1/2}e_2=e_2$.
By taking differentiation at $t=0$, we obtain that $\overline{\widehat{\Delta}_{\rho, 1/2} }\widehat{\cL}\overline{\widehat{\Delta}_{\rho, 1/2}}
    =\overline{\widehat{\cL}}$ and $\widehat{\cL} \overline{\widehat{\Delta}_{\rho, 1/2} }= \widehat{\cL}\widehat{\Delta}_{\rho, 1/2}^{-1}$.
Finally, by the fact that $\Phi_t$ is in equilibrium with respect to $\rho$, we see that $\Phi$ is bimodule KMS symmetric with respect to $\widehat{\Delta}_{\rho, 1/2}$.
\end{proof}

\begin{proposition}\label{prop:bimoduleder}
Suppose $\{\Phi_t\}_{t\geq 0}$ is bimodule KMS symmetric with respect to $\widehat{\Delta}$.
Then
\begin{align*}
\widehat{\cL}_0 \overline{\widehat{\Delta}}= \widehat{\cL}_0\widehat{\Delta}^{-1}, \quad \overline{\widehat{\cL}_0} 
=\overline{\widehat{\Delta}}\widehat{\cL}_0 \overline{\widehat{\Delta}},
\quad 
\overline{\overline{\widehat{\Delta}^{1/2}}\widehat{\cL}_0\overline{\widehat{\Delta}^{1/2}}}  = \widehat{\Delta}^{-1/2} \widehat{\cL}_0 \widehat{\Delta}^{-1/2}.
\end{align*}
and 
\begin{align*}
\overline{\widehat{\cL}_1}^*=\widehat{\cL}_1\overline{\widehat{\Delta}}.  
\end{align*}
\end{proposition}
\begin{proof}
Multiplying $1-e_2$ from the left and right hand sides to the equalities $\widehat{\cL}\overline{\widehat{\Delta}}= \widehat{\cL}\widehat{\Delta}^{-1}$ and $\overline{\widehat{\cL}}
= \overline{\widehat{\Delta}}\widehat{\cL} \overline{\widehat{\Delta}}$, we obtain that $\widehat{\cL}_0 \overline{\widehat{\Delta}}= \widehat{\cL}_0\widehat{\Delta}^{-1}$ and $\overline{\widehat{\cL}_0}
= \overline{\widehat{\Delta}}\widehat{\cL}_0 \overline{\widehat{\Delta}}$.
Multiplying $e_2$ from the left hand side and $1-e_2$ from the right hand side to $\overline{\widehat{\cL}}
= \overline{\widehat{\Delta}}\widehat{\cL} \overline{\widehat{\Delta}}$, we obtain that $\overline{\widehat{\cL}_1}^*=\widehat{\cL}_1\overline{\widehat{\Delta}}$.

Note that 
\begin{align*}
\overline{\overline{\widehat{\Delta}^{1/2}}\widehat{\cL}_0\overline{\widehat{\Delta}^{1/2}}}  =  \widehat{\Delta}^{1/2}\overline{\widehat{\cL}_0} \widehat{\Delta}^{1/2}
=\widehat{\Delta}^{1/2} \overline{\widehat{\Delta}}\widehat{\cL}_0 \overline{\widehat{\Delta}}\widehat{\Delta}^{1/2}=\widehat{\Delta}^{-1/2} \widehat{\cL}_0 \widehat{\Delta}^{-1/2}. 
\end{align*}
This completes the computation.
\end{proof}

\begin{proposition}\label{prop:kms1}
Suppose that $\{\Phi_t\}_{t\geq 0}$ is bimodule quantum Markov semigroup such that $\Phi_t$ is equilibrium for some faithful normal state for all $t \ge 0$ and $0<\widehat{\Delta}\in \cM'\cap \cM_2$ with $\widehat{\Delta}e_2=e_2$ and $ \overline{\widehat{\Delta}}= \widehat{\Delta}^{-1}$.
Then $\{\Phi_t\}_{t\geq 0}$ is bimodule KMS symmetric with respect to $\widehat{\Delta}$ if and only if $\overline{\widehat{\cL}_0} 
=\overline{\widehat{\Delta}}\widehat{\cL}_0 \overline{\widehat{\Delta}}$ and 
\begin{align*}
 \vcenter{\hbox{\begin{tikzpicture}[scale=0.65]
\draw [blue] (-0.5, -1.5)--(-0.5, 0) .. controls +(0, 0.6) and +(0,0.6).. (0.5, 0)--(0.5, -1.5);
\begin{scope}[shift={(0.5, -0.3)}]
\draw [fill=white] (-0.3, -0.3) rectangle (0.3, 0.3);
\node at (0, 0) {\tiny $\mathbf{y}^*$};
\end{scope}
\begin{scope}[shift={(0, -1)}]
\draw [fill=white] (-0.7, -0.3) rectangle (0.7, 0.3);
\node at (0, 0) {\tiny $\overline{\widehat{\Delta}}$};    
\end{scope}
\end{tikzpicture}}}= \vcenter{\hbox{\begin{tikzpicture}[scale=0.65]
\draw [blue] (-0.5, -0.8)--(-0.5, 0) .. controls +(0, 0.6) and +(0,0.6).. (0.5, 0)--(0.5, -0.8);
\begin{scope}[shift={(0.5, -0.3)}]
\draw [fill=white] (-0.3, -0.3) rectangle (0.3, 0.3);
\node at (0, 0) {\tiny $\mathbf{y}$};
\end{scope}
\end{tikzpicture}}},    
\end{align*}
where $\displaystyle \mathbf{y}=\frac{1}{2}( 1*\widehat{\cL}_0)+i \Im \bE_{\cM}(\fF^{-1}(\widehat{\cL}_1))$.
\end{proposition}
\begin{proof}
Note that     
\begin{align*}
    \widehat{\cL}=\vcenter{\hbox{\begin{tikzpicture}[scale=0.65]
    \begin{scope}[shift={(0,1.5)}]
    \draw [blue] (-0.5, 0.8)--(-0.5, 0) .. controls +(0, -0.6) and +(0,-0.6).. (0.5, 0)--(0.5, 0.8);    
\begin{scope}[shift={(0.5, 0.3)}]
\end{scope}
    \end{scope}
\draw [blue] (-0.5, -0.8)--(-0.5, 0) .. controls +(0, 0.6) and +(0,0.6).. (0.5, 0)--(0.5, -0.8);
\begin{scope}[shift={(0.5, -0.3)}]
\draw [fill=white] (-0.3, -0.3) rectangle (0.3, 0.3);
\node at (0, 0) {\tiny $\mathbf{y}^*$};
\end{scope}
\end{tikzpicture}}}
+
\vcenter{\hbox{\begin{tikzpicture}[scale=0.65]
    \begin{scope}[shift={(0,1.5)}]
    \draw [blue] (-0.5, 0.8)--(-0.5, 0) .. controls +(0, -0.6) and +(0,-0.6).. (0.5, 0)--(0.5, 0.8);    
\begin{scope}[shift={(0.5, 0.3)}]
\draw [fill=white] (-0.3, -0.3) rectangle (0.3, 0.3);
\node at (0, 0) {\tiny $\mathbf{y}$};
\end{scope}
    \end{scope}
\draw [blue] (-0.5, -0.8)--(-0.5, 0) .. controls +(0, 0.6) and +(0,0.6).. (0.5, 0)--(0.5, -0.8);
\begin{scope}[shift={(0.5, -0.3)}]
\end{scope}
\end{tikzpicture}}}
-  \vcenter{\hbox{\begin{tikzpicture}[scale=0.6]
        \draw [blue] (0.2, -0.8) --(0.2, 0.8);
        \draw [blue] (-0.2, -0.8) --(-0.2, 0.8);
        \draw [fill=white] (-0.5, -0.4) rectangle (0.5, 0.4);
    \node at (0, 0) {\tiny $\widehat{\cL}_0$};
        \end{tikzpicture}}}.
\end{align*} 

Suppose that $\{\Phi_t\}_{t\geq 0}$ is bimodule KMS symmetric.
The assumption $\overline{\widehat{\cL}}
= \overline{\widehat{\Delta}}\widehat{\cL} \overline{\widehat{\Delta}}$ implies that 
\begin{align*}
  \vcenter{\hbox{\begin{tikzpicture}[scale=0.65]
    \begin{scope}[shift={(0,1.5)}]
    \draw [blue] (-0.5, 0.5)--(-0.5, 0) .. controls +(0, -0.6) and +(0,-0.6).. (0.5, 0)--(0.5, 0.5);    
\begin{scope}[shift={(0.5, 0.3)}]
\end{scope}
    \end{scope}
\draw [blue] (-0.5, -1.5)--(-0.5, 0) .. controls +(0, 0.6) and +(0,0.6).. (0.5, 0)--(0.5, -1.5);
\begin{scope}[shift={(0.5, -0.3)}]
\draw [fill=white] (-0.3, -0.3) rectangle (0.3, 0.3);
\node at (0, 0) {\tiny $\mathbf{y}^*$};
\end{scope}
\begin{scope}[shift={(0, -1)}]
\draw [fill=white] (-0.7, -0.3) rectangle (0.7, 0.3);
\node at (0, 0) {\tiny $\overline{\widehat{\Delta}}$};    
\end{scope}
\end{tikzpicture}}}
+
\vcenter{\hbox{\begin{tikzpicture}[scale=0.65]
    \begin{scope}[shift={(0,1.5)}]
    \draw [blue] (-0.5, 1.5)--(-0.5, 0) .. controls +(0, -0.6) and +(0,-0.6).. (0.5, 0)--(0.5, 1.5);    
\begin{scope}[shift={(0.5, 0.3)}]
\draw [fill=white] (-0.3, -0.3) rectangle (0.3, 0.3);
\node at (0, 0) {\tiny $\mathbf{y}$};
\end{scope}
\begin{scope}[shift={(0, 1)}]
\draw [fill=white] (-0.7, -0.3) rectangle (0.7, 0.3);
\node at (0, 0) {\tiny $\overline{\widehat{\Delta}}$};    
\end{scope}
    \end{scope}
\draw [blue] (-0.5, -0.5)--(-0.5, 0) .. controls +(0, 0.6) and +(0,0.6).. (0.5, 0)--(0.5, -0.5);
\begin{scope}[shift={(0.5, -0.3)}]
\end{scope}
\end{tikzpicture}}}
=
\vcenter{\hbox{\begin{tikzpicture}[scale=0.65]
    \begin{scope}[shift={(0,1.5)}]
    \draw [blue] (-0.5, 0.8)--(-0.5, 0) .. controls +(0, -0.6) and +(0,-0.6).. (0.5, 0)--(0.5, 0.8);    
\begin{scope}[shift={(0.5, 0.3)}]
\end{scope}
    \end{scope}
\draw [blue] (-0.5, -0.8)--(-0.5, 0) .. controls +(0, 0.6) and +(0,0.6).. (0.5, 0)--(0.5, -0.8);
\begin{scope}[shift={(0.5, -0.3)}]
\draw [fill=white] (-0.3, -0.3) rectangle (0.3, 0.3);
\node at (0, 0) {\tiny $\mathbf{y}$};
\end{scope}
\end{tikzpicture}}}
+
\vcenter{\hbox{\begin{tikzpicture}[scale=0.65]
    \begin{scope}[shift={(0,1.5)}]
    \draw [blue] (-0.5, 0.8)--(-0.5, 0) .. controls +(0, -0.6) and +(0,-0.6).. (0.5, 0)--(0.5, 0.8);    
\begin{scope}[shift={(0.5, 0.3)}]
\draw [fill=white] (-0.3, -0.3) rectangle (0.3, 0.3);
\node at (0, 0) {\tiny $\mathbf{y}^*$};
\end{scope}
    \end{scope}
\draw [blue] (-0.5, -0.8)--(-0.5, 0) .. controls +(0, 0.6) and +(0,0.6).. (0.5, 0)--(0.5, -0.8);
\begin{scope}[shift={(0.5, -0.3)}]
\end{scope}
\end{tikzpicture}}}.
\end{align*}
By taking the Fourier transform, we have that 
\begin{align*}
\vcenter{\hbox{\begin{tikzpicture}[scale=0.65]
\draw [blue] (1, -1.8)--(1, 0.7);
\draw [blue] (-1.2, 0.7)--(-1.2, -1.3)..controls+(0, -0.5) and +(0, -0.5)..(-0.5, -1.3)--(-0.5, 0) .. controls +(0, 0.6) and +(0,0.6).. (0.5, 0)--(0.5, -1.8);
\begin{scope}[shift={(0.5, -0.3)}]
\draw [fill=white] (-0.3, -0.3) rectangle (0.3, 0.3);
\node at (0, 0) {\tiny $\mathbf{y}^*$};
\end{scope}
\begin{scope}[shift={(0, -1)}]
\draw [fill=white] (-0.7, -0.3) rectangle (0.7, 0.3);
\node at (0, 0) {\tiny $\overline{\widehat{\Delta}}$};    
\end{scope}
\end{tikzpicture}}}
+
\vcenter{\hbox{\begin{tikzpicture}[scale=0.65]
    \begin{scope}[shift={(0,1.5)}]
    \draw [blue] (-1, 1.8)--(-1, -0.7);
    \draw [blue] (-0.5, 1.8)--(-0.5, 0) .. controls +(0, -0.6) and +(0,-0.6).. (0.5, 0)--(0.5, 1.3)..controls +(0, 0.5) and +(0, 0.5).. (1.2,1.3 )--(1.2, -0.7);    
\begin{scope}[shift={(0.5, 0.3)}]
\draw [fill=white] (-0.3, -0.3) rectangle (0.3, 0.3);
\node at (0, 0) {\tiny $\mathbf{y}$};
\end{scope}
\begin{scope}[shift={(0, 1)}]
\draw [fill=white] (-0.7, -0.3) rectangle (0.7, 0.3);
\node at (0, 0) {\tiny $\overline{\widehat{\Delta}}$};    
\end{scope}
    \end{scope}
\end{tikzpicture}}}
=
\vcenter{\hbox{\begin{tikzpicture}[scale=0.65]
    \draw [blue] (-1, 1.2)--(-1, -1.2) (-0.3, 1.2)--(-0.3, -1.2);
    \begin{scope}[shift={(-1, 0)}]
\draw [fill=white] (-0.3, -0.3) rectangle (0.3, 0.3);
\node at (0, 0) {\tiny $\mathbf{y}$};
\end{scope}
\end{tikzpicture}}}
+\vcenter{\hbox{\begin{tikzpicture}[scale=0.65]
    \draw [blue] (-1, 1.2)--(-1, -1.2) (-0.3, 1.2)--(-0.3, -1.2);
    \begin{scope}[shift={(-0.3, 0)}]
\draw [fill=white] (-0.3, -0.3) rectangle (0.3, 0.3);
\node at (0, 0) {\tiny $\overline{\mathbf{y}}^*$};
\end{scope}
\end{tikzpicture}}}.
\end{align*}
Taking the conditional expectation and by the fact that $\tau(\Im \bE_{\cM}(\fF^{-1}(\widehat{\cL}_1))) =0$, we obtain that $  \vcenter{\hbox{\begin{tikzpicture}[scale=0.65]
\draw [blue] (-0.5, -1.5)--(-0.5, 0) .. controls +(0, 0.6) and +(0,0.6).. (0.5, 0)--(0.5, -1.5);
\begin{scope}[shift={(0.5, -0.3)}]
\draw [fill=white] (-0.3, -0.3) rectangle (0.3, 0.3);
\node at (0, 0) {\tiny $\mathbf{y}^*$};
\end{scope}
\begin{scope}[shift={(0, -1)}]
\draw [fill=white] (-0.7, -0.3) rectangle (0.7, 0.3);
\node at (0, 0) {\tiny $\overline{\widehat{\Delta}}$};    
\end{scope}
\end{tikzpicture}}}= \vcenter{\hbox{\begin{tikzpicture}[scale=0.65]
\draw [blue] (-0.5, -0.8)--(-0.5, 0) .. controls +(0, 0.6) and +(0,0.6).. (0.5, 0)--(0.5, -0.8);
\begin{scope}[shift={(0.5, -0.3)}]
\draw [fill=white] (-0.3, -0.3) rectangle (0.3, 0.3);
\node at (0, 0) {\tiny $\mathbf{y}$};
\end{scope}
\end{tikzpicture}}}$.

Suppose that $\overline{\widehat{\cL}_0} 
=\overline{\widehat{\Delta}}\widehat{\cL}_0 \overline{\widehat{\Delta}}$ and $ \vcenter{\hbox{\begin{tikzpicture}[scale=0.65]
\draw [blue] (-0.5, -1.5)--(-0.5, 0) .. controls +(0, 0.6) and +(0,0.6).. (0.5, 0)--(0.5, -1.5);
\begin{scope}[shift={(0.5, -0.3)}]
\draw [fill=white] (-0.3, -0.3) rectangle (0.3, 0.3);
\node at (0, 0) {\tiny $\mathbf{y}^*$};
\end{scope}
\begin{scope}[shift={(0, -1)}]
\draw [fill=white] (-0.7, -0.3) rectangle (0.7, 0.3);
\node at (0, 0) {\tiny $\overline{\widehat{\Delta}}$};    
\end{scope}
\end{tikzpicture}}}= \vcenter{\hbox{\begin{tikzpicture}[scale=0.65]
\draw [blue] (-0.5, -0.8)--(-0.5, 0) .. controls +(0, 0.6) and +(0,0.6).. (0.5, 0)--(0.5, -0.8);
\begin{scope}[shift={(0.5, -0.3)}]
\draw [fill=white] (-0.3, -0.3) rectangle (0.3, 0.3);
\node at (0, 0) {\tiny $\mathbf{y}$};
\end{scope}
\end{tikzpicture}}}$.
We see that $\overline{\widehat{\cL}} 
=\overline{\widehat{\Delta}}\widehat{\cL} \overline{\widehat{\Delta}}$.
This implies that $\{\Phi_t\}_{t\geq 0}$ is bimodule KMS symmetric.
\end{proof}

\begin{corollary}\label{cor:kms2}
Suppose that $\{\Phi_t\}_{t\geq 0}$ is a bimodule quantum Markov semigroup such that $\Phi_t$ admits a faithful normal equilibrium state for all $t \ge 0$ and $1*\widehat{\cL}_0$ is a multiple of the identity.
Suppose $0<\widehat{\Delta}\in \cM'\cap \cM_2$ with $\widehat{\Delta}e_2=e_2$ and $ \overline{\widehat{\Delta}}= \widehat{\Delta}^{-1}$.
Then $\{\Phi_t\}_{t\geq 0}$ is bimodule KMS symmetric with respect to $\widehat{\Delta}$ if and only if $\overline{\widehat{\cL}_0} 
=\overline{\widehat{\Delta}}\widehat{\cL}_0 \overline{\widehat{\Delta}}$ and 
\begin{align*}
 \vcenter{\hbox{\begin{tikzpicture}[scale=0.65]
\draw [blue] (-0.5, -1.5)--(-0.5, 0) .. controls +(0, 0.6) and +(0,0.6).. (0.5, 0)--(0.5, -1.5);
\begin{scope}[shift={(0.5, -0.3)}]
\draw [fill=white] (-0.3, -0.3) rectangle (0.3, 0.3);
\node at (0, 0) {\tiny $\mathbf{y}_1$};
\end{scope}
\begin{scope}[shift={(0, -1)}]
\draw [fill=white] (-0.7, -0.3) rectangle (0.7, 0.3);
\node at (0, 0) {\tiny $\overline{\widehat{\Delta}}$};    
\end{scope}
\end{tikzpicture}}}
=- \ \vcenter{\hbox{\begin{tikzpicture}[scale=0.65]
\draw [blue] (-0.5, -0.8)--(-0.5, 0) .. controls +(0, 0.6) and +(0,0.6).. (0.5, 0)--(0.5, -0.8);
\begin{scope}[shift={(0.5, -0.3)}]
\draw [fill=white] (-0.3, -0.3) rectangle (0.3, 0.3);
\node at (0, 0) {\tiny $\mathbf{y}_1$};
\end{scope}
\end{tikzpicture}}},    
\end{align*}
where $\displaystyle \mathbf{y}_1= \Im \bE_{\cM}(\fF^{-1}(\widehat{\cL}_1))$.    
\end{corollary}
\begin{proof}
It follows from Proposition \ref{prop:kms1} by taking $1*\widehat{\cL}_0=\kappa$ for some $\kappa>0$.
\end{proof}

\begin{corollary}
Suppose that $0\leq \widehat{\cL}_0\in \cM'\cap \cM_2$ such that $\widehat{\cL}_0e_2=0$.
Suppose $0< \widehat{\Delta}\in \cM'\cap \cM_2$ with $\widehat{\Delta}e_2=e_2$ and $ \overline{\widehat{\Delta}}= \widehat{\Delta}^{-1}$.
If $1*\widehat{\cL}_0=1$, $q*\overline{\widehat{\cL}_0}=q$ for some $0<q$ in $\cN'\cap \cM$ and $\overline{\widehat{\cL}_0} 
=\overline{\widehat{\Delta}}\widehat{\cL}_0 \overline{\widehat{\Delta}}$.
Then the bimodule quantum Markov semigroup $\{\Phi_t\}_{t\geq 0}$ given by the generator $\cL$ as $\widehat{\cL}=\lambda^{-1/2}e_2-\widehat{\cL}_0$ is bimodule KMS symmetric with respect to $\widehat{\Delta}$.
\end{corollary}
\begin{proof}
By a direct computation, we have that $\overline{\widehat{\cL}} 
=\overline{\widehat{\Delta}}\widehat{\cL} \overline{\widehat{\Delta}}$. 
By the assumption, we obtain that $\Phi_t=e^{t(\cL_0-I)}$, where $\cL_0$ is the quantum channel given by $\widehat{\cL}_0$.
By the assumption that $q * \overline{\widehat{\cL}_0}=q$, we see that $\widehat{\cL} * \overline{q}=(\lambda^{-1/2} e_2- \widehat{\cL}_0)*\overline{q}=0$. 
Then $\Phi_t^*(q)=q$ and $\Phi_t$ admits a common faithful normal equilibrium state for all $t \ge 0$.
Hence $\{\Phi_t\}_{t\geq 0}$ is bimodule KMS symmetric with respect to $\widehat{\Delta}$.
\end{proof}

\begin{remark}
Suppose that $\bC\subset M_n(\bC)$ is the inclusion.  
We shall interpret the example in \cite{CarMaa17} of KMS symmetric semigroup in our framework as follows.
Let
\begin{align*}
    H_0=\sum_{j=1}^m \vcenter{\hbox{\begin{tikzpicture}[scale=0.65]
    \begin{scope}[shift={(0,1.5)}]
    \draw [blue] (-0.5, 0.8)--(-0.5, 0) .. controls +(0, -0.6) and +(0,-0.6).. (0.5, 0)--(0.5, 0.8);    
\begin{scope}[shift={(0.5, 0.3)}]
\end{scope}
    \end{scope}
\draw [blue] (-0.5, -0.8)--(-0.5, 0) .. controls +(0, 0.6) and +(0,0.6).. (0.5, 0)--(0.5, -0.8);
\begin{scope}[shift={(0.5, -0.3)}]
\draw [fill=white] (-0.3, -0.3) rectangle (0.3, 0.3);
\node at (0, 0) {\tiny $v_j^*$};
\end{scope}
\begin{scope}[shift={(0.5, 1.8)}]
\draw [fill=white] (-0.3, -0.3) rectangle (0.3, 0.3);
\node at (0, 0) {\tiny $v_j$};
\end{scope}
\end{tikzpicture}}}
\end{align*}
with $\displaystyle \sum_{j=1}^m v_j^*v_j=1$ and $\displaystyle \sum_{j=1}^m v_j \rho v_j^*=\rho$, where $\rho$ is a density matrix in $M_n(\bC)$, i.e. $1*H_0=1$ and $H_0*\overline{\rho}=\displaystyle \sum_{j=1}^m\overline{v_j^*}\,\overline{\rho} \,\overline{v_j}= \overline{\displaystyle \sum_{j=1}^mv_j \rho v_j^*}=\overline{\rho}$.
Let $\widehat{\Delta}=\vcenter{\hbox{\begin{tikzpicture}[scale=1.3]
        \draw [blue] (0, -0.5)--(0, 0.5) (0.5, -0.5)--(0.5, 0.5);
        \draw [fill=white] (-0.2, -0.2) rectangle (0.2, 0.2);
        \node at (0, 0) {\tiny $\overline{\rho^{\frac{1}{2}}}$};
        \begin{scope}[shift={(0.5, 0)}]
         \draw [fill=white] (-0.2, -0.2) rectangle (0.2, 0.2);
        \node at (0, 0) {\tiny $\rho^{-\frac{1}{2}}$};   
        \end{scope}
    \end{tikzpicture}}}$ and 
\begin{align*}
       \widehat{\cL} =\lambda^{-1/2}e_2- (\widehat{\Delta} \overline{H_0}\widehat{\Delta}) * H_0.
\end{align*}
We see that $\overline{\widehat{\cL}}=\lambda^{-1/2} e_2-\overline{H_0}*(\overline{\widehat{\Delta}} H_0 \overline{\widehat{\Delta}})$.
Note that $1* (\widehat{\Delta} \overline{H_0}\widehat{\Delta})=1$ and $(\widehat{\Delta} \overline{H_0}\widehat{\Delta})*\overline{\rho}=\overline{\rho}$.
Now, we obtain that 
\begin{align*}
    1*\widehat{\cL}=0, \quad \rho*\overline{\widehat{\cL}}=0,
\end{align*}
and 
\begin{align*}
 \overline{\widehat{\Delta}} \left( (\widehat{\Delta} \overline{H_0}\widehat{\Delta}) * H_0 \right) \overline{\widehat{\Delta}} =   \overline{H_0}*(\overline{\widehat{\Delta}} H_0 \overline{\widehat{\Delta}}).
\end{align*}
This shows that the semigroup with the generator $\cL$ is (bimodule) KMS symmetric with respect to $\widehat{\Delta}$.
In this case, we have that 
\begin{align*}
    \widehat{\cL}_0=(1-e_2)\left( (\widehat{\Delta} \overline{H_0}\widehat{\Delta}) * H_0 \right)(1-e_2),
\end{align*}
which is not easy to be read its spectral decomposition.

\end{remark}

\begin{remark}
Suppose that $\cN \subset \cM$ is irreducible and $0\leq H_0=\overline{H_0} \in \cM'\cap \cM_2$ with $H_0 e_2=0$.
Let $0< \widehat{\Delta}\in \cM'\cap \cM_2$ with $\widehat{\Delta}e_2=e_2$ and $\overline{\widehat{\Delta}}=\widehat{\Delta}^{-1}$.
By taking $\widehat{\cL}_0=\widehat{\Delta}^{1/2} H_0\widehat{\Delta}^{1/2}$, we see that the semigroup $\{\Phi_t\}_{t\geq 0}$ with generator $\cL$ given by $\displaystyle \cL(x)=\frac{1}{2}(1*\widehat{\cL}_0)x +\frac{1}{2}x(1*\widehat{\cL}_0)-x* \widehat{\cL}_0$ is bimodule KMS symmetric, where $\Phi_t$ is equilibrium with respect to the trace $\tau$.
In this case, $\cL_w=0$.
\end{remark}

\begin{proposition}
Suppose that $0\leq H_0\in \cM'\cap \cM_2$ such that $H_0e_2=0$ and $\overline{H_0}=H_0$.
Let $0<\widehat{\Delta}\in \cM'\cap \cM_2$ with $\widehat{\Delta}e_2=e_2$ and $\overline{\widehat{\Delta}}=\widehat{\Delta}^{-1}$. 
Then the semigroup with generator $\cL$ with $\widehat{\cL}=\lambda^{-1/2}e_2-\widehat{\Delta}^{1/2} H_0 \widehat{\Delta}^{1/2}$ is bimodule KMS symmetric if $\bE_{\cM_1}(\widehat{\Delta}^{-1/2} H_0 \widehat{\Delta}^{-1/2})=\lambda^{1/2}$ and the associated semigroup is equilibrium with respect to normal faithful state for $t\geq 0$.
\end{proposition}
\begin{proof}
Note that 
\begin{align*}
\overline{\widehat{\Delta}} \widehat{\cL}\overline{\widehat{\Delta}}    
= \lambda^{-1/2}e_2-\overline{\widehat{\Delta}^{1/2}} H_0 \overline{\widehat{\Delta}^{1/2}}
=\overline{\widehat{\cL}}.
\end{align*}
By the assumption $\bE_{\cM_1}(\widehat{\Delta}^{-1/2} H_0 \widehat{\Delta}^{-1/2})=\lambda^{1/2}$, we have that
\begin{align*}
 1*(\widehat{\Delta}^{1/2} H_0 \widehat{\Delta}^{1/2} )=\lambda^{-5/2}\bE_{\cM}(e_2e_1  \widehat{\Delta}^{1/2} H_0 \widehat{\Delta}^{1/2} e_1e_2)=1.
\end{align*}
Hence the semigroup is bimodule KMS symmetric.
\end{proof}

In what follows, we explore the structure of the bimodule KMS symmetric semigroups.
We expect that the balanced Lindbladian $\widehat{\cL_{\Delta}}:=\overline{\widehat{\Delta}^{1/2}} \widehat{\cL}_0 \overline{\widehat{\Delta}^{1/2}}$ is tracially symmetric.
However, we will impose certain conditions to guarantee that this holds.

\begin{proposition}\label{kmsprop:balancedlind}
Suppose that $\{\Phi_t\}_{t\geq 0}$ is bimodule KMS symmetric quantum Markov semigroup with respect to $\widehat{\Delta}\in \cM'\cap \cM_2$.
If one of the following conditions holds:
\begin{enumerate}[(1)]
\item $\cR(\overline{\widehat{\cL}_0})=\cR(\widehat{\cL}_0)$,
\item $\overline{\widehat{\Delta}}=\widehat{\Delta}^{-1}$,
\end{enumerate}
then $\overline{\widehat{\cL}_{\Delta}}= \widehat{\cL}_{\Delta}$.
\end{proposition}
\begin{proof}
Suppose that $(1)$ holds.
We obtain that 
\begin{align*}
\cR(\widehat{\cL}_0^{1/2})=\cR(\widehat{\cL}_0)=\cR(\overline{\widehat{\cL}_0})=\cR(\overline{\widehat{\Delta}}\widehat{\cL}_0\overline{\widehat{\Delta}})=\cR(\overline{\widehat{\Delta}}\widehat{\cL}_0^{1/2}).
\end{align*}
This shows that $\overline{\widehat{\Delta}}\Ran(\widehat{\cL}_0^{1/2})=\Ran(\widehat{\cL}_0^{1/2})$, where $\Ran(x)$ is the range of $x$.
Hence $\cR(\overline{\widehat{\Delta}}^n\widehat{\cL}_0^{1/2})=\cR(\widehat{\cL}_0^{1/2})$ for all $n \geq 1$.
Note that $\cR(\widehat{\cL}_0^{1/2}\overline{\widehat{\Delta}}^n)$ is equivalent to $\cR(\overline{\widehat{\Delta}}^n\widehat{\cL}_0^{1/2})$ and $\cR(\widehat{\cL}_0^{1/2}\overline{\widehat{\Delta}}^n)\leq \cR(\widehat{\cL}_0^{1/2})$.
We see that  $\cR(\widehat{\cL}_0^{1/2}\overline{\widehat{\Delta}}^n)=\cR(\widehat{\cL}_0^{1/2})$ for all $n \geq 1$.
Multiplying $\widehat{\cL}_0^{1/2}$, we have that $\cR(\widehat{\cL}_0\overline{\widehat{\Delta}}^n)=\cR(\widehat{\cL}_0)$ for all $n \geq 1$.
Similarly, we have that $\cR(\widehat{\cL}_0\widehat{\Delta}^{-n})=\cR(\widehat{\cL}_0)$ for all $n \geq 1$.
This implies that 
\begin{align*}
\widehat{\cL}_0 \overline{\widehat{\Delta}}^{n} &=\widehat{\cL}_0\overline{\widehat{\Delta}}^{n-1}_0 \overline{\widehat{\Delta}} =\widehat{\cL}_0\overline{\widehat{\Delta}}^{n-1}_0 \cR(\widehat{\cL}_0) \overline{\widehat{\Delta}}\\
&=\widehat{\cL}_0\overline{\widehat{\Delta}}^{n-1}_0 \cR(\widehat{\cL}_0) \widehat{\Delta}^{-1}=\widehat{\cL}_0\overline{\widehat{\Delta}}^{n-1}_0 \widehat{\Delta}^{-1}\\
&=\widehat{\cL}_0\overline{\widehat{\Delta}}^{n-2}_0 \widehat{\Delta}^{-2}=\cdots=\widehat{\cL}_0 \widehat{\Delta} ^{-n}
\end{align*}
for all $n \geq 1$ and $\widehat{\cL}_0 \overline{\widehat{\Delta}}^{1/2}=\widehat{\cL}_0 \widehat{\Delta} ^{-1/2}$.
By Proposition \ref{prop:bimoduleder}, the operator $\widehat{\cL}_{\Delta}=\overline{\widehat{\Delta}^{1/2}}\widehat{\cL}_0\overline{\widehat{\Delta}^{1/2}}$ satisfies $\overline{\widehat{\cL}_{\Delta}}=\widehat{\cL}_{\Delta}$.

Suppose that $(2)$ holds.
Then $\overline{\widehat{\Delta}}^{1/2}=\widehat{\Delta}^{-1/2}$.
This shows that $\overline{\widehat{\cL}_{\Delta}}=\widehat{\cL}_{\Delta}$.
\end{proof}

In the following, we assume that $\cR(\overline{\widehat{\cL}_0})=\cR(\widehat{\cL}_0)$.
By Proposition \ref{kmsprop:balancedlind}, we see that $\cR(\widehat{\cL}_0)\widehat{\Delta} \cR(\widehat{\cL}_0) =\widehat{\Delta} \cR(\widehat{\cL}_0)$.
This implies that $\widehat{\Delta} \cR(\widehat{\cL}_0)=  \cR(\widehat{\cL}_0)\widehat{\Delta}$.
Let $\cA$ be the von Neumann subalgebra of $\cM'\cap \cM_2$ generated by $\widehat{\cL}_{\Delta}$ and $\widehat{\Delta}\cR(\widehat{\cL}_0)$.
Since $\mathcal{M}' \cap \mathcal{M}_2$ is finite-dimensional, $\mathcal{A}$ decomposes into a direct sum of matrix algebras:
\begin{align*}
\cA=\bigoplus_{\ell=1}^m M^{(\ell)},
\end{align*}
where $M^{(\ell)}$ are matrix algebras and there exists an involution $*$ on the index set $\{1, \ldots, m\}$ by the fact that $\cA$ is invariant under the contragredient.
This is due to the fact that $\overline{\widehat{\cL}_\Delta}=\widehat{\cL}_\Delta$ and 
\begin{align*}
\overline{\widehat{\Delta}\cR(\widehat{\cL}_0)}=\overline{\cR(\widehat{\cL}_0)} \, \overline{\widehat{\Delta}}=\cR(\overline{\widehat{\cL}_0}) \widehat{\Delta}^{-1}=\cR(\widehat{\cL}_0) \widehat{\Delta}^{-1} .
\end{align*}
Moreover, $M^{(\ell)}$ is isomorphic to $M^{(\ell^*)}$ for all $\ell=1, \ldots, m$.

\begin{lemma}\label{lem:con1}
Suppose that $E, F\in \cM'\cap \cM_2$ are minimal projections such that $\overline{E}=E$, $\overline{F}=F$ and $EF=0$.
If there exists a partial isometry $V\in \cM'\cap \cM_2$ such that $V^*V=E$ and $VV^*=F$, then $V$ can be taken satisfying $V^*=\overline{V}$.
\end{lemma}
\begin{proof}
Let $V_0\in \cM'\cap \cM_2$ such that $V_0^*V_0=E$ and $V_0V_0^*=F$.
Then $\overline{V_0}\ \overline{V_0}^*=E$,  $\overline{V_0}^* \ \overline{V_0}=F$.
Note that $E, F$ are minimal projections.
We see that $V_0\overline{V_0}=e^{i\theta} F$ for some $\theta\in \bR$.
Hence we can replace $V_0$ by $e^{i\theta/2}V_0$ and assume that $V_0\overline{V_0}=F$.
In this case, we have that $\overline{V_0} V_0=\overline{V_0} \ \overline{V_0}^* V_0^*V_0=E$.
Let $\displaystyle V=\frac{1}{2}( V_0+\overline{V_0}^*)$.
We have that 
\begin{align*}
VV^*=\frac{1}{4}( V_0V_0^* +V_0\overline{V_0} + \overline{V_0}^*V_0^* + \overline{V_0}^* \overline{V_0})=F.
\end{align*}
Similarly, we have that $VV^*=E$.
Note that $V^*=\overline{V}$.
We see that $V$ is the partial isometry satisfying the requirement.
\end{proof}

\begin{lemma}\label{lem:con2}
Suppose that $E\in \cM'\cap \cM_2$ is minimal projections such that $\overline{E}E=0$.
If there exists a partial isometry $V\in \cM\cap \cM_2$ such that $V^*V=E$ and $VV^*=\overline{E}$,  then there exist minimal projections $E_0, F_0\in \cM'\cap \cM_2$ such that $\overline{E_0}=E_0$, $\overline{F_0}=F_0$ and $E+\overline{E}=E_0+F_0$.
\end{lemma}
\begin{proof}
Let $V$ be a partial isometry such that $V^*V=E$ and $VV^*=\overline{E}$.
We see that $V\overline{V}^*=\overline{E}$ and $V^*=\overline{V}^*$.
Let $\displaystyle E_0=\frac{1}{2} (E+\overline{E}+V +V^*)$.
We obtain that $E_0$ is a minimal projection and $\overline{E_0}=E_0$.
Similarly, we take $\displaystyle F_0=\frac{1}{2} (E+\overline{E}-V -V^*)$ and $F_0$ is a minimal projection with $\overline{F_0}=F_0$.
\end{proof}

By Lemmas \ref{lem:con1}, \ref{lem:con2} and the fact that $\overline{\widehat{\cL}_{\Delta}}=\widehat{\cL}_{\Delta}$, we can choose   a system $\{E_{j,k}^{(\ell)}\}_{j,k=1}^{d_\ell}$ of matrix units of $M^{(\ell)}$ such that
$\widehat{\cL}_{\Delta}$ is diagonal with respect to the system of matrix units, where $d_\ell$ is the size of the matrix algebra $M^{(\ell)}$ and $\overline{E_{j,k}^{(\ell)}}=E_{k, j}^{(\ell^*)}$, for $j,k=1, \ldots, d_{\ell}$.

Let
\begin{align*}
\widehat{\cL}_{\Delta}=& \sum_{\ell=1}^m \sum_{j=1}^{d_\ell} \omega_j^{(\ell)} E_{j,j}^{(\ell)}, 
\end{align*}
where $\omega_j^{(\ell)} \geq 0$ for all $j=1, \ldots, d_\ell$ and $\ell=1, \ldots, m$.
By the fact that $\overline{\widehat{\cL}_{\Delta}}=\widehat{\cL}_{\Delta}$, we have that $\omega_j^{(\ell^*)}=\omega_j^{(\ell)}$ for $j=1, \ldots, d_\ell$.
By the assumption that $\cR(\overline{\widehat{\cL}_0})=\cR(\widehat{\cL}_0)$,
we have that $\cR(\widehat{\cL}_{\Delta})=\cR(\widehat{\Delta}\widehat{\cL}_{\Delta} \widehat{\Delta})$, i.e. $\cR(\widehat{\cL}_{\Delta})=\cR(\widehat{\Delta}\widehat{\cL}_{\Delta})$.
This implies that $\cR(\widehat{\cL}_{\Delta})=\cR(\widehat{\Delta}^s\widehat{\cL}_{\Delta} \widehat{\Delta}^s)$ for $s\in \bR$.
Hence $\widehat{\cL}_\Delta$ is invertible in $\cA$.
Let $B_\ell =\diag(\omega_1^{(\ell)}, \ldots ,\omega_{d_{\ell}}^{(\ell)})\in M_{d_\ell}(\bC)$ for $\ell=1, \ldots, m$, which is a diagonal matrix.
We see that $B_\ell$ is invertible and $B_{\ell^*}=B_{\ell}$.

Let
\begin{align*}
    \widehat{\Delta}^{1/2}\cR(\widehat{\cL}_0)=& \sum_{\ell=1}^m \sum_{j,k=1}^{d_\ell} a_{jk}^{(\ell)} E_{j,k}^{(\ell)}=\sum_{\ell=1}^m \widehat{\Delta}^{(\ell)},
\end{align*}
where $\widehat{\Delta}^{(\ell)}\in M^{(\ell)}$.
By the fact that $\cR(\widehat{\cL}_0)\overline{\widehat{\Delta}}=\cR(\widehat{\cL}_0)\widehat{\Delta}^{-1}$, we have that $\overline{\widehat{\Delta}^{(\ell)}}= \widehat{\Delta}^{(\ell^*)-1}$, i.e.
\begin{align*}
(a_{kj}^{(\ell^*)})_{j,k=1}^{d_\ell} = \left( (a_{jk}^{(\ell)})_{j,k=1}^{d_\ell}\right)^{-1}\in M_{d_{\ell}}(\bC).
\end{align*}

For each $\ell$ with $\ell^*\neq \ell$, we can choose a unitary element $U^{(\ell)} \in M^{(\ell)}$ such that $\widehat{\Delta}^{(\ell)}$ is diagonal with respect to the system $\{U^{(\ell)} E_{j,k}^{(\ell)} U^{(\ell)*}\}_{j,k=1}^{d_\ell}$ of matrix units and $U^{(\ell^*)}=\overline{U^{(\ell)*}}$.
We have that 
\begin{align*}
\overline{F_{j,k}^{(\ell)}}=\overline{U^{(\ell)} E_{j,k}^{(\ell)} U^{(\ell)*}} 
=\overline{U^{(\ell)*}}\  \overline{ E_{j,k}^{(\ell)}}\ \overline{ U^{(\ell)}}
=U^{(\ell^*)} E_{k,j}^{(\ell^*)} U^{(\ell^*)*}=F_{k,j}^{(\ell^*)}.
\end{align*}
Let 
\begin{align*}
U^{(\ell)*}\widehat{\Delta}^{(\ell)} U^{(\ell)}=\sum_{j=1}^{d_\ell}\mu_j^{(\ell)}E_{j,j}^{(\ell)}.
\end{align*}
Then  $\mu_j^{(\ell^*)}={\mu_j^{(\ell)}}^{-1}$ since $\overline{\widehat{\Delta}^{(\ell)}}= \widehat{\Delta}^{(\ell^*)-1}$.
Let $\displaystyle U^{(\ell)}= \sum_{j,k=1}^{d_\ell} u_{j,k}^{(\ell)}F_{j,k}^{(\ell)}$ and $u_{j,k}^{(\ell^*)}=\overline{u_{j,k}^{(\ell)}}$, where the $\overline{a}$ is the complex conjugation of the complex number $a$.
We denote by $F_\ell =\diag(\mu_1^{(\ell)}, \ldots ,\mu_{d_{\ell}}^{(\ell)})\in M_{d_\ell}(\bC)$ for $\ell\neq \ell^*$.


For each $\ell$ with $\ell^*=\ell$, we still can choose a unitary element $U^{(\ell)} \in M^{(\ell)}$ such that $\widehat{\Delta}^{(\ell)}$ is diagonal with respect to the system $\left\{U^{(\ell)} E_{j,k}^{(\ell)} U^{(\ell)*}\right\}_{j,k=1}^{d_\ell}$ of matrix units and the set $\{F_{j,j}^{(\ell)}\}_{j=1}^{d_\ell}$ is invariant under the contragredient, where $F_{j,k}^{(\ell)} = U^{(\ell)} E_{j,k}^{(\ell)} U^{(\ell)*}$.
Hence there exists an involution $*$ on the set $\{1, 2, \ldots, d_\ell\}$.
Note that the involution is different from the involution on $\ell$, but they are induced from the contragredient.
We shall use the same notation for convenience.
\begin{lemma}
Suppose that $\ell=\ell^*$.
We can choose $U^{(\ell)}\in M^{(\ell)}$ such that $\overline{F_{j,k}^{(\ell)}}=F_{k^*,j^*}^{(\ell)}$.
\end{lemma}
\begin{proof}
Note that 
\begin{align*}
\overline{F_{j,k}^{(\ell)}}=\overline{F_{j,j}^{(\ell)}F_{j,k}^{(\ell)} F_{k,k}^{(\ell)}}
=\overline{F_{k,k}^{(\ell)}}\  \overline{F_{j,k}^{(\ell)} } \ \overline{F_{j,j}^{(\ell)}}
=F_{k^*,k^*}^{(\ell)} \  \overline{F_{j,k}^{(\ell)} } \  F_{j^*,j^*}^{(\ell)}
= \alpha_{j,k} F_{k^*,j^*}^{(\ell)},
\end{align*}
where $|\alpha_{j,k}|=1$ and $\alpha_{j,j}=1$.
Hence $F_{j,k}^{(\ell)}=\alpha_{j,k} \alpha_{k^*, j^*} F_{j,k}^{(\ell)}$, i.e. $\alpha_{k^*, j^*}=\alpha_{j,k}^{-1}$ and 
\begin{align*}
\alpha_{j,t} \overline{ F_{t^*,j^*}^{(\ell)}} =F_{j,t}^{(\ell)}= F_{j,k}^{(\ell)}F_{k, t}^{(\ell)}=\alpha_{j,k} \alpha_{k,t}\overline{ F_{k^*,j^*}^{(\ell)}}\ \overline{ F_{t^*,k^*}^{(\ell)}}=\alpha_{j,k} \alpha_{k,t} \overline{ F_{t^*,j^*}^{(\ell)}},
\end{align*}
i.e. $\alpha_{j,t}= \alpha_{j,k} \alpha_{k,t} $.
Suppose that $\alpha_{j,k}=e^{i\theta_{j,k}}$ for $\theta_{j,k}\in \bR$ such that $\theta_{j,t}=\theta_{j,k}+\theta_{k,t}$ and $\theta_{k^*, j^*}=-\theta_{j,k}$.
Let $\alpha_{j,k}^{1/2}= e^{i\theta_{j,k}/2}$.
Then $\left\{\alpha_{j,k}^{-1/2} F_{j,k}^{(\ell)}\right\}_{j,k=1}^{d_\ell}$ is a system of matrix units for $M^{(\ell)}$ satisfying $\overline{\alpha_{j,k}^{-1/2} F_{j,k}^{(\ell)}}= \alpha_{j^*,k^*}^{-1/2} F_{j^*,k^*}^{(\ell)}$.
By taking the unitary element $U^{(\ell)}$ in $M^{(\ell)}$ such that $U^{(\ell)} E_{j,k}^{(\ell)} U^{(\ell)*}=\alpha_{j,k}^{-1/2} F_{j,k}^{(\ell)}$ and replacing the system $\{F_{j,k}^{(\ell)}\}_{j,k=1}^{d_\ell}$ by $\left\{\alpha_{j,k}^{-1/2} F_{j,k}^{(\ell)}\right\}_{j,k=1}^{d_\ell}$, we see that the lemma is true.
\end{proof}

Let $\displaystyle V^{(\ell)}=\sum_{j=1}^{d_\ell} F_{j, j^*}^{(\ell)}$ and $\displaystyle U^{(\ell)}= \sum_{j,k=1}^{d_\ell} u_{j,k}^{(\ell)}F_{j,k}^{(\ell)}$.
Then $V^{(\ell)} F_{j,j}^{(\ell)} V^{(\ell)*}=F_{j^*, j^*}^{(\ell)}$ and $V^{(\ell)} F_{j,k}^{(\ell)} V^{(\ell)*}= F_{j^*, k^*}^{(\ell)}$.
This implies that 
\begin{align*}
\overline{F_{k,j}^{(\ell)}}=F_{j^*, k^*}^{(\ell)}=V^{(\ell)} F_{j,k}^{(\ell)} V^{(\ell)*}= V^{(\ell)} U^{(\ell)} E_{j,k}^{(\ell)} U^{(\ell)*} V^{(\ell)*}.
\end{align*}
On the other hand, we have that 
\begin{align*}
\overline{F_{k,j}^{(\ell)}}=\overline{U^{(\ell)} E_{k,j}^{(\ell)} U^{(\ell)*} } =\overline{U^{(\ell)*}} E_{j, k}^{(\ell)} \overline{U^{(\ell)}}.
\end{align*}
Combining the above two equalities, we obtain that 
\begin{align*}
\overline{U^{(\ell)*}}=
V^{(\ell)} {U^{(\ell)}}
=\sum_{j=1}^{d_\ell} F_{j,j^*}^{(\ell)} U^{(\ell)},
\end{align*}
i.e. $u_{j,k^*}^{(\ell)}=\overline{u_{j,k}^{(\ell)}}$.
Suppose that $\displaystyle U^{(\ell)*}\widehat{\Delta}^{(\ell)} U^{(\ell)}=\sum_{j=1}^{d_\ell}\mu_j^{(\ell)}E_{j,j}^{(\ell)}$ for $\ell=\ell^*$ and denote by $F_\ell =\diag(\mu_1^{(\ell)}, \ldots ,\mu_{d_{\ell}}^{(\ell)})\in M_{d_\ell}(\bC)$ for $\ell= \ell^*$.
Then $\mu_{j^*}^{(\ell)}={\mu_j^{(\ell)}}^{-1}$ for $\ell=\ell^*$ and 
\begin{align*}
\overline{\widehat{\Delta}^{(\ell)}}=\overline{\sum_{j=1}^{d_\ell} \mu_j^{(\ell)} F_{j,j}^{(\ell)}} =\sum_{j=1}^{d_\ell}  \mu_j^{(\ell)} F_{j^*,j^*}^{(\ell)}=\widehat{\Delta}^{(\ell)-1}.
\end{align*}

\begin{notation}
For convenience, we denote by $U_\ell=\left(u_{j,k}^{(\ell)}F_{j,k}^{(\ell)}\right)_{j,k=1}^{d_\ell} \in (\cM'\cap \cM_2)\otimes M_{d_\ell}(\bC)$.
\end{notation}




\begin{remark}
In \cite{VerWir23}, the $\mathscr{V}$-transform of a quantum channel introduced by Vernooij and Wirth is the key to characterize the generator $\cL$. 
In contrast to their method, our approach leverages the finiteness of the inclusion.
\end{remark}

\begin{remark}
Suppose that $\overline{\widehat{\Delta}}=\widehat{\Delta}^{-1}$.
The algebra $\mathcal{A}$ will be the algebra generated by $\widehat{\cL}_{\Delta}$ and $\widehat{\Delta}( \cR(\widehat{\cL}_0)\vee \cR(\overline{\widehat{\cL}_0}))$.
In this case, the matrix $B_\ell$ might not have full rank.
\end{remark}

\section{Divergence and Gradient}
In this section, we shall introduce the directional derivations and divergence, gradients etc.
For each $k, \ell$, we are interesting in directions $F_{k,k}^{(\ell)}$.
The associated directional derivation $\partial_k^{(\ell)}: \cM \to \cM_1$ is given by 
\begin{align*}
    (\partial_k^{(\ell)} x) e_2=  \lambda^{-1/2} F_{k,k}^{(\ell)} [x, e_1]e_2 
 = \left[x , \mathfrak{F}^{-1}(F_{k,k}^{(\ell)})\right]e_2, \quad \text{ for all } x\in \cM.
\end{align*}
Moreover, we see that for any $x\in \cM_1$, the adjoint of the directional derivations is
\begin{equation}\label{eq:adder}
\begin{aligned}
\partial_{j}^{(\ell)*} x=& \bE_{\cM} \left(\left[x, \mathfrak{F}^{-1}(F_{j,j}^{(\ell)})^*\right]\right)=\left\{\begin{array}{ll} 
 \bE_{\cM} \left(\left[x, \mathfrak{F}^{-1}(F_{j,j}^{(\ell^*)})\right]\right) = \bE_{\cM}\left(\partial_{j}^{(\ell^*)} x\right), & \ell^* \neq \ell, \\
 \bE_{\cM} \left(\left[x, \mathfrak{F}^{-1}(F_{j^*,j^*}^{(\ell)})\right]\right) = \bE_{\cM}\left(\partial_{j^*}^{(\ell)} x\right), & \ell^* =\ell,
\end{array} \right.
\end{aligned}
\end{equation}
where the domain of $\partial_k^{(\ell)}$ extends to $\cM_1$.

\begin{proposition}\label{prop:support}
Suppose that $\cR(\widehat{\cL}_0)=\cR(\overline{\widehat{\cL}_0})$.
For any $x\in \cM$, the following statements are equivalent:
\begin{enumerate}[(1)]
\item $\left[x, \mathfrak{F}^{-1}(F_{k,k}^{(\ell)})\right]=0$ for all $k, \ell$.
Equivalently, $F_{k,k}^{(\ell)}[x, e_1]e_2=0$ for all $k, \ell$.
\item $\left[x, \mathfrak{F}^{-1}(E_{k,k}^{(\ell)})\right]=0$ for all $k, \ell$. Equivalently, $E_{k,k}^{(\ell)}[x, e_1]e_2=0$ for all $k, \ell$. 
\item $[x, \mathfrak{F}^{-1}(\widehat{\cL}_\Delta^{1/2})]=0$.
\item $[x, \mathfrak{F}^{-1}(\widehat{\cL}_0^{1/2})]=0$.
\end{enumerate}
\end{proposition}
\begin{proof}
(1)$\Rightarrow$(2):
The assumption implies that $F_{j,k}^{(\ell)}[x, e_1]e_2=0$ for all $j,k=1, \ldots, d_\ell$.
Note that $E_{j,k}^{(\ell)}$ is a linear sum of $F_{j,k}^{(\ell)}$.
We see that $E_{k,k}^{(\ell)}[x, e_1]e_2=0$ for all $k=1, \ldots, d_\ell$.
This proves (2).

(2)$\Rightarrow$(1): 
The assumption implies that $E_{j,k}^{(\ell)}[x, e_1]e_2=0$ for all $j,k=1, \ldots, d_\ell$.
Note that $\Span\{E_{j,k}^{(\ell)}\}_{j,k=1}^{d_\ell}=\Span\{F_{j,k}^{(\ell)}\}_{j,k=1}^{d_\ell}$.
We see that $F_{k,k}^{(\ell)}[x, e_1]e_2=0$ for all $k=1, \ldots, d_\ell$.
This proves (1).

(3)$\Leftrightarrow$(2): It follows from the assumption.

(3)$\Leftrightarrow$(4): By the assumption, we have that $\cR(\widehat{\cL}_{\Delta})=\cR(\widehat{\cL}_0)$.
This shows that the equivalence of $(3)$ and $(4)$.
\end{proof}

\begin{definition}[Relative Ergodicity]
Suppose that $\{\Phi_t\}_{t\geq 0}$ is a bimodule KMS symmetric quantum Markov semigroup.
We say $\{\Phi_t\}_{t\geq 0}$ is relatively ergodic if $\displaystyle \bigvee_{k \geq 1} \cR(\widehat{\cL}_0^{*k})=1$.
\end{definition}

\begin{remark}
The relative ergodicity of bimodule KMS symmetric semigroup is different from the usual ergodicity which requires that the fixed point space $\mathscr{M}(\{\Phi_t\}_{t\geq 0})=\cN$.
\end{remark}

\begin{lemma}\label{lem:kmssupport}
Suppose $\{\Phi_t\}_{t\geq 0}$ is bimodule KMS symmetric with respect to $\widehat{\Delta}$ and relatively ergodic.
Then
\begin{align}
    \Ker\Gamma=\Ker\cL = \nfix(\{\Phi_t\}_{t\geq 0}) =\fix(\{\Phi_t\}_{t\geq 0})=\cN.
\end{align}
\end{lemma}
\begin{proof}
The ergodicity implies $\Ker\Gamma=\cN$.
By the fact that $\Phi_t$ is equilibrium, we have that  $\cN\subset \fix(\{\Phi_t\}_{t\geq 0})\subset \nfix(\{\Phi_t\}_{t\geq 0})\subset \Ker\Gamma=\cN$.
This concludes the proof of the lemma.
\end{proof}

\begin{theorem}\label{thm:kmslimit}
 Suppose $\{\Phi_t\}_{t\geq 0}$ is bimodule KMS symmetric with respect to $\widehat{\Delta}$ and relatively ergodic.
 Then $\displaystyle \lim_{t\to \infty} \widehat{\Phi_t}$ exists.
If $\Phi_t$ is bimodule KMS symmetric with respect to $\widehat{\Delta}$ for all $t\geq 0$, then 
 \begin{align}\label{eq:limit1}
 \left(\lim_{t\to \infty} \widehat{\Phi_t}\right)*\gamma_{1, +}^{-1}\left( \mathfrak{F}^{-1}(\widehat{\Delta})^* \mathfrak{F}^{-1}(\widehat{\Delta})\right)=1.
 \end{align}
\end{theorem}
\begin{proof}
By Lemma \ref{lem:kmssupport}, we have that  $\Ker\Gamma=\Ker\cL= \nfix(\{\Phi_t\}_{t\geq 0})=\cN$.
By Theorem 5.40 in \cite{WuZha25}, we obtain that $\displaystyle \lim_{t\to \infty} \widehat{\Phi_t}=\widehat{\bE}_{\Phi}$.
By the fact that $\{\Phi_t\}_{t\geq 0}$ is relatively ergodic, i.e. $\Ker \cL=\cN$, we obtain that $\displaystyle \lim_{t\to \infty} \bE_{\cN} {\Phi_t}(x)=\lim_{t\to\infty} {\Phi_t}(x)$ for any $x\in \cM$.
We see that $ \widehat{\bE}_{\Phi} * \widehat{\bE}_{\cN}=\widehat{\bE}_{\Phi}$.
This implies that $\widehat{\bE}_{\Phi}= \bE_{\cM_1}(\widehat{\bE}_{\Phi})$, so $\widehat{\mathbb{E}}_{{\Phi}}\in \cM'\cap \cM_1$.
If $\Phi_t$ is bimodule KMS symmetric with respect to $\widehat{\Delta}$, then $\overline{\widehat{\bE}_{\Phi}}= \overline{\widehat{\Delta}} \widehat{\bE}_{\Phi}\overline{\widehat{\Delta}}$. 
Pictorially, we have
\begin{align*}
    \overline{\widehat{\Delta}} \widehat{\mathbb{E}}_{\Phi}\overline{\widehat{\Delta}} = \vcenter{\hbox{\begin{tikzpicture}[scale=1.2]
        \draw [blue] (0, -1)--(0, 1) (0.5, -1)--(0.5, 1);
        \draw [fill=white] (-0.3, -0.2) rectangle (0.3, 0.2);
        \node at (0, 0) {\tiny $\widehat{\mathbb{E}}_{\Phi}$};
        \draw [fill=white] (-0.2, 0.3) rectangle (0.7, 0.7);
        \node at (0.25, 0.5) {\tiny $\overline{\widehat{\Delta}}$};
        \begin{scope}[shift={(0, -1)}]
            \draw [fill=white] (-0.2, 0.3) rectangle (0.7, 0.7);
        \node at (0.25, 0.5) {\tiny $\overline{\widehat{\Delta}}$};
        \end{scope}
    \end{tikzpicture}}} 
    = \vcenter{\hbox{\begin{tikzpicture}[scale=1.2]
        \draw [blue] (0, -0.5)--(0, 1) (0.5, -0.5)--(0.5, 1);
        \begin{scope}[xshift=0.5cm]
            \draw [fill=white] (-0.3, -0.2) rectangle (0.3, 0.2);
        \node at (0, 0) {\tiny $\overline{\widehat{\mathbb{E}}_{\Phi}}$};
        \end{scope}
    \end{tikzpicture}}} = \overline{\widehat{\mathbb{E}}_{\Phi}}. 
\end{align*}
Hence 
\begin{align*}
\widehat{\bE}_{\Phi}*\gamma_{1, +}^{-1}\left( \mathfrak{F}(\overline{\widehat{\Delta}})^* \mathfrak{F}(\overline{\widehat{\Delta}})\right) 
= \lambda^{-1/2}\mathbb{E}_{\cM_1}( \overline{\widehat{\Delta}} \widehat{\mathbb{E}}_{\Phi}\overline{\widehat{\Delta}}) 
= \lambda^{-1/2}\mathbb{E}_{\cM_1}(\overline{\widehat{\mathbb{E}_{\Phi}}}) =1.
\end{align*}
This completes the proof of the theorem.
\end{proof}

\begin{remark}
Suppose that $\cN \subset \cM$ is irreducible and $\{\Phi_t\}_{t\geq 0}$ is a bimodule KMS symmetric quantum Markov semigroup with $\widehat{\Delta}\neq 1$.
Then $\{\Phi_t\}_{t\geq 0}$ is not KMS symmetric with respect to any normal faithful state $\rho$ that preserves the quantum symmetry, i.e. $\sigma_i^{\rho}(e_1)$ is bounded.
In fact, suppose that $\{\Phi_t\}_{t\geq 0}$ is KMS symmetric, we see Equation \eqref{eq:limit1} holds.
Note that $\widehat{\bE}_{\Phi}$ is a multiple of scalar.
We see that $\widehat{\Delta}=1$.
This leads a contradiction.
\end{remark}

We denote by $\displaystyle F_\Phi=\lim_{t\to \infty} \widehat{\Phi_t}$ the limit of the semigroup $\{\Phi_t\}_{t\geq 0}$ when it exists.
In the following, we shall assume that the semigroup is relatively ergodic and give the divergences and gradients.
Let 
\begin{align*}
    \cM_1^{\oplus} =\left\{(x_{k}^{(\ell)})_{k, \ell=1}^{d_\ell, m}: x_{k}^{(\ell)}\in \cM_1,  x_{k}^{(\ell)} e_2=F_{k,k}^{(\ell)} x_{k}^{(\ell)} e_2 \right\} \subset \bigoplus_{\ell=1}^m \cM_1^{\oplus d_\ell},
\end{align*}
with the inner product $\langle \cdot, \cdot\rangle$ given by
\begin{align*}
\left\langle  (x_{k}^{(\ell)})_{k, \ell},  (y_{k}^{(\ell)})_{k, \ell}\right\rangle =\sum_{k, \ell}\tau_1\left( y_{k}^{(\ell)*} x_{k}^{(\ell)}\right)
\end{align*}
whenever $(x_{k}^{(\ell)})_{k, \ell=1}^{d_\ell, m}, (y_{k}^{(\ell)})_{k, \ell=1}^{d_\ell, m}\in \cM_1^{\oplus}$.
The space $\cM_1^{\oplus}$ is a closed space with respect to the norm.
Note that $\cM_1^{\oplus}$ is linear isomorphic to $\cM_1^{\oplus}e_2=\left\{(x_k^{(\ell)})_{k, \ell}: (x_k^{(\ell)})_{k, \ell}\in \cM_1^{\oplus}\right\}$.
Let 
\begin{align*}
\oM=\left\{(x_{k, r}^{(\ell)})_{k, r, \ell} : x_{k,r}^{(\ell)} \in \cM_1, x_{k,r}^{(\ell)}e_2
= F_{k,r}^{(\ell)} y_{k,r}^{(\ell)} e_2, \text{ for some } y_{k,r}^{(\ell)} \in \cM_1  \right\}.
\end{align*}
Let $\oM e_2=\left\{(x_{k, r}^{(\ell)}e_2)_{k, r, \ell} :(x_{k, r}^{(\ell)})_{k, r, \ell} \in \oM \right\}$ with the inner product $\langle \cdot, \cdot \rangle$ defined by 
\begin{align*}
\left\langle (x_{k, r}^{(\ell)}e_2)_{k, r, \ell} , (y_{k, r}^{(\ell)}e_2)_{k, r, \ell} \right\rangle 
= & \sum_{\ell=1}^m \lambda^{-1} \left(\tau_2\otimes \Tr_{d_\ell}\right)\left((y_{k, r}^{(\ell)}e_2)_{k, r}^* (x_{k, r}^{(\ell)}e_2)_{k, r} \right) \\
= & \sum_{k, r, \ell} \lambda^{-1}\tau_2\left(e_2y_{k, r}^{(\ell) *} x_{k, r}^{(\ell)} e_2\right),
\end{align*}
where $(x_{k, r}^{(\ell)}e_2)_{k, r, \ell} , (y_{k, r}^{(\ell)}e_2)_{k, r, \ell} \in \oM e_2.$
We see that $\cM_1^{\oplus}$ is a subspace of $\oM$.
Let $\iota, \widehat{\iota}: \cM_1^{\oplus} \to \oM$ be the inclusions defined as follows:
\begin{align*}
\iota(x_k^{(\ell)})_{k, \ell} =\bigoplus_{\ell=1}^m \begin{pmatrix}
x_{1}^{(\ell)}  &   &   \\
   &\ddots &   \\
 &   & x_{d_{\ell}}^{(\ell)} 
\end{pmatrix}, 
\quad 
\widehat{\iota}(x_k^{(\ell)})_{k, \ell} =\bigoplus_{\ell=1}^m \begin{pmatrix}
x_{1}^{(\ell)}  & \cdots &  x_{1}^{(\ell)}  \\
\vdots  & & \vdots \\
x_{d_\ell}^{(\ell)} & \cdots & x_{d_{\ell}}^{(\ell)}
\end{pmatrix},
\end{align*}
where $(x_k^{(\ell)})_{k, \ell}\in \cM_1^{\oplus}$.
We see that 
\begin{align*}
\iota^* \bigoplus_{\ell=1}^m \begin{pmatrix}
x_{11}^{(\ell)}  & \cdots &  x_{1d_{\ell}}^{(\ell)}  \\
\vdots  & & \vdots \\
x_{d_\ell 1}^{(\ell)} & \cdots & x_{d_{\ell} d_{\ell}}^{(\ell)}
\end{pmatrix}
=(x_{kk}^{(\ell)})_{k,\ell},
\quad 
\widehat{\iota}^* \bigoplus_{\ell=1}^m \begin{pmatrix}
x_{11}^{(\ell)}  & \cdots &  x_{1d_{\ell}}^{(\ell)}  \\
\vdots  & & \vdots \\
x_{d_\ell 1}^{(\ell)} & \cdots & x_{d_{\ell} d_{\ell}}^{(\ell)}
\end{pmatrix}=
\left(\sum_{j=1}^{d_\ell} x_{k j}^{(\ell)}\right)_{k, \ell}.
\end{align*}
We denote by $\iota^{(\ell)}$, $\widehat{\iota}^{(\ell)}$ the $\ell$-components of $\iota$, $\widehat{\iota}$ respectively.

The gradient $\nabla_0: \cM\to \cM_1^{\oplus}$ is defined to be
\begin{align*}
    \nabla_0 x=(\partial_k^{(\ell)} x)_{k, \ell=1}^{d_\ell, m}, \quad x\in \cM.
\end{align*}
This indicates that
\begin{align*}
    ( \nabla_0 x)e_2=\lambda^{-1/2}\left(F_{k,k}^{(\ell)} [x, e_1] e_2\right)_{k, \ell=1}^{d_\ell, m}.
\end{align*}
The kernel of $\nabla_0$ is $\cN$ by the relative ergodicity of the bimodule KMS quantum semigroup.

The divergence $\Div_0: \cM_1^{\oplus}\to \cM$ is defined as follows
\begin{align*}
  \Div_0(x_k^{(\ell)})_{k,\ell=1}^{d_\ell, m} =\sum_{k, \ell=1}^{d_\ell, m} \partial_k^{(\ell)*} x_k^{(\ell)}, \quad (x_k^{(\ell)})_{k, \ell}^{d_\ell, m} \in  \cM_1^{\oplus}.
\end{align*}

\begin{proposition}\label{prop:supportrange}
We have that $\nabla_0^*=\Div_0$.
Moreover, the range $\Ran(\nabla_0)$ is closed and 
\begin{align*}
  \Div_0(x_k^{(\ell)})_{k,\ell=1}^{d_\ell, m}= \lambda^{-1} \sum_{\ell=1}^m\sum_{k=1}^{d_\ell} \bE_{\cM} \left( \left[x_k^{(\ell)}e_2, \fF^{-1}(F_{k,k}^{(\ell)})^*\right]\right)
 \end{align*} 
   for $(x_k^{(\ell)})_{k, \ell} \in  \cM_1^{\oplus}$.
If the semigroup is relatively ergodic, then $\Ker \nabla_0=\cN$ and the closure of the range $\Ran(\Div_0)=\{x\in \cM: \bE_{\cN}(x)=0\}$.
\end{proposition}
\begin{proof}
For any $ (x_k^{(\ell)})_{k, \ell} \in  \cM_1^{\oplus}$ and $y\in \cM$, we have
\begin{align*}
\left\langle (x_k^{(\ell)})_{k, \ell}, \nabla_0y \right\rangle
=& \sum_{\ell=1}^m \sum_{k=1}^{d_\ell}\tau_1\left(x_k^{(\ell)}(\partial_k^{(\ell)} y)^*\right)
=\sum_{\ell=1}^m \sum_{k=1}^{d_\ell}\tau_1\left((\partial_k^{(\ell)*}x_k^{(\ell)})y^*\right)\\
=& \left\langle \Div_0 (x_k^{(\ell)})_{k, \ell}, y \right\rangle.
\end{align*}
This shows that $\nabla_0^*=\Div_0$.
Note that the range of the map $x\mapsto [x, e_1]$ from $\cM$ to $\cM_1$ is closed.
We see that $\Ran \nabla_0$ is closed.

For any $ (x_k^{(\ell)})_{k, \ell} \in  \cM_1^{\oplus}$
\begin{align*}
 \Div_0(x_k^{(\ell)})_{k,\ell=1}^{d_\ell, m}=  \sum_{k, \ell} \bE_{\cM}\left(\left[x_k^{(\ell)}, \fF^{-1}(F_{k,k}^{(\ell)})^*\right]\right)
= \lambda^{-1} \sum_{k, \ell} \bE_{\cM}\left(\left[x_k^{(\ell)}e_2, \fF^{-1}(F_{k,k}^{(\ell)})^*\right]\right). 
\end{align*}

If the semigroup is relatively ergodic, then $\left\{x\in \cM: [x, \mathfrak{F}^{-1}(\widehat{\cL}_0^{1/2})]=0\right\}=\cN$.
By Proposition \ref{prop:support}, we have that $\Ker\nabla_0=\cN$.

Next, we shall proceed to prove that $\Ran (\nabla_0)$ is closed.
Suppose there exists a sequence $\{x_m\}_{m \ge 1}$ in $\cM$ such that $\{\nabla_0 x_m\}_{m \ge 1}$ converges in $\oM$. 
This implies that for any $k$ and $\ell$, the sequence $\{F_{k,k}^{(\ell)}[x_m,e_1] e_2\}_{m \ge 1}$ converges in $\oM e_2$. 
Consequently, 
\begin{align*}
\left\{ \widehat{\cL}_0[x_m,e_1]e_2 = \lambda^{1/2}[x_m, \mathfrak{F}^{-1}(\widehat{\cL}_0)]e_2 \right\}_{m \ge 1}
\end{align*}
converges, which in turn yields that $\left\{ \widehat{\cL}_0^{*n}[x_m,e_1]e_2 \right\}_{m \ge 1}$ converges for all $n \ge 1$. 
Utilizing the assumption that the convolution support of $\widehat{\cL}_0$ is $1$, we deduce that $\lim\limits_{m \rightarrow \infty}[x_m,e_1] = y$ for some $y \in \cM_1$.

The mapping $:x \mapsto [x,e_1]$ for $x\in\cM$ has a closed range.
Without loss of generality, we may assume $\bE_\cN(x_m)=0$ for all $ m \ge 1$. 
Let $z_1=(1-e_1) y e_1$, $z_2 =e_1 y (1-e_1)$, $z_1'=\lambda^{-1} \bE_\cM(z_1)$, and $z_2'=\lambda^{-1} \bE_\cM(z_2)$. 
Setting $\displaystyle z=\frac{z_1'-z_2'}{2}$, we obtain $y=\left[z,e_1\right]$. Thus, $y$ lies in the image of the derivation, which ultimately proves that $\Ran \nabla_0$ is a closed set.
Hence $(\Ker\nabla_0)^{\perp}= \Ran(\Div_0)=\{x\in \cM: \bE_{\cN}(x)=0\}$ followed from the duality.
\end{proof}

\section{Laplacian}

In this section, we undertake a detailed computation of the dual Laplacian associated with a bimodule quantum Markov semigroup. 
Our approach is based on expressing the generators and directional derivatives in terms of their spectral decompositions. 
This perspective not only clarifies the structure of the dual operator but also reveals how the noncommutative geometry of the underlying bimodule influences the form of the Laplacian and its adjoint.

We shall assume that $\cR(\widehat{\cL}_0)=\cR(\overline{\widehat{\cL}_0})$ throughout this section.
The following proposition will express the dual of Laplacian in the form of the divergence.
\begin{proposition}\label{prop:laplacian1}
For any $y\in \cM$, we have that 
\begin{align*}
 \cL_a^*(y) = 
 \frac{\lambda^{-1/2}}{2} \sum_{j,k,r, \ell} \omega_j^{(\ell)} \partial_{k}^{(\ell)*}\left(  \mu_k^{(\ell)-1} \mu_r^{(\ell)-1} y \fF^{-1}\left( E_{j,j}^{(\ell)} F_{r, r}^{(\ell)} \right)   - \mu_k^{(\ell)} \mu_r^{(\ell)} \fF^{-1}\left( E_{j,j}^{(\ell)} F_{r, r}^{(\ell)} \right) y \right).
\end{align*}
\end{proposition}

\begin{proof}
By Lemma 5.23 in \cite{WuZha25}, for any $x\in \cM$, we have that
\begin{align*}
  \cL_a(x)= & -\frac{\lambda^{-5/2}}{2} \bE_{\cM}(e_2e_1 \widehat{\cL}_0 [x, e_1]e_2) + \frac{\lambda^{-5/2}}{2}  \bE_{\cM}( e_2[x, e_1] \widehat{\cL}_0 e_1e_2 )\\
  =&  -\frac{\lambda^{-5/2}}{2} \bE_{\cM}\left(e_2e_1 \widehat{\Delta}^{1/2}\widehat{\cL}_\Delta \widehat{\Delta}^{1/2}[x, e_1]e_2\right) 
  + \frac{\lambda^{-5/2}}{2}  \bE_{\cM}\left( e_2[x, e_1] \widehat{\Delta}^{1/2}\widehat{\cL}_\Delta \widehat{\Delta}^{1/2} e_1e_2\right ).
\end{align*}
By taking the spectral decompositions of $\widehat{\Delta}^{1/2}$ and $\widehat{\cL}_\Delta  $, we obtain that
\begin{align}
 &   2\lambda^{5/2} \tau(y^*\cL_a(x)) \\
 =& - \sum_{j,k,r,  \ell} \omega_j^{(\ell)} \mu_k^{(\ell)} \mu_r^{(\ell)}\tau\left(y^*\bE_{\cM}\left(e_2e_1 F_{r, r}^{(\ell)} E_{j,j}^{(\ell)}F_{k,k}^{(\ell)}[x, e_1]e_2\right)\right)  \label{eq:sum1}\\
 & +\sum_{j,k, r, \ell}\omega_j^{(\ell)} \mu_k^{(\ell)} \mu_{r}^{(\ell)}\tau\left(y^*\bE_{\cM}\left(e_2 [x, e_1]F_{k,k}^{(\ell)}E_{j,j}^{(\ell)} F_{r,r}^{(\ell)} e_1e_2\right)\right) \label{eq:sum2}.
\end{align}

The contragredient on $\cM'\cap \cM_2$ affects different indices of the spectral decompositions. 
For the summand \eqref{eq:sum1}, we split it into two cases: the case $\ell^* \neq \ell$ and the case $\ell^*=\ell$.
Then for the case $\ell^*\neq \ell$, we have that
\begin{align*}
& \sum_{j,k,r,  \ell^*\neq \ell} \omega_j^{(\ell)} \mu_k^{(\ell)} \mu_r^{(\ell)}\tau\left(y^*\bE_{\cM}\left(e_2e_1 F_{r, r}^{(\ell)} E_{j,j}^{(\ell)}F_{k,k}^{(\ell)}[x, e_1]e_2\right)\right)\\
=&\lambda^{1/2} \sum_{j,k,r,  \ell^*\neq \ell} \omega_j^{(\ell)} \mu_k^{(\ell)} \mu_r^{(\ell)}\tau\left(y^*\bE_{\cM}\left(e_2e_1 F_{r, r}^{(\ell)} E_{j,j}^{(\ell)}(\partial_k^{(\ell)} x)e_2\right)\right)\\
=&\lambda  \sum_{j,k,r,  \ell^*\neq \ell} \omega_j^{(\ell)} \mu_k^{(\ell)} \mu_r^{(\ell)}\tau\left(y^*\bE_{\cM}\left(e_2\fF^{-1}\left( E_{j,j}^{(\ell^*)} F_{r, r}^{(\ell^*)} \right) (\partial_k^{(\ell)} x)e_2\right)\right)\\
=&\lambda^2  \sum_{j,k,r,  \ell^*\neq \ell} \omega_j^{(\ell)} \mu_k^{(\ell)} \mu_r^{(\ell)}\tau\left(y^* \fF^{-1}\left( E_{j,j}^{(\ell^*)} F_{r, r}^{(\ell^*)} \right) (\partial_k^{(\ell)} x)\right).
\end{align*}

For the case $\ell^*=\ell$, we have that 
\begin{align*}
& \sum_{j,k,r,  \ell^*= \ell} \omega_j^{(\ell)} \mu_k^{(\ell)} \mu_r^{(\ell)}\tau\left(y^*\bE_{\cM}\left(e_2e_1 F_{r, r}^{(\ell)} E_{j,j}^{(\ell)}F_{k,k}^{(\ell)}[x, e_1]e_2\right)\right)\\
=&\lambda^{1/2} \sum_{j,k,r,  \ell^*= \ell} \omega_j^{(\ell)} \mu_k^{(\ell)} \mu_r^{(\ell)}\tau\left(y^*\bE_{\cM}\left(e_2e_1 F_{r, r}^{(\ell)} E_{j,j}^{(\ell)}(\partial_k^{(\ell)} x)e_2\right)\right)\\
=&\lambda  \sum_{j,k,r,  \ell^*= \ell} \omega_j^{(\ell)} \mu_k^{(\ell)} \mu_r^{(\ell)}\tau\left(y^*\bE_{\cM}\left(e_2\fF^{-1}\left( E_{j,j}^{(\ell)} F_{r^*, r^*}^{(\ell)} \right) (\partial_k^{(\ell)} x)e_2\right)\right)\\
=&\lambda^2  \sum_{j,k,r,  \ell^*= \ell} \omega_j^{(\ell)} \mu_k^{(\ell)} \mu_r^{(\ell)}\tau\left(y^* \fF^{-1}\left( E_{j,j}^{(\ell)} F_{r^*, r^*}^{(\ell)} \right) (\partial_k^{(\ell)} x)\right).
\end{align*}

For the summand \eqref{eq:sum2}, we split it into two cases: the case $\ell^* \neq \ell$ and the case $\ell^*=\ell$.
For the case $\ell^* \neq \ell$, we have that 
\begin{align*}
& \sum_{j,k, r, \ell^* \neq \ell}\omega_j^{(\ell)} \mu_k^{(\ell)} \mu_{r}^{(\ell)}\tau\left(y^*\bE_{\cM}\left(e_2 [x, e_1]F_{k,k}^{(\ell)}E_{j,j}^{(\ell)} F_{r,r}^{(\ell)} e_1e_2\right)\right) \\
=& \lambda^{1/2}\sum_{j,k, r, \ell^* \neq \ell}\omega_j^{(\ell)} \mu_k^{(\ell)} \mu_{r}^{(\ell)}\tau\left(y^*\bE_{\cM}\left(e_2 (\partial_{k}^{(\ell^*)}x)E_{j,j}^{(\ell)} F_{r,r}^{(\ell)} e_1e_2\right)\right) \\
=&\lambda \sum_{j,k, r, \ell^* \neq \ell}\omega_j^{(\ell)} \mu_k^{(\ell)} \mu_{r}^{(\ell)}\tau\left(y^*\bE_{\cM}\left(e_2 (\partial_{k}^{(\ell^*)}x)\fF^{-1}\left(E_{j,j}^{(\ell)} F_{r,r}^{(\ell)} \right)e_2\right)\right) \\
=&\lambda \sum_{j,k, r, \ell^* \neq \ell}\omega_j^{(\ell)} \mu_k^{(\ell)-1} \mu_{r}^{(\ell)-1}\tau\left(y^*\bE_{\cM}\left(e_2 (\partial_{k}^{(\ell)}x)\fF^{-1}\left(E_{j,j}^{(\ell^*)} F_{r,r}^{(\ell^*)} \right)e_2\right)\right) \\
=&\lambda^2 \sum_{j,k, r, \ell^* \neq \ell} \omega_j^{(\ell)} \mu_k^{(\ell)-1} \mu_{r}^{(\ell)-1}\tau_2\left(y^*( \partial_{k}^{(\ell)}x)\fF^{-1}\left(E_{j,j}^{(\ell^*)} F_{r,r}^{(\ell^*)}  \right)\right).
\end{align*}

For the case $\ell^*=\ell$, we have that 
\begin{align*}
& \sum_{j,k, r, \ell^* = \ell}\omega_j^{(\ell)} \mu_k^{(\ell)} \mu_{r}^{(\ell)}\tau\left(y^*\bE_{\cM}\left(e_2 [x, e_1]F_{k,k}^{(\ell)}E_{j,j}^{(\ell)} F_{r,r}^{(\ell)} e_1e_2\right)\right) \\
=& \lambda^{1/2}\sum_{j,k, r, \ell^* = \ell}\omega_j^{(\ell)} \mu_k^{(\ell)} \mu_{r}^{(\ell)}\tau\left(y^*\bE_{\cM}\left(e_2 (\partial_{k^*}^{(\ell)}x)E_{j,j}^{(\ell)} F_{r,r}^{(\ell)} e_1e_2\right)\right) \\
=&\lambda \sum_{j,k, r, \ell^* = \ell}\omega_j^{(\ell)} \mu_k^{(\ell)-1} \mu_{r}^{(\ell)}\tau\left(y^*\bE_{\cM}\left(e_2 (\partial_{k}^{(\ell)}x)\fF^{-1}(E_{j,j}^{(\ell)} F_{r,r}^{(\ell)} )e_2\right)\right) \\
=&\lambda \sum_{j,k, r, \ell^* = \ell}\omega_j^{(\ell)} \mu_k^{(\ell)-1} \mu_{r}^{(\ell)-1}\tau\left(y^*\bE_{\cM}\left(e_2 (\partial_{k}^{(\ell)}x)\fF^{-1}\left(E_{j,j}^{(\ell)} F_{r^*,r^*}^{(\ell)} \right)e_2\right)\right) \\
=&\lambda^2 \sum_{j,k, r, \ell^* = \ell} \omega_j^{(\ell)} \mu_k^{(\ell)-1} \mu_{r}^{(\ell)-1}\tau_2\left(y^*( \partial_{k}^{(\ell)}x)\fF^{-1}\left(E_{j,j}^{(\ell)} F_{r^*,r^*}^{(\ell)}  \right)\right).
\end{align*}

We have that 
\begin{align*}
 &   2\lambda^{5/2} \tau(y^*\cL_a(x))\\
 =& - \lambda^2  \sum_{j,k,r,  \ell^*\neq \ell} \omega_j^{(\ell)} \mu_k^{(\ell)} \mu_r^{(\ell)}\tau\left(y^* \fF^{-1}\left( E_{j,j}^{(\ell^*)} F_{r, r}^{(\ell^*)} \right) (\partial_k^{(\ell)} x)\right)\\
 & + \lambda^2 \sum_{j,k, r, \ell^* \neq \ell} \omega_j^{(\ell)} \mu_k^{(\ell)-1} \mu_{r}^{(\ell)-1}\tau_2\left(y^*( \partial_{k}^{(\ell)}x)\fF^{-1}\left(E_{j,j}^{(\ell^*)} F_{r,r}^{(\ell^*)}  \right)\right) \\
 & - \lambda^2  \sum_{j,k,r,  \ell^*= \ell} \omega_j^{(\ell)} \mu_k^{(\ell)} \mu_r^{(\ell)}\tau\left(y^* \fF^{-1}\left( E_{j,j}^{(\ell)} F_{r^*, r^*}^{(\ell)} \right) (\partial_k^{(\ell)} x)\right)\\
 & + \lambda^2 \sum_{j,k, r, \ell^* \neq \ell} \omega_j^{(\ell)} \mu_k^{(\ell)-1} \mu_{r}^{(\ell)-1}\tau_2\left(y^*( \partial_{k}^{(\ell)}x)\fF^{-1}\left(E_{j,j}^{(\ell)} F_{r^*,r^*}^{(\ell)}  \right)\right)\\
 =&\lambda^2  \sum_{j,k,r,  \ell^*\neq \ell} \omega_j^{(\ell)} \tau_1\left(  \mu_k^{(\ell)-1} \mu_r^{(\ell)-1} \fF^{-1}\left( E_{j,j}^{(\ell^*)} F_{r, r}^{(\ell^*)} \right) y^* (\partial_k^{(\ell)} x)  \right. \\
& \left. - \mu_k^{(\ell)} \mu_r^{(\ell)}y^* \fF^{-1}\left( E_{j,j}^{(\ell^*)} F_{r, r}^{(\ell^*)} \right) (\partial_k^{(\ell)} x)\right)\\
&  +\lambda^2  \sum_{j,k,r,  \ell^*= \ell} \omega_j^{(\ell)} \tau_1\left(  \mu_k^{(\ell)-1} \mu_r^{(\ell)-1} \fF^{-1}\left( E_{j,j}^{(\ell)} F_{r^*, r^*}^{(\ell)} \right) y^* (\partial_k^{(\ell)} x)\right.\\
&\left.  - \mu_k^{(\ell)} \mu_r^{(\ell)}y^* \fF^{-1}\left( E_{j,j}^{(\ell)} F_{r^*, r^*}^{(\ell)} \right) (\partial_k^{(\ell)} x)\right).
\end{align*}
Hence
\begin{align*}
& 2\lambda^{1/2}\cL_a^*(y) \\
=& \sum_{j,k,r,  \ell^*\neq \ell} \omega_j^{(\ell)} \partial_{k}^{(\ell)*}\left(  \mu_k^{(\ell)-1} \mu_r^{(\ell)-1} y \fF^{-1}\left( E_{j,j}^{(\ell^*)} F_{r, r}^{(\ell^*)} \right)^*   - \mu_k^{(\ell)} \mu_r^{(\ell)} \fF^{-1}\left( E_{j,j}^{(\ell^*)} F_{r, r}^{(\ell^*)} \right)^* y \right) \\
& + \sum_{j,k,r,  \ell^*= \ell} \omega_j^{(\ell)} \partial_{k}^{(\ell)*}\left(  \mu_k^{(\ell)-1} \mu_r^{(\ell)-1} y \fF^{-1}\left( E_{j,j}^{(\ell)} F_{r^*, r^*}^{(\ell)} \right)^*   - \mu_k^{(\ell)} \mu_r^{(\ell)} \fF^{-1}\left( E_{j,j}^{(\ell)} F_{r^*, r^*}^{(\ell)} \right)^* y \right) \\
=& \sum_{j,k,r,  \ell^*\neq \ell} \omega_j^{(\ell)} \partial_{k}^{(\ell)*}\left(  \mu_k^{(\ell)-1} \mu_r^{(\ell)-1} y \fF^{-1}\left( E_{j,j}^{(\ell)} F_{r, r}^{(\ell)} \right)   - \mu_k^{(\ell)} \mu_r^{(\ell)} \fF^{-1}\left( E_{j,j}^{(\ell)} F_{r, r}^{(\ell)} \right) y \right) \\
& + \sum_{j,k,r,  \ell^*= \ell} \omega_j^{(\ell)} \partial_{k}^{(\ell)*}\left(  \mu_k^{(\ell)-1} \mu_r^{(\ell)-1} y \fF^{-1}\left( E_{j,j}^{(\ell)} F_{r, r}^{(\ell)} \right)   - \mu_k^{(\ell)} \mu_r^{(\ell)} \fF^{-1}\left( E_{j,j}^{(\ell)} F_{r, r}^{(\ell)} \right) y \right)\\
=& \sum_{j,k,r, \ell} \omega_j^{(\ell)} \partial_{k}^{(\ell)*}\left(  \mu_k^{(\ell)-1} \mu_r^{(\ell)-1} y \fF^{-1}\left( E_{j,j}^{(\ell)} F_{r, r}^{(\ell)} \right)   - \mu_k^{(\ell)} \mu_r^{(\ell)} \fF^{-1}\left( E_{j,j}^{(\ell)} F_{r, r}^{(\ell)} \right) y \right).
\end{align*}
This completes the computation of the proposition.
\end{proof}

\begin{remark}
If the semigroup is (bimodule) GNS symmetric, then $E_{j,j}^{(\ell)}=F_{j,j}^{(\ell)}$ and $j$ is forced to be $1$.
The Laplacian is 
\begin{align*}
2\lambda^{1/2}\cL_a^*(y)= \sum_{ \ell=1}^m \omega_k^{(\ell)} \partial_{k}^{(\ell)*}\left(  \mu_k^{(\ell)-2}  y \fF^{-1}\left(  F_{k, k}^{(\ell)} \right)   - \mu_k^{2(\ell)} \fF^{-1}\left( F_{k, k}^{(\ell)} \right) y \right),
\end{align*}
where $k=1$.
\end{remark}

\begin{proposition}
For any $y\in \cM$, we have that 
\begin{align*}
\cL_a^*(y) 
=\frac{\lambda^{-1/2}}{2}\sum_{j,k,r, \ell} \omega_j^{(\ell)} \partial_{k}^{(\ell)*}\left(  \mu_k^{(\ell)-1} \mu_r^{(\ell)-1} y \fF^{-1}\left( F_{k,k}^{(\ell)}E_{j,j}^{(\ell)} F_{r, r}^{(\ell)} \right)   - \mu_k^{(\ell)} \mu_r^{(\ell)} \fF^{-1}\left( F_{k,k}^{(\ell)} E_{j,j}^{(\ell)} F_{r, r}^{(\ell)} \right) y \right).
\end{align*}
\end{proposition}
\begin{proof}
For any $x\in \cM$, we see that 
\begin{align*}
& 2\lambda^{1/2} \tau(\cL_a^*(y) x^*)=\langle 2\lambda^{1/2}\cL_a^*(y), x \rangle \\
=&  \sum_{j,k,r, \ell} \omega_j^{(\ell)} \left\langle \partial_{k}^{(\ell)*}\left(  \mu_k^{(\ell)-1} \mu_r^{(\ell)-1} y \fF^{-1}\left( E_{j,j}^{(\ell)} F_{r, r}^{(\ell)} \right)   - \mu_k^{(\ell)} \mu_r^{(\ell)} \fF^{-1}\left( E_{j,j}^{(\ell)} F_{r, r}^{(\ell)} \right) y \right), x\right \rangle \\
=&  \sum_{j,k,r, \ell} \omega_j^{(\ell)} \left\langle \left(  \mu_k^{(\ell)-1} \mu_r^{(\ell)-1} y \fF^{-1}\left( E_{j,j}^{(\ell)} F_{r, r}^{(\ell)} \right)   - \mu_k^{(\ell)} \mu_r^{(\ell)} \fF^{-1}\left( E_{j,j}^{(\ell)} F_{r, r}^{(\ell)} \right) y \right), \partial_{k}^{(\ell)} x\right \rangle \\
=& \lambda^{-1} \sum_{j,k,r, \ell} \omega_j^{(\ell)} \left\langle \left(  \mu_k^{(\ell)-1} \mu_r^{(\ell)-1} y \fF^{-1}\left( E_{j,j}^{(\ell)} F_{r, r}^{(\ell)} \right)   - \mu_k^{(\ell)} \mu_r^{(\ell)} \fF^{-1}\left( E_{j,j}^{(\ell)} F_{r, r}^{(\ell)} \right) y \right)e_2, ( \partial_{k}^{(\ell)} x) e_2\right \rangle \\
=& \lambda^{-1} \sum_{j,k,r, \ell} \omega_j^{(\ell)} \left\langle \left(  \mu_k^{(\ell)-1} \mu_r^{(\ell)-1} y \fF^{-1}\left( F_{k,k}^{(\ell)}E_{j,j}^{(\ell)} F_{r, r}^{(\ell)} \right)   - \mu_k^{(\ell)} \mu_r^{(\ell)} \fF^{-1}\left( F_{k,k}^{(\ell)}E_{j,j}^{(\ell)} F_{r, r}^{(\ell)} \right) y \right)e_2, ( \partial_{k}^{(\ell)} x) e_2\right \rangle \\
=&  \sum_{j,k,r, \ell} \omega_j^{(\ell)} \left\langle \left(  \mu_k^{(\ell)-1} \mu_r^{(\ell)-1} y \fF^{-1}\left( F_{k,k}^{(\ell)}E_{j,j}^{(\ell)} F_{r, r}^{(\ell)} \right)   - \mu_k^{(\ell)} \mu_r^{(\ell)} \fF^{-1}\left( F_{k,k}^{(\ell)}E_{j,j}^{(\ell)} F_{r, r}^{(\ell)} \right) y \right), \partial_{k}^{(\ell)} x \right \rangle.
\end{align*}
This proves the proposition.
\end{proof}

Now we would like to further write the formulation in Proposition \ref{prop:laplacian1} in terms of gradients.
Let
\begin{align*}
  \widehat{\Delta}_{k,r}^{(\ell)} =  \sum_{j=1}^{d_\ell}  \omega_j^{(\ell)} F_{k,k}^{(\ell)}E_{j,j}^{(\ell)} F_{r,r}^{(\ell)}.
\end{align*}
Then 
\begin{align}\label{eq:lap1}
 \cL_a^*(y) 
 = \frac{\lambda^{-1/2}}{2} \sum_{k, r, \ell}  \partial_k^{(\ell)*} \left(\mu_k^{(\ell)-1} \mu_r^{(\ell)-1}  y \mathfrak{F}^{-1}(\widehat{\Delta}_{k,r}^{(\ell)}) - \mu_k^{(\ell)}\mu_r^{(\ell)}  \mathfrak{F}^{-1}(\widehat{\Delta}_{k,r}^{(\ell)}) y\right).
\end{align}
Note that  $F_{1,1}^{(\ell)}, \ldots, F_{d_\ell, d_\ell}^{(\ell)}$ are the directions for the gradients.
Let
\begin{align*}
    \Pi^{(\ell)} =   U_\ell^* B_\ell U_\ell \in ( \cM'\cap \cM_2)\otimes M_{d_\ell}(\bC).
\end{align*}
Note that $\displaystyle U^{(\ell)}=\sum_{j,k=1}^{d_\ell} u_{j,k}^{(\ell)} F_{j,k}^{(\ell)}$.
We have that 
\begin{align*}
\widehat{\Delta}_{k,r}^{(\ell)} =&   \sum_{j=1}^{d_\ell}  \omega_j^{(\ell)} F_{k,k}^{(\ell)}E_{j,j}^{(\ell)} F_{r,r}^{(\ell)} \\
=&  \sum_{j=1}^{d_\ell}  \omega_j^{(\ell)} F_{k,k}^{(\ell)} U^{(\ell)*}F_{j,j}^{(\ell)} U^{(\ell)}F_{r,r}^{(\ell)} \\
=&  \sum_{j=1}^{d_\ell} \overline{u_{j,k}^{(\ell)}} \omega_j^{(\ell)} u_{j,r}^{(\ell)} F_{k,r}^{(\ell)} . 
\end{align*}
This shows that 
\begin{align*}
\Pi^{(\ell)} = \left(\widehat{\Delta}_{k,r}^{(\ell)} \right)_{k, r=1}^{d_\ell} =U_\ell^* B_\ell U_\ell.
\end{align*}

The matrix $\Pi^{(\ell)}$, $\ell=1, \ldots, m$ induce a linear map from $ \oM $ to $ \oM$ defined by 
\begin{align*}
\Pi(x_{k, r}^{(\ell)})_{k,r, \ell} 
=& \bigoplus_{\ell} \lambda^{-1} (\bE_{\cM_1}\otimes I) \left(\Pi^{(\ell)} (x_{k,r}^{(\ell)} e_2)_{k, r=1}^{d_\ell}\right).
\end{align*}
where $\bE_{\cM_1}^{\oplus d_\ell}$ is the direct sum of $d_\ell$ copies of the conditional expectations $\bE_{\cM_1}$.
It is worth to note that
\begin{align*}
  \left( U_\ell^* B_\ell U_\ell\right) \left( U_\ell^* B_\ell^{-1} U_\ell \right)   =\diag(F_{1,1}^{(\ell)}, \ldots , F_{d_\ell, d_\ell}^{(\ell)}).
\end{align*}

\begin{proposition}\label{prop:invertible}
The map $\Pi: \oM \to \oM$ is invertible.
\end{proposition}
\begin{proof}
We define the map $\Pi^{-1}: \oM  \to \oM$ as follows:
\begin{align*}
\Pi^{-1}(x_{k,r}^{(\ell)})_{k, r, \ell} =
 \bigoplus_{\ell=1}^m  \lambda^{-1} (\bE_{\cM_1}\otimes I) \left( \left( U_\ell^* B_\ell^{-1} U_\ell\right)   
(x_{k,r}^{(\ell)} e_2)_{k, r=1}^{d_\ell}\right).
\end{align*}
It suffices to show that $\Pi\, \Pi^{-1}$ and  $\Pi^{-1} \Pi  $  are the identity map on $\oM$.
Note that 
\begin{align*}
\Pi\,\left( \Pi^{-1}(x_{k,r}^{(\ell)})_{k, r, \ell} \right)e_2 
= & \bigoplus_{\ell}  \left( U_\ell^* B_\ell U_\ell\right)   \Pi^{-1}(x_{k,r}^{(\ell)})_{k, r=1}^{d_\ell} e_2\\
=& (x_{k,r}^{(\ell)})_{k,r,  \ell} e_2
\end{align*}
by the fact $\lambda^{-1} \bE_{\cM_1}(x e_2)e_2=xe_2,\, x \in \cM_2$. This shows that $\Pi\, \Pi^{-1}$ is the identity map.
Similarly, we have that $ \Pi^{-1}\Pi$ is the identity map.
Hence $\Pi^{-1}$ is the inverse of $\Pi$.
\end{proof}


Now we see that the matrix $\displaystyle \bigoplus_{\ell=1}^m \Pi^{(\ell)}$ plays a key role in the study of the Laplacian of the semigroups, denoted by $\Pi$ again if there is no confusion.
We shall call it as the directional matrix of the semigroup $\{\Phi_t\}_{t\geq 0}$.

Now we recall the weight transform $\K_D$ introduced by Carlen and Maas \cite{CarMaa17}.
Suppose that $D\in \cM$ is strictly positive and $\mu>0$. 
We define the linear map $\K_{D, \mu}:\cM_1\to \cM_1$ as follows: 
\begin{align*}
\K_{D, \mu} (x)=\int_0^1 \mu^{1-2s}D^{s} x D^{1-s} ds, \quad x\in \cM_1.   
\end{align*}
Then for any $v\in \cM_1$,
\begin{align}\label{eq:log}
 \K_{D, \mu}\left((\log\mu^{-1} D ) v- v\log \mu D \right) 
 = \mu^{-1} D v  - \mu v D.
\end{align}
The inverse of $\K_{D, \mu}$ is known to be 
\begin{align*}
  \K_{D, \mu}^{-1}(x) = \int_0^\infty (s+\mu D)^{-1} x  (s+\mu^{-1} D)^{-1} ds, \quad x\in \cM_1.
\end{align*}
For any $k, r, \ell$, we let $\K_{D, k, r, \ell}=\K_{D, \mu_k^{(\ell)}\mu_r^{(\ell)}}$ and $\K_{D}: \oM \to \oM$ defined by
\begin{align*}
& \K_{D}(x_{k,r}^{(\ell)})_{k,r, \ell} =\bigoplus_{\ell=1}^m \left( \K_{D, k, r, \ell} (x_{k,r}^{(\ell)})\right)_{k,r=1}^{d_{\ell}}\\
=&  \int_0^1
\bigoplus_{\ell=1}^m  \begin{pmatrix} \mu_1^{(\ell)(1-2s)} D^s & &  \\  & \ddots &  \\ & & \mu_{d_\ell}^{(\ell)(1-2s)}D^s \end{pmatrix}
(x_{k,r}^{(\ell)})_{k,r}  \begin{pmatrix}
  D^{1-s} \mu_1^{(\ell)(1-2s)}& &  \\ & \ddots &  \\ & &   D^{1-s}\mu_{d_\ell}^{(\ell)(1-2s)}
\end{pmatrix} ds\\
=&  \int_0^1
\bigoplus_{\ell=1}^m F_\ell^{1-2s} \left(D^s x_{k,r}^{(\ell)} D^{1-s}\right)_{k,r} F_\ell^{1-2s} ds
\end{align*}
for all $(x_{k,r}^{(\ell)})_{k,r, \ell}\in \oM$.
We denote by $\K_D^{(\ell)}$ the $\ell$-component of $\K_D$.

%
\begin{proposition}\label{prop:kd}
The linear map $\K_{D}: \oM \to \oM$ is invertible.
Moreover, for any $x\in \cM_1$, 
\begin{align*}
 \left(\K_{D, k, r, \ell}(x)\right)^*=\left\{\begin{array}{ll}
\K_{D, k, r, \ell^*}(x^*), & \ell^* \neq \ell \\
\K_{D, k^*, r^*, \ell}(x^*). & \ell^*=\ell.
\end{array}\right.
\end{align*}
\end{proposition}
\begin{proof}
We define the linear map $\K_{D}^{-1}: \oM \to \oM$ as follows:
\begin{align*}
\K_{D}^{-1} (x_{k,r}^{(\ell)})_{k, r,\ell} =\bigoplus_{\ell=1}^m \left( \K_{D, k, r, \ell}^{-1} (x_{k,r}^{(\ell)})\right)_{k,r=1}^{d_{\ell}}.
\end{align*}
for all $(x_{k,r}^{(\ell)})_{k,r, \ell}\in \oM$.
Then $\K_D\K_D^{-1}$ and $\K_D^{-1}\K_D$ are identity map on $\oM$.
Moreover, for any $x\in \cM_1$, we have that 
\begin{align*}
 \left(\K_{D, k, r, \ell}(x)\right)^*
 =&  \int_0^1
  (\mu_k^{(\ell)}\mu_r^{(\ell)})^{1-2s}  D^{1-s} x^* D^{s} ds \\
=&  \int_0^1 (\mu_k^{(\ell)}\mu_r^{(\ell)})^{-1+2s}  D^{s} x^* D^{1-s}  ds \\
=&  \int_0^1 (\mu_k^{(\ell)-1}\mu_r^{(\ell)-1})^{1-2s}  D^{s} x^* D^{1-s}  ds \\
=&\left\{\begin{array}{ll}
\K_{D, k, r, \ell^*}(x^*), & \ell^* \neq \ell, \\
\K_{D, k^*, r^*, \ell}(x^*), & \ell^*=\ell.
\end{array}\right.
\end{align*}
This completes the proof of the proposition.
\end{proof}

\begin{theorem}\label{thm:laplacian2}
Suppose that $D\in \cM$ is strictly positive.
Then 
\begin{align*}
\cL_a^*(D)=\frac{\lambda^{-1/2}}{2}\Div_0 \widehat{\iota}^* \K_D \Pi \iota (\nabla_0 \log D)-
 \lambda^{-1/2}\Div_0 \widehat{\iota}^* \K_D \bigoplus_{\ell=1}^m\left ((\log \mu_k^{(\ell)}\mu_r^{(\ell)})  \mathfrak{F}^{-1}(\widehat{\Delta}_{k,r}^{(\ell)})\right)_{k, r=1}^{d_\ell}.
\end{align*}
\end{theorem}
\begin{proof}
By Equations \eqref{eq:lap1} and \eqref{eq:log}, we have that
\begin{equation}\label{eq:lap23}
\begin{aligned}
 \cL_a^*(D) =  &\frac{\lambda^{-1/2}}{2} \sum_{k, r, \ell}  \partial_k^{(\ell)*} \left(\mu_k^{(\ell)-1} \mu_r^{(\ell)-1}  D \mathfrak{F}^{-1}(\widehat{\Delta}_{k,r}^{(\ell)}) - \mu_k^{(\ell)}\mu_r^{(\ell)}  \mathfrak{F}^{-1}(\widehat{\Delta}_{k,r}^{(\ell)}) D\right)\\   
 =& \frac{\lambda^{-1/2}}{2} \sum_{k, r, \ell}    \partial_k^{(\ell)*} \K_{D, k, r, \ell} \left( (\log \mu_k^{(\ell)-1}\mu_r^{(\ell)-1}  D) \mathfrak{F}^{-1}(\widehat{\Delta}_{k,r}^{(\ell)}) -  \mathfrak{F}^{-1}(\widehat{\Delta}_{k,r}^{(\ell)}) (\log \mu_k^{(\ell)} \mu_r^{(\ell)}D)\right)\\
  =& \frac{\lambda^{-1/2}}{2} \sum_{k,r, \ell}  \partial_k^{(\ell)*} \K_{D, k, r, \ell} \left( (\log D) \mathfrak{F}^{-1}(\widehat{\Delta}_{k,r}^{(\ell)}) -  \mathfrak{F}^{-1}(\widehat{\Delta}_{k,r}^{(\ell)}) (\log  D)\right)\\
 & -  \lambda^{-1/2} \sum_{k, r, \ell}  (\log \mu_k^{(\ell)} \mu_r^{(\ell)})   \partial_k^{(\ell)*} \K_{D, k, r, \ell} \left( \mathfrak{F}^{-1}(\widehat{\Delta}_{k,r}^{(\ell)})\right).
\end{aligned}
\end{equation}
 
For the first sum in Equation \eqref{eq:lap23}, we have that 
\begin{align*}
& \sum_{k,r,  \ell} \partial_k^{(\ell)*} \K_{D, k,r, \ell} ( (\log D) \mathfrak{F}^{-1}(\widehat{\Delta}_{k,r}^{(\ell)}) -  \mathfrak{F}^{-1}(\widehat{\Delta}_{k,r}^{(\ell)}) (\log D))\\
 =& \lambda^{-3/2}  \sum_{k, r, \ell} \partial_k^{(\ell)*} \K_{D, k, r, \ell} \bE_{\cM_1}\left( \widehat{\Delta}_{k,r}^{(\ell)} [\log D, e_1] e_2\right)\\
 =& \lambda^{-3/2} \Div_0 \widehat{\iota}^* \bigoplus_{\ell=1}^m \big( ( \K_{D, k,r,  \ell} \bE_{\cM_1}(\widehat{\Delta}_{k,r}^{(\ell)} [\log D, e_1] e_2))\big)_{k, r=1}^{d_{\ell}}\\
  =& \lambda^{-3/2} \Div_0 \widehat{\iota}^* \K_D\bigoplus_{\ell=1}^m \big( \bE_{\cM_1}(\widehat{\Delta}_{k,r}^{(\ell)} [\log D, e_1] e_2))\big)_{k, r=1}^{d_{\ell}}\\
 =& \lambda^{-3/2} \Div_0 \widehat{\iota}^* \K_D\bigoplus_{\ell=1}^m  (\bE_{\cM_1}\otimes I )\Pi^{(\ell)}\iota^{(\ell)}\left(F_{k,k}^{(\ell)} [\log D, e_1] e_2\right)_{k=1}^{d_{\ell}}\\
 =& \lambda^{-1} \Div_0 \widehat{\iota}^* \K_D\bigoplus_{\ell=1}^m  (\bE_{\cM_1}\otimes I )\Pi^{(\ell)}\iota^{(\ell)}\left( (\partial_k^{(\ell)}\log D ) e_2\right)_{k=1}^{d_{\ell}}\\
  =& \Div_0 \widehat{\iota}^* \K_D \Pi \iota (\nabla_0 \log D).
\end{align*}   
For the second sum in Equation \eqref{eq:lap23}, we have that 
\begin{align*}
& \sum_{k, r, \ell}  (\log \mu_k^{(\ell)} \mu_r^{(\ell)})  \partial_k^{(\ell)*} \K_{D, k, r, \ell} ( \mathfrak{F}^{-1}(\widehat{\Delta}_{k,r}^{(\ell)}))\\
=& \Div_0 \widehat{\iota}^* \K_D \bigoplus_{\ell=1}^m\left ((\log \mu_k^{(\ell)}\mu_r^{(\ell)})  \mathfrak{F}^{-1}(\widehat{\Delta}_{k,r}^{(\ell)})\right)_{k, r=1}^{d_\ell}.
\end{align*}
This completes the computation.
\end{proof}

\begin{remark}
Suppose that $\overline{\widehat{\Delta}}=\widehat{\Delta}^{-1}$.
Theorem \ref{thm:laplacian2} is still true.
Proposition \ref{prop:invertible} might not hold, but it is still true by adjusting the divergence and gradients.

\end{remark}

Now we shall switch the $\iota$ to $\widehat{\iota}$.
To do so, we shall use some notations.
Let $\mathbf{D}^{(\ell)} =\diag(\underbrace{D, \ldots, D}_{d_\ell} )$ and define a linear map $\K_{D, \Pi}: \oM \to \oM$ as follows:
\begin{align*}
\K_{D, \Pi}(x_{k,r}^{(\ell)})_{k,r, \ell}=\bigoplus_{\ell=1}^m \frac{1}{d_\ell} \int_0^1 F_\ell^{1-2s} \mathbf{D}^{(\ell) s} \Pi^{(\ell)} \left(F_\ell^{1-2s}  (x_{k,r}^{(\ell)})_{k,r=1}^{d_\ell}\right) \mathbf{D}^{(\ell) (1-s)}  ds,
\end{align*}
where $(x_{k,r}^{(\ell)})_{k,r, \ell}\in \oM$, $\Pi^{(\ell)}$ is viewed as a linear map from $\oM$ to $\oM$.
We denote by $\K_{D, \Pi}^{(\ell)}$ the $\ell$-component of $\K_{D, \Pi}$.

\begin{proposition}
The linear map $\K_{D, \Pi}$ is invertible and positive acting on $\oM$.
\end{proposition}
\begin{proof}
Suppose $\displaystyle (x_{k,r}^{(\ell)})_{k,r, \ell} \in \oM$ and $\mathbf{X}^{(\ell)} =(x_{k,r}^{(\ell)} e_2)_{k,r=1}^{d_\ell} $ we have that 
\begin{align*}
& \left\langle \K_{D, \Pi}  (x_{k,r}^{(\ell)})_{k,r, \ell} ,  (x_{k,r}^{(\ell)})_{k,r, \ell}  \right\rangle \\
=& \sum_{\ell=1}^m \lambda^{-1}\left\langle \K_{D, \Pi}^{(\ell)}  \mathbf{X}^{(\ell)} ,  \mathbf{X}^{(\ell)}  \right\rangle \\
=& \lambda^{-1} \sum_{\ell=1}^m \frac{1}{d_\ell}\int_0^1 (\tau_2\otimes \Tr_{d_\ell} )\left( F_\ell^{1-2s} \mathbf{D}^{(\ell) s} \Pi^{(\ell)} F_\ell^{1-2s}   \mathbf{X}^{(\ell)} \mathbf{D}^{(\ell) (1-s)}  \mathbf{X}^{(\ell)*} \right)ds.
\end{align*}
Note that $\mathbf{D}^{(\ell) s} \Pi^{(\ell)}= \Pi^{(\ell)} \mathbf{D}^{(\ell) s}$.
We see that $F_\ell^{1-2s} \mathbf{D}^{(\ell) s} \Pi^{(\ell)} F_\ell^{1-2s}$ is positive semidefinite and
\begin{align*}
F_\ell^{1-2s} \mathbf{D}^{(\ell) s} \Pi^{(\ell)} F_\ell^{1-2s}
\geq  \left(\min_{1\leq j \leq d_{\ell}}\{\omega_j^{(\ell)}\}\right) \left(\min_{1\leq j \leq d_{\ell}}\{\mu_j^{(\ell)(2-4s)}\}\right)\|D^{-1}\|^{-s} \diag(F_{11}^{(\ell)}, \ldots,F_{d_{\ell}d_\ell}^{(\ell)}).
\end{align*}
By the fact that $\mathbf{D}^{(\ell) (1-s)} \geq   \|D^{-1}\|^{-1+s} \diag(F_{11}^{(\ell)}, \ldots,F_{d_{\ell}d_\ell}^{(\ell)})$,
we obtain that 
\begin{align*}
& \left\langle \K_{D, \Pi}  (x_{k,r}^{(\ell)})_{k,r, \ell} ,  (x_{k,r}^{(\ell)})_{k,r, \ell}  \right\rangle \\
\geq &   \min_{1\leq \ell\leq m}  \left( \frac{1}{d_\ell} \left(\min_{1\leq j \leq d_{\ell}}\{\omega_j^{(\ell)}\}\right)  \left(\min_{1\leq j \leq d_{\ell}}\{\mu_j^{(\ell) \pm 2}\}\right)\right) \|D^{-1}\|^{-1}  \left\langle  (x_{k,r}^{(\ell)})_{k,r, \ell} ,  (x_{k,r}^{(\ell)})_{k,r, \ell}  \right\rangle.
\end{align*}
This shows that $\K_{D, \Pi}$ is invertible and positive acting on $\oM$.
\end{proof}

Note that $\K_{D, \Pi}\neq \K_D \Pi$.
Now we shall write the Laplacian in terms of $\K_{D, \Pi}$ as follows.

\begin{proposition}\label{prop:reform}
We have that 
\begin{align*}
\widehat{\iota}^* \K_D \Pi\  \iota= \widehat{\iota}^* \K_{D, \Pi}\  \widehat{\iota}.
\end{align*}
\end{proposition}

\begin{proof}
Recall that $\mathbbm{1}^{(\ell)} =  \begin{pmatrix} 1 & \cdots &1 \\ \vdots & \ddots & \vdots \\ 1&\cdots & 1
\end{pmatrix} \in M_{d_\ell} (\bC)$.
Suppose that $(x_k^{(\ell)})_{k, \ell} , (y_k^{(\ell)})_{k, \ell} \in \cM_1^{\oplus}$.
Let $\mathbf{X}^{(\ell)} =  \diag(x_1^{(\ell)}e_2, \ldots, x_{d_\ell}^{(\ell)}e_2)$ and $\mathbf{Y}^{(\ell)} =  \diag(y_1^{(\ell)}e_2, \ldots, y_{d_\ell}^{(\ell)}e_2)$.

Then we have that 
\begin{align*}
& \left\langle \widehat{\iota}^* \K_D \Pi \iota(x_k^{(\ell)})_{k, \ell}, (y_k^{(\ell)})_{k, \ell} \right\rangle \\
=& \lambda^{-1}  \left\langle \widehat{\iota}^* \K_D \Pi \iota(x_k^{(\ell)})_{k, \ell} e_2, (y_k^{(\ell)})_{k, \ell} e_2\right\rangle \\
=& \lambda^{-1} 
\sum_{\ell=1}^m \left\langle \K_D \Pi^{(\ell)}\mathbf{X}^{(\ell)} ,
\mathbf{Y}^{(\ell)} \mathbbm{1}^{(\ell)} \right\rangle\\
=& \lambda^{-1} \int_0^1 \sum_{\ell=1}^m \left \langle F_\ell^{1-2s} \mathbf{D}^{(\ell) s} \Pi^{(\ell)} \mathbf{X}^{(\ell)}\mathbf{D}^{(\ell) (1-s)} F_\ell^{1-2s}, \mathbf{Y}^{(\ell)} \mathbbm{1}^{(\ell)} \right\rangle  ds\\
=&  \lambda^{-1} \sum_{\ell=1}^m \int_0^1  (\tau_2\otimes \Tr_{d_\ell}) \left( F_\ell^{1-2s} \mathbf{D}^{(\ell) s} \Pi^{(\ell)} \mathbf{X}^{(\ell)}\mathbf{D}^{(\ell) (1-s)} F_\ell^{1-2s} \mathbbm{1}^{(\ell)}  \mathbf{Y}^{(\ell)*}  \right) ds \\
=&  \lambda^{-1} \sum_{\ell=1}^m \int_0^1 (\tau_2\otimes \Tr_{d_\ell}) \left( F_\ell^{1-2s} \mathbf{D}^{(\ell) s} \Pi^{(\ell)} F_\ell^{1-2s}  \mathbf{X}^{(\ell)}\mathbf{D}^{(\ell) (1-s)} \mathbbm{1}^{(\ell)}  \mathbf{Y}^{(\ell)*}  \right) ds\\
=&  \lambda^{-1} \sum_{\ell=1}^m \frac{1}{d_\ell} \int_0^1 (\tau_2\otimes \Tr_{d_\ell}) \left( F_\ell^{1-2s} \mathbf{D}^{(\ell) s} \Pi^{(\ell)} F_\ell^{1-2s}  \mathbf{X}^{(\ell)}  \mathbbm{1}^{(\ell)} \mathbf{D}^{(\ell) (1-s)} \mathbbm{1}^{(\ell)}  \mathbf{Y}^{(\ell)*}  \right) ds\\
=&  \sum_{\ell=1}^m  (\tau_1\otimes \Tr_{d_\ell}) \left(  \K_{D, \Pi}^{(\ell)} (\widehat{\iota}^{(\ell)} (x_k^{(\ell)})_{k, \ell} )  (\widehat{\iota}^{(\ell)} (y_k^{(\ell)})_{k, \ell} )^*\right)\\
=&  \left\langle \K_{D, \Pi}  (\widehat{\iota} (x_k^{(\ell)})_{k, \ell} ) , \widehat{\iota} (y_k^{(\ell)})_{k, \ell} \right\rangle\\
=&  \left\langle  \widehat{\iota}^* \K_{D, \Pi}  (\widehat{\iota} (x_k^{(\ell)})_{k, \ell} ) ,  (y_k^{(\ell)})_{k, \ell} \right\rangle.
\end{align*}
This completes the computation.
\end{proof}

\begin{theorem}\label{thm:laplacian3}
Suppose that $D\in \cM$ is strictly positive.
Then 
\begin{align*}
\cL_a^*(D)=\frac{\lambda^{-1/2}}{2}\Div_0 \widehat{\iota}^* \K_{D, \Pi} \widehat{\iota }(\nabla_0 \log D)-
 \lambda^{-1/2}\Div_0 \widehat{\iota}^* \K_D \bigoplus_{\ell=1}^m\left ((\log \mu_k^{(\ell)}\mu_r^{(\ell)})  \mathfrak{F}^{-1}(\widehat{\Delta}_{k,r}^{(\ell)})\right)_{k, r=1}^{d_\ell}.
\end{align*}
\end{theorem}
\begin{proof}
It follows from Proposition \ref{prop:reform}.
\end{proof}

Let 
\begin{align*}
\Div= 2^{-1/2} \lambda^{-1/4}\Div_0\widehat{\iota}^*, \quad \nabla=2^{-1/2} \lambda^{-1/4}\widehat{\iota}\nabla_0.
\end{align*}

We have that 
\begin{theorem}\label{thm:laplacian4}
Suppose that $D\in \cM$ is strictly positive.
Then 
\begin{align*}
\cL_a^*(D)=\Div \K_{D, \Pi}(\nabla \log D)-
2^{1/2} \lambda^{-1/4}\Div \K_D \bigoplus_{\ell=1}^m\left ((\log \mu_k^{(\ell)}\mu_r^{(\ell)})  \mathfrak{F}^{-1}(\widehat{\Delta}_{k,r}^{(\ell)})\right)_{k, r=1}^{d_\ell}.
\end{align*}
\end{theorem}
\begin{proof}
It follows the notation of divergences and gradients.
\end{proof}

\section{Gradient Flow}

The gradient flow equation for the bimodule KMS symmetric quantum Markov semigroups is quite different from the one for the bimodule GNS symmetric quantum Markov semigroups.
In this section, we shall setup the gradient flow equation under the assumption that $\cR(\widehat{\cL}_0)=\cR(\overline{\widehat{\cL}_0})$.

Suppose that $X=(x_{k,r}^{(\ell)})_{k, r, \ell}, Y=(y_{k,r}^{(\ell)})_{k, r, \ell}\in \oM$.
The sesquilinear form $\langle X, Y \rangle_{D, \widehat{\Delta}}$ on $\oM$ is defined as follows:
\begin{align*}
    \langle X, Y\rangle_{D, \widehat{\Delta}} = \langle \K_{D,\Pi}  X, Y\rangle,
\end{align*}
where $\langle \cdot , \cdot\rangle$ is the usual inner product on $\oM$ induced by $\tau_1$.

\begin{proposition}
The sesquilinear form $\langle \cdot, \cdot \rangle_{D, \widehat{\Delta}}$ is an inner product on $\oM$.
\end{proposition}
\begin{proof}
It follows from Proposition \ref{prop:reform} that $\K_{D, \Pi}$ is positive definite.
\end{proof}

We shall denote by $\|\cdot\|_{D, \widehat{\Delta}}$ is the norm induced from the inner product  $\langle \cdot, \cdot\rangle_{D, \widehat{\Delta}}$ .

\begin{proposition}\label{prop:dualgradient}
Suppose that $x\in \cM$.
Then 
\begin{align*}
\left(\Div \K_{D, \Pi} \nabla x\right)^* =\Div  \K_{D, \Pi} \nabla x^*.
\end{align*}
\end{proposition}
\begin{proof}
By Theorem \ref{thm:laplacian2} and Propositions \ref{prop:kd}, \ref{prop:reform}, we have that 
\begin{align*}
& (\Div  \K_{D, \Pi} \nabla x^*)^* \\
= & \frac{\lambda^{-1/2}}{2} \left(\sum_{k,r, \ell}  \partial_k^{(\ell)*} \K_{D, k, r, \ell} \left( x^* \mathfrak{F}^{-1}(\widehat{\Delta}_{k,r}^{(\ell)}) -  \mathfrak{F}^{-1}(\widehat{\Delta}_{k,r}^{(\ell)}) x^* \right)\right)^* \\
= & \frac{\lambda^{-1/2}}{2} \sum_{k,r, \ell^* \neq \ell} \bE_{\cM} \partial_k^{(\ell)} \K_{D, k, r, \ell^*} \left( x \mathfrak{F}^{-1}(\widehat{\Delta}_{k,r}^{(\ell)})^* -  \mathfrak{F}^{-1}(\widehat{\Delta}_{k,r}^{(\ell)})^* x \right) \quad \text{Equation \eqref{eq:adder}}\\
&+  \frac{\lambda^{-1/2}}{2} \sum_{k,r, \ell^* = \ell}  \bE_{\cM} \partial_{k}^{(\ell)} \K_{D, k^*, r^*, \ell} \left( x \mathfrak{F}^{-1}(\widehat{\Delta}_{k,r}^{(\ell)})^* -  \mathfrak{F}^{-1}(\widehat{\Delta}_{k,r}^{(\ell)})^* x \right) \\
= & \frac{\lambda^{-1/2}}{2} \sum_{k,r, \ell^* \neq \ell} \bE_{\cM} \partial_k^{(\ell)} \K_{D, k, r, \ell^*} \left( x \mathfrak{F}^{-1}(\widehat{\Delta}_{k,r}^{(\ell^*)})-  \mathfrak{F}^{-1}(\widehat{\Delta}_{k,r}^{(\ell^*)}) x \right) \\
&+  \frac{\lambda^{-1/2}}{2} \sum_{k,r, \ell^* = \ell}  \bE_{\cM} \partial_{k}^{(\ell)} \K_{D, k^*, r^*, \ell} \left( x \mathfrak{F}^{-1}(\widehat{\Delta}_{k^*,r^*}^{(\ell)}) -  \mathfrak{F}^{-1}(\widehat{\Delta}_{k^*,r^*}^{(\ell)}) x \right) \\
= & \frac{\lambda^{-1/2}}{2} \sum_{k,r, \ell^* \neq \ell} \bE_{\cM} \partial_k^{(\ell^*)} \K_{D, k, r, \ell} \left( x \mathfrak{F}^{-1}(\widehat{\Delta}_{k,r}^{(\ell)})-  \mathfrak{F}^{-1}(\widehat{\Delta}_{k,r}^{(\ell)}) x \right) \\
&+  \frac{\lambda^{-1/2}}{2} \sum_{k,r, \ell^* = \ell}  \bE_{\cM} \partial_{k^*}^{(\ell)} \K_{D, k, r, \ell} \left( x \mathfrak{F}^{-1}(\widehat{\Delta}_{k,r}^{(\ell)}) -  \mathfrak{F}^{-1}(\widehat{\Delta}_{k,r}^{(\ell)}) x \right) \\
= & \frac{\lambda^{-1/2}}{2} \sum_{k,r, \ell}  \partial_k^{(\ell)*} \K_{D, k, r, \ell} \left( x\mathfrak{F}^{-1}(\widehat{\Delta}_{k,r}^{(\ell)}) -  \mathfrak{F}^{-1}(\widehat{\Delta}_{k,r}^{(\ell)}) x \right) \\
=& \Div  \K_{D, \Pi} \nabla x.
\end{align*}
This completes the proof of the proposition.
\end{proof}

\begin{theorem}\label{thm:gradient}
Suppose that the bimodule KMS quantum Markov semigroup $\{\Phi_t\}_{t \geq 0}$ is relative ergodic such that $\cR(\widehat{\cL}_0)=\cR(\overline{\widehat{\cL}_0})$ and $\{D_s\}_{s\in (-\epsilon, \epsilon)}$ is a continuously differentiable family of density operators in $\cM$ with $\bE_{\cN}(D_s)$ is independent of $s$ with $D_0=D$ and $\epsilon>0$.
Then there exists $X\in \oM$ such that
\begin{align*}
  \left. \frac{d D_s}{ds}\right|_{s=0}=\dot{D}_0=\Div \K_{D, \Pi} X.
\end{align*}
If $X$ is minimal subject to the induced norm $\|\cdot\|_{D, \widehat{\Delta}}$, then there exists a unique self-adjoint element $x$ in $\cM$ such that $\bE_{\cN}(x)=0$ and $X=\nabla x$.
\end{theorem}
\begin{proof}
Note that $\bE_{\cN}(\dot{D}_0)=0$.
By Proposition \ref{prop:supportrange}, We see that $\dot{D}_0\in \Ran(\Div)$ and there exists $V\in \oM$ such that $\dot{D}_0=\Div V$.
Let $X_0=\K_{D, \Pi}^{-1}V$.
We see that $\dot{D}_0=\Div \K_{D, \Pi} X_0.$

Suppose that $X\in \oM$ is minimal subject to the norm $\|\cdot\|_{D, \widehat{\Delta}}$ and $A\in \oM$ such that $\Div A=0$. 
Let $Y\in \oM$ such that $ Y= X+r \K_{D, \Pi}^{-1} A$ with $r\in \bR$.
We have that 
\begin{align*}
\Div  \K_{D, \Pi } Y 
 =& \Div \K_{D, \Pi}  X+r\Div\K_{D ,\Pi}  \K_{D, \Pi}^{-1} A \\
 =&  \dot{D}_0 +r\Div A \\
 =& \dot{D}_0 
 \end{align*}
  and $\langle X, X\rangle_{D, \widehat{\Delta}}\leq \langle Y, Y\rangle_{D, \widehat{\Delta}}$.
This implies that for any $r\in \bR$, 
\begin{align*}
0 \leq 2 \Re \langle  X, A \rangle +r \langle A, A \rangle_{D, \widehat{\Delta}}. 
\end{align*}
Hence $0=\langle  X, A\rangle$ and $X\perp \Ker \Div$, i.e. $X\in \Ran \nabla$.
Therefore, there exists $x_1\in \cM$ such that $ \nabla x_1=X$.
Let $\displaystyle \mathbb{P}=\lim_{m \to\infty} (\nabla \Div)^{1/m}$ and $\widetilde{X}=\mathbb{P}(X_0)$.
We see that $\mathbb{P}: \oM \to \oM$, $\widetilde{X}\in \oM$ and $\widetilde{X} \perp \Ker\Div$.
This shows that $X=\widetilde{X}\in \oM$.

By Proposition \ref{prop:dualgradient}, we see that 
\begin{align*}
 \dot{D}_0=\dot{D}_0 ^*=\left( \Div \K_{D ,\Pi} \nabla x_1 \right)^*
    =\Div  \K_{D ,\Pi} \nabla x_1^* 
\end{align*}
i.e.  $\displaystyle \dot{D}_0= \Div  \K_{D, \Pi} (\nabla x_1^*)$.
Let $\displaystyle x_2=\frac{1}{2}(x_1+x_1^*)$.
We have that 
\begin{align*}
\dot{D}_0= \Div \K_{D, \Pi} \nabla x_2 .
\end{align*}
By letting $x=x_2-\bE_{\cN}(x_2)$, we obtain that $\displaystyle \dot{D}_0= \Div \K_{D, \Pi}  \nabla x$.

Suppose that $y$ is a solution to the equation subject to $y^*=y$ and $\bE_{\cN}(y)=0$.
Then $\nabla (x-y)=0$ and $x-y\in \cN$.
This shows that $x=y$.
\end{proof}


\begin{lemma}\label{lem:kmsgradient}
Suppose that $\bfx\in \cN'\cap \cM$.
Then we have that $i[\bfx, D]\in \Ran \Div$.
\end{lemma}
\begin{proof}
For any $v\in \cN$, we have that 
\begin{align*}
\tau(v\bE_{\cN}( \bfx D-D\bfx))=\tau(v\bfx D-vD\bfx)=\tau(\bfx v D-v D \bfx)=0.
\end{align*}
This shows that $\bE_{\cN}( \bfx D-D\bfx)=0$, i.e. $ \bfx D-D\bfx \in \Ran\Div$.
\end{proof}

\begin{proposition}
There exists a density operator $D_\Delta\in \cM$ such that 
\begin{align*}
 i[\bfy_1, D]  +2^{1/2}\lambda^{-1/4}\Div \K_D\bigoplus_{\ell=1}^m \left((\log \mu_k^{(\ell)} \mu_r^{(\ell)})  \mathfrak{F}^{-1}(\widehat{\Delta}_{k,r}^{(\ell)})\right)_{k, r=1}^{d_\ell} 
 =\Div  \K_{D,\Pi} \nabla \log D_\Delta .
\end{align*}
\end{proposition}
\begin{proof}
By Theorem \ref{thm:laplacian2} and Lemma \ref{lem:kmsgradient}, we have that the left hand side is self-adjoint followed from Proposition \ref{prop:dualgradient} and it vanishes on the conditional expectation $\bE_{\cN}$.
Denote the left hand side by $\widetilde{x}\in \cM$.
There exists $A\in \oM$ such that $\widetilde{x}=\Div A$.
By Theorem \ref{thm:gradient}, there exists self-adjoint element $x_{\Delta}\in \cM$ such that $ \displaystyle 2^{1/2}\lambda^{-1/4}\Div \K_D\bigoplus_{\ell=1}^m \left((\log \mu_k^{(\ell)} \mu_r^{(\ell)})  \mathfrak{F}^{-1}(\widehat{\Delta}_{k,r}^{(\ell)})\right)_{k, r=1}^{d_\ell} 
 =\Div  \K_{D,\Pi} \nabla x_{\Delta}.$
Let $D_\Delta =e^{x_{\Delta}}\in \cM$.
Then the proposition is true.
\end{proof}

We shall call $D_\Delta$ the hidden density for the bimodule quantum semigroup $\{\Phi_t\}_{t\geq 0}$.
Hence we have the following theorem.
\begin{theorem}\label{thm:sqrt}
Suppose that $D\in \cM$ is strictly positive.
Then 
\begin{align*}
 \cL^*(D)  =  \Div \K_{D, \Pi}  \left( \nabla \log D -\nabla \log D_{\Delta}\right).
\end{align*} 
\end{theorem}
\begin{proof}
    It follows from Theorem \ref{thm:laplacian4} and the previous computation.
\end{proof}


\begin{remark}
Suppose that $\overline{\widehat{\Delta}}=\widehat{\Delta}^{-1}$.
We have that $\overline{\widehat{\cL}_{\Delta}}=\widehat{\cL}_{\Delta}$.
In this case $\cR(\overline{\cL}_0)$ might not be equal to $\cR(\overline{\cL}_0)$ and the element $U^{(\ell)}$ is not a unitary.
To obtain a gradient flow equation in Theorem \ref{thm:gradient}, one should modify the divergence and gradient such that they are supported on a range projection equivalent to the range projection $\cR(\widehat{\cL}_0)$.
We shall discuss the case in the future.

\end{remark}

\section{The inclusion $\bC^n \subset  M_n(\bC)$}

In this section, we shall obtain some bimodule KMS symmetric semigroups for the inclusion $\bC^n \subset M_n(\bC)$.
We denote the Schur product of the matrices by $\diamond$ for convenience.
 
Let $(E_{j,k})_{j,k=1}^n$ be a system of matrix units of $M_n(\bC)$ such that $E_{k,k}\Omega$, $k=1, \ldots, n$ spans $\bC^n$.
Let $|j\rangle=\sqrt{n}E_{j,j}\Omega$.
Then $\{|j \rangle\}_{j=1}^n$ is an orthonormal basis of $\bC^n$.
We shall identify $L^2(\cM)$ with $\bC^n \otimes \bC^n$ via the following unitary transformation:
\begin{align*}
 E_{j,k}\Omega \mapsto n^{-1/2} |j\rangle   |k \rangle. 
\end{align*}
The left action of $E_{j,k}$ on $L^2(\cM)$ is identified as $E_{j,k}\otimes I$.
The Jones projection $\displaystyle e_1=\sum_{j=1}^n E_{j,j}\otimes E_{j,j}$ and its trace is $\displaystyle \frac{1}{n}$.
The basic construction $\cM_1$ is spanned by $\{ E_{t,s}\otimes E_{k,k}: t,s, k=1, \ldots, n\}$.
The conjugation $J$ satisfies $J (E_{j,k}\otimes E_{t,s})^* J =E_{s,t}\otimes E_{k,j}$.
The relative commutant $\cN'\cap \cM_1$ is spanned by $\{E_{j,j}\otimes E_{k,k}: j=1, \ldots, n\}$.
We shall identify $L^2(\cM_1)$ with $\bC^n \otimes \bC^n \otimes \bC^n$ via the following unitary transformation: 
\begin{align*}
E_{t,s}\otimes E_{k,k}\Omega_1\mapsto n^{-1} |t\rangle |s\rangle|k\rangle .
\end{align*}
The left action of $E_{t,s}\otimes E_{k,k}$ on $L^2(\cM_1)$ is $E_{t, s}\otimes I \otimes E_{k,k}$.
The Jones projection $\displaystyle e_2=I\otimes I \otimes \left( \frac{1}{n} \sum_{j,k=1}^n E_{j,k}\right)$.
The basic construction $\cM_2$ is spanned by $\{E_{j,k}\otimes I \otimes E_{t,s}: j,k,t,s=1, \ldots, n\}$.
The relative commutant $\cM'\cap \cM_2$ is spanned by $\{I\otimes I\otimes E_{j,k}: j,k=1, \ldots, n\}$.
The conditional expectation $\bE_{\cM_1'}$ is defined as 
\begin{align*}
    \bE_{\cM_1'}(X)=\sum_{k=1}^n (I\otimes I \otimes E_{k,k}) X (I\otimes I \otimes E_{k,k}),
\end{align*}
where $X\in \cM_2$.
The conjugation $J_1$ satisfies $J_1(E_{j,k}\otimes I \otimes E_{t,s})^* J_1=I\otimes E_{k,j}\otimes E_{s,t}$.

The Fourier transform $\fF: \cN'\cap \cM_1 \to \cM'\cap \cM_2$ is 
\begin{align*}
\fF(E_{j,j}\otimes I \otimes E_{k,k}) =\frac{1}{\sqrt{n}} I\otimes I \otimes E_{k,j}
\end{align*}
and the convolution $I\otimes I \otimes E_{t,s}$ and $I\otimes I \otimes E_{j,k}$ is as follows
\begin{align*}
    (I\otimes I \otimes E_{j,k})*( I\otimes I \otimes E_{t,s})=\sqrt{n}\delta_{j,t}\delta_{k,s}E_{j,k}.
\end{align*}

\begin{lemma}
Suppose that $\Phi: \cM\to\cM$ is a bimodule map.
Then there exists a matrix $A\in M_n(\bC)$ such that $\Phi(X)=X \diamond A$.
\end{lemma}
\begin{proof}
Suppose that $\widehat{\Phi}=I\otimes I \otimes Y$ for some $Y=(Y_{j,k})_{j,k=1}^n\in M_n(\bC)$.
We have that
\begin{align*}
& n^{5/2}\bE_{\cM}(e_2e_1 \widehat{\Phi} (E_{t,s}\otimes I \otimes I) e_1e_2 )\\
=& n^{1/2} \bE_{\cM}\left( \left(\sum_{j, k=1}^n E_{j,j}\otimes I \otimes E_{k,j}\right) (E_{t,s}\otimes I \otimes Y) \left(\sum_{j, k=1}^n E_{j,j}\otimes I \otimes E_{j,k}\right)\right)\\
=& n^{1/2} \bE_{\cM}\left( \left(\sum_{k=1}^n E_{t,t}\otimes I \otimes E_{k,t}\right) (E_{t,s}\otimes I \otimes Y) \left(\sum_{k=1}^n E_{s,s}\otimes I \otimes E_{s,k}\right)\right)\\
=& n^{1/2} E_{t,s} Y_{t,s}.
\end{align*}
This implies that $\Phi(X)=\sqrt{n} X\diamond Y$.
By taking $A=\sqrt{n}Y$, we have that $\Phi(X)=X\diamond A$.

\end{proof}

\begin{remark}
Suppose that $\displaystyle \widehat{\Phi}=I\otimes I \otimes \sum_{j,k=1}^n \widehat{\Phi}_{j,k}E_{j,k}$.
We have that 
\begin{align*}
1*\widehat{\Phi}= n^{1/2} \sum_{j=1}^n \widehat{\Phi}_{j,j} E_{j,j}.    
\end{align*}
\end{remark}

\begin{remark}
Suppose that $\Phi: \cM \to \cM$ given by $\Phi(X)=X\diamond A$ for some $A\in M_n(\bC)$.
Suppose that $D=\diag(t_1, \ldots, t_n)$ is the diagonal matrix in $M_n(\bC)$.
We have that for any $X\in M_n(\bC)$,
\begin{align*}
\Phi(DX)=(DX)\diamond A=D(X\diamond A)=D\Phi(X).
\end{align*}
Similarly, we have that $\Phi(XD)=\Phi(X)D$.
This implies that $\Phi$ is a bimodule map for the inclusion $\bC^n \subset M_n(\bC)$.
\end{remark}

\begin{remark}
Suppose that $\Phi$ is a bimodule map.
Then $\Phi$ is positive if and only if $\Phi$ is completely positive if and only if $\widehat{\Phi}\geq 0$.
We see that $\Phi$ is unital if and only if $\displaystyle \widehat{\Phi}_{j,j}=\frac{1}{\sqrt{n}}$, $j=1, \ldots, n$.    
The contragredient $\overline{\widehat{\Phi}}$ of $\widehat{\Phi}$ is the transpose of $\widehat{\Phi}$.
We obtain that $\Phi^*$ is unital if $\Phi$ is unital.
Hence the unital bimodule map $\Phi$ is equilibrium for the trace $\tau$.
\end{remark}

\begin{proposition}
Suppose $\bC^n \subset M_n(\bC)$ is the finite inclusion and $\rho$ is a faithful state on $M_n(\bC)$.
We have that 
\begin{align*}
\widehat{\Delta}_\rho=\Delta_\rho^{\mathsf{T}} \diamond \Delta_\rho^{-1},
\end{align*}
where $\Delta_\rho^{\mathsf{T}} \diamond \Delta_\rho^{-1}$ is viewed as an element in $\cM'\cap \cM_2$.
\end{proposition}
\begin{proof}
Recall that $\widehat{\Delta}_\rho=\lambda^{-1/2}\mathfrak{F}( \bE_{\cN'} (\Delta_{\rho} e_1 \Delta_{\rho}^{-1} ))$.
Suppose that $\displaystyle \Delta_\rho=\sum_{j,k=1}^n a_{j,k} E_{j, k}$ and $\displaystyle \Delta_\rho^{-1}=\sum_{j,k=1}^n b_{j,k} E_{j, k}$.
We have that 
\begin{align*}
\Delta_{\rho} e_1 \Delta_{\rho}^{-1} 
=& ((\sum_{j,k=1}^n a_{j,k} E_{j, k})\otimes I \otimes I ) (\sum_{k=1}^n E_{k,k}\otimes I \otimes E_{k,k}) ((\sum_{j,k=1}^n b_{j,k} E_{j, k}) \otimes I \otimes I ) \\
=&\sum_{t,s, k=1}^n  a_{t,k} b_{k,s}E_{t,s}\otimes I \otimes E_{k,k},
\end{align*}
and
\begin{align*}
\bE_{\cN'}(\Delta_{\rho} e_1 \Delta_{\rho}^{-1})= \sum_{j,k=1}^n a_{j,k}b_{k,j} E_{j,j}\otimes I \otimes E_{k,k}.
\end{align*}
Hence
\begin{align*}
\widehat{\Delta}_\rho
=& n^{1/2}\mathfrak{F}( \bE_{\cN'} (\Delta_{\rho} e_1 \Delta_{\rho}^{-1} ))\\
=&  \sum_{j,k=1}^n a_{j,k}b_{k,j} I \otimes I \otimes E_{k,j}\\
=& \Delta_\rho^{\mathsf{T}} \diamond \Delta_\rho^{-1}.
\end{align*}
This completes the computation.
\end{proof}

In the following, we shall present an example of KMS symmetric semigroup on $M_3(\bC).$
Let $\omega=e^{2\pi i/3}$.
Suppose that 
\begin{align*}
    \widehat{\Delta}=\frac{1}{3}\begin{pmatrix}
        1 & 1 & 1 \\
        1 & 1 & 1 \\
        1& 1 & 1
    \end{pmatrix} + \frac{\mu}{3}\begin{pmatrix}
        1 & \omega & \omega^2 \\
        \omega^2 & 1 & \omega \\
        \omega & \omega^2 & 1
    \end{pmatrix}
    +\frac{\mu^{-1}}{3}\begin{pmatrix}
        1 & \omega^2 & \omega \\
        \omega & 1 & \omega^2 \\
        \omega^2 & \omega & 1
    \end{pmatrix}\in \cM'\cap \cM_2,
\end{align*}
where $\mu>0$ and $\mu\neq 1$.
We have that 
\begin{align*}
\overline{\widehat{\Delta}}=\widehat{\Delta}^{-1}, \quad \widehat{\Delta} e_2=e_2.
\end{align*}
Let
\begin{align*}
H=\frac{1}{2}\begin{pmatrix}
        1 & -1 & 0 \\
        -1 & 1 & 0 \\
        0 & 0 & 0
    \end{pmatrix}
    + \frac{1}{2}\begin{pmatrix}
        1 & 0 & -1 \\
        0 & 0 & 0 \\
        -1 & 0 & 1
    \end{pmatrix}\in \cM'\cap \cM_2,
\end{align*}
and 
\begin{align*}
\widehat{\cL}_0=\widehat{\Delta}^{1/2}H\widehat{\Delta}^{1/2}.
\end{align*}
Precisely,
\begin{align*}
\widehat{\cL}_0=\frac{\mu}{3} \begin{pmatrix} 1 & \omega^2 & \omega \\ \omega & 1 & \omega^2 \\ \omega^2 & \omega & 1 \end{pmatrix}
+ \frac{\mu^{-1}}{3} \begin{pmatrix} 1 & \omega & \omega^2 \\ \omega^2 & 1 & \omega \\ \omega & \omega^2 & 1 \end{pmatrix}
+ \frac{1}{6} \begin{pmatrix} 2 & -1 & -1 \\ -1 & -1 & 2 \\ -1 & 2 & -1 \end{pmatrix}.
\end{align*}
We have that 
\begin{align*}
He_2=0, \quad \overline{H}=H, \quad \overline{\widehat{\cL}_0}  =\overline{\widehat{\Delta}} \widehat{\cL}_0 \overline{\widehat{\Delta}}= \widehat{\Delta}^{-1/2} H\widehat{\Delta}^{-1/2}.
\end{align*}
Now taking the convolution, we obtain that 
\begin{align*}
    1*\widehat{\cL}_0=\sqrt{3} \left( \frac{1}{3} E_{1,1} -\frac{1}{6}E_{2,2}-\frac{1}{6}E_{3,3} +\frac{\mu+\mu^{-1}}{3}I\right).
\end{align*}
Let $\displaystyle h=\frac{1}{2}\frac{\mu^{-1/2}-\mu^{1/2}}{\mu^{-1/2}+\mu^{1/2}} (E_{2,2}-E_{3,3})$.
Note that 
\begin{align*}
    & \overline{\widehat{\Delta}} (\mu^{-1/2}(E_{1,1}+\omega E_{2,2}+\omega^2 E_{3,3})+\mu^{1/2}(E_{1,1}+\omega^2 E_{2,2}+\omega E_{3,3})) e_2 \\
    =& (\mu^{1/2}(E_{1,1}+\omega E_{2,2}+\omega^2 E_{3,3})+\mu^{-1/2}(E_{1,1}+\omega^2 E_{2,2}+\omega E_{3,3})) e_2.
\end{align*}
We see that
\begin{align*}
  & \overline{\widehat{\Delta}} ((\mu^{1/2}+\mu^{-1/2})(E_{1,1}-\frac{1}{2}E_{2,2}-\frac{1}{2}E_{3,3})+ i\sqrt{3}(\mu^{-1/2}-\mu^{1/2})(E_{2,2}-E_{3,3})) e_2 \\
    =&  ((\mu^{1/2}+\mu^{-1/2})(E_{1,1}-\frac{1}{2}E_{2,2}-\frac{1}{2}E_{3,3})- i\sqrt{3}(\mu^{-1/2}-\mu^{1/2})(E_{2,2}-E_{3,3}))  e_2.   
\end{align*}
Hence
\begin{align*}
   \overline{ \widehat{\Delta}}\left( \frac{1}{2}(1*\widehat{\cL}_0)+ ih\right)e_2=\left(\frac{1}{2}(1*\widehat{\cL}_0)- ih\right)e_2.
\end{align*}
By Proposition \ref{prop:kms1}, we see that the semigroup arising from the generator $\cL$ given by
\begin{align*}
    \widehat{\cL}=\vcenter{\hbox{\begin{tikzpicture}[scale=0.65]
    \begin{scope}[shift={(0,1.5)}]
    \draw [blue] (-0.5, 0.8)--(-0.5, 0) .. controls +(0, -0.6) and +(0,-0.6).. (0.5, 0)--(0.5, 0.8);    
\begin{scope}[shift={(0.5, 0.3)}]
\end{scope}
    \end{scope}
\draw [blue] (-0.5, -0.8)--(-0.5, 0) .. controls +(0, 0.6) and +(0,0.6).. (0.5, 0)--(0.5, -0.8);
\begin{scope}[shift={(0.5, -0.3)}]
\draw [fill=white] (-0.3, -0.3) rectangle (0.3, 0.3);
\node at (0, 0) {\tiny $\mathbf{y}^*$};
\end{scope}
\end{tikzpicture}}}
+
\vcenter{\hbox{\begin{tikzpicture}[scale=0.65]
    \begin{scope}[shift={(0,1.5)}]
    \draw [blue] (-0.5, 0.8)--(-0.5, 0) .. controls +(0, -0.6) and +(0,-0.6).. (0.5, 0)--(0.5, 0.8);    
\begin{scope}[shift={(0.5, 0.3)}]
\draw [fill=white] (-0.3, -0.3) rectangle (0.3, 0.3);
\node at (0, 0) {\tiny $\mathbf{y}$};
\end{scope}
    \end{scope}
\draw [blue] (-0.5, -0.8)--(-0.5, 0) .. controls +(0, 0.6) and +(0,0.6).. (0.5, 0)--(0.5, -0.8);
\begin{scope}[shift={(0.5, -0.3)}]
\end{scope}
\end{tikzpicture}}}
-  \vcenter{\hbox{\begin{tikzpicture}[scale=0.6]
        \draw [blue] (0.2, -0.8) --(0.2, 0.8);
        \draw [blue] (-0.2, -0.8) --(-0.2, 0.8);
        \draw [fill=white] (-0.5, -0.4) rectangle (0.5, 0.4);
    \node at (0, 0) {\tiny $\widehat{\cL}_0$};
        \end{tikzpicture}}},
\end{align*} 
where $\displaystyle \mathbf{y}=\frac{1}{2}(1*\widehat{\cL}_0)-ih$, is a bimodule KMS symmetric quantum semigroup.
The diagonals of $\widehat{\cL}$ are zero. 
The off-diagonal entries are 
\begin{align*}
\widehat{\cL}_{1,2} &= \frac{\sqrt{3}(4\mu+4\mu^{-1}+1) + 2(\mu+\mu^{-1}+1)}{12}  -i \left( \frac{\sqrt{3}(\mu-\mu^{-1})}{6} - \frac{1-\mu}{2(1+\mu)} \right), \\[1em]
\widehat{\cL}_{1,3} &= \frac{\sqrt{3}(4\mu+4\mu^{-1}+1) + 2(\mu+\mu^{-1}+1)}{12} + i \left( \frac{\sqrt{3}(\mu-\mu^{-1})}{6} - \frac{1-\mu}{2(1+\mu)} \right) ,\\[1em]
\widehat{\cL}_{2,1} &= \frac{\sqrt{3}(4\mu+4\mu^{-1}+1) + 2(\mu+\mu^{-1}+1)}{12} + i \left( \frac{\sqrt{3}(\mu-\mu^{-1})}{6} + \frac{1-\mu}{2(1+\mu)} \right), \\[1em]
\widehat{\cL}_{3,1} &= \frac{\sqrt{3}(4\mu+4\mu^{-1}+1) + 2(\mu+\mu^{-1}+1)}{12} - i \left( \frac{\sqrt{3}(\mu-\mu^{-1})}{6} + \frac{1-\mu}{2(1+\mu)} \right), \\[1em]
\widehat{\cL}_{2,3} &= \frac{\sqrt{3}(2\mu+2\mu^{-1}-1) + \mu+\mu^{-1}-2}{6} - i\frac{\sqrt{3}(\mu-\mu^{-1})}{6}, \\[1em]
\widehat{\cL}_{3,2} &= \frac{\sqrt{3}(2\mu+2\mu^{-1}-1) + \mu+\mu^{-1}-2}{6} +i\frac{\sqrt{3}(\mu-\mu^{-1})}{6}.
\end{align*}
Note that $\cR(\widehat{\cL}_0)=\cR(H)=\cR(\overline{\widehat{\cL}_0})=1-e_2$.
We see that this semigroup is relatively ergodic.
By Theorem \ref{thm:kmslimit}, we see that this semigroup is not KMS symmetry due to the fact that the limit is a multiple of the identity.

Let $\cA$ be the algebra generated by $H$ and $\widehat{\Delta}(1-e_2)$.
Then $\cA$ is isomorphic to $M_2(\bC)$.
Note that $\displaystyle H= \frac{3}{2}P_{1,1} + \frac{1}{2}P_{2,2}$, where
\begin{align*}
  P_{1,1} = \begin{pmatrix}
  \frac{2}{3} & -\frac{1}{3} & -\frac{1}{3} \\
  -\frac{1}{3} & \frac{1}{6} & \frac{1}{6} \\
  -\frac{1}{3} & \frac{1}{6} & \frac{1}{6}
  \end{pmatrix}, \quad 
    P_{2,2} = \begin{pmatrix}
  0 & 0 & 0 \\[2pt]
  0 & \frac{1}{2} & -\frac{1}{2} \\
  0 & -\frac{1}{2} & \frac{1}{2}
  \end{pmatrix}.
\end{align*}
Let
\begin{align*}
Q_{1,1}= \frac{1}{3}\begin{pmatrix}
        1 & \omega & \omega^2 \\
        \omega^2 & 1 & \omega \\
        \omega & \omega^2 & 1
    \end{pmatrix},\quad 
Q_{1,2} =\frac{1}{3}\begin{pmatrix}
        1 & \omega^2 & \omega \\
        \omega^2 &  \omega & 1 \\
        \omega & 1 & \omega^2 
    \end{pmatrix},
     \quad  Q_{2,2}= \frac{1}{3}\begin{pmatrix}
        1 & \omega^2 & \omega \\
        \omega & 1 & \omega^2 \\
        \omega^2 & \omega & 1
    \end{pmatrix}.
\end{align*}
We obtain that 
\begin{align*}
P_{1,1} =& \frac{1}{2} (Q_{1,1}+Q_{1,2}+Q_{2,1}+Q_{2,2}), \\
P_{2,2} =& \frac{1}{2} (Q_{1,1}-Q_{1,2}-Q_{2,1}+Q_{2,2}).
\end{align*}
Let $\displaystyle U=\frac{1}{\sqrt{2}} (Q_{1,1}+Q_{1,2}-Q_{2,1}+Q_{2,2})$.
Hence the directional matrix $\Pi$ is 
\begin{align*}
\Pi =\begin{pmatrix}
  Q_{1,1} & \frac{1}{2}Q_{1,2} \\ 
\frac{1}{2} Q_{2,1} &Q_{2,2}
\end{pmatrix},
\end{align*}
and 
   \begin{align*}
&\K_{D, \Pi}(x_{k,r})_{k,r=1}^2\\
=& \frac{1}{2} \int_0^1 F^{1-2s} \mathbf{D}^{ s} \Pi \left(F_\ell^{1-2s}  (x_{k,r})_{k,r=1}^{2}\right) \mathbf{D}^{ (1-s)}  ds\\
=&\frac{3}{2} (\bE_{\cM_1} \otimes I )
\int_0^1 
\begin{pmatrix}
  \mu^{2-4s} D^sQ_{1,1}  & \frac{1}{2}D^sQ_{1,2}  \\ 
\frac{1}{2} D^sQ_{2,1}& \mu^{-(2-4s)}D^sQ_{2,2}
\end{pmatrix}
\begin{pmatrix}
    x_{1,1}e_2D^{1-s} & x_{1,2} e_2D^{1-s} \\
    x_{2,1}e_2D^{1-s} & x_{2,2} e_2D^{1-s}
\end{pmatrix}
\,ds
.
\end{align*}

\begin{remark}
By taking $\displaystyle H=\frac{2}{\mu+\mu^{-1}} \begin{pmatrix}
        1 & -1/2 & -1/2 \\
        -1/2 & 1 & -1/2 \\
        -1/2 & -1/2 & 1
    \end{pmatrix}$, we shall obtain a bimodule GNS symmetric quantum semigroup.
\end{remark}

\begin{remark}
The above construction of bimodule KMS symmetric semigroups also works for any $n\geq 3$.
\end{remark}

\section{Relative Entropy}

In this section, we study the relative entropy functional with respect to the hidden density under the assumption that $\cR(\widehat{\cL}_0)=\cR(\overline{\widehat{\cL}_0})$. 
For a strictly positive density operator $D\in \cM_+$, we define 
\begin{align*}
    \cM_{+}(D) = \{\rho\in \cM_+\vert \bE_{\cN}(\rho) = \bE_{\cN}(D)\}.
\end{align*}
\begin{definition}[Riemannian Metric]
Let $D$ be a strictly positive density operator in $\cM_+$, and suppose $\{D_s\}_{s\in (a,b)}$ is a continuous path in $\cM_+(D)$ with $D_0 = D$, where $a<0< b$. 
The Riemannian metric $g_{\cL}$ on $\cM_+$ with respect to a bimodule KMS symmetric semigroup $\{\Phi_t\}_{t\geq 0}$ is defined as 
\begin{align*}
    \|\dot{D}_0\|_{g_{\cL}}=\min\left\{\|X\|_{D, \widehat{\Delta}}: \dot{D}_0= \Div \K_{D, \Pi} X \right\}.
\end{align*}
\end{definition}

Suppose $f:\cM_+\to \bR$ is a differentiable function.
For any self-adjoint $x\in \cM$ with $\tau(x)=0$, there exists $\displaystyle \frac{df}{dD}\in \cM$ such that 
\begin{align*}
    \lim_{s\to 0} \frac{f(D+sx)-f(D)}{s} =\tau\left(\frac{df}{dD} x\right).
\end{align*}

\begin{definition}[Gradient Vector Field]
Suppose $f$ is a differentiable function on $\cM_+$.
The gradient vector field $\grad_{g, D} f\in \left\{\nabla x: x=x^*\in \cM, \bE_{\cN}(x)=0\right\}\subset \oM$ is defined as
\begin{align*}
    \left.\frac{d}{ds}f(D_s)\right|_{s=0}
    = \left\langle \grad_{g, D} f , \nabla x \right\rangle_{D, \widehat{\Delta}},
\end{align*}
for any smooth path $\{D_s\}_{(-\epsilon, \epsilon)}$ with $D_0=D$ such that $\displaystyle \dot{D}_0= \Div\K_{D, \Pi} \nabla x$.
Note that the gradient vector field is unique.
\end{definition}

\begin{lemma}\label{lem:vecterfieldform}
Suppose $f$ is a differentiable function. 
Then 
\begin{align*}
   \grad_{g, D}f= \nabla\frac{df}{dD},
\end{align*}
and 
\begin{align*}
    \| \grad_{g, D}f\|_{D, \widehat{\Delta}}^2= \left\langle \Div \K_{D, \Pi} \nabla \frac{df}{dD} , \frac{df}{dD}\right\rangle.
\end{align*}
\end{lemma}
\begin{proof}
For any differentiable path $\{D_s\}_s$ of density operators with $D_0=D$, we have that 
\begin{align*}
 \left.\frac{d}{ds}f(D_s)\right|_{s=0}
 =& \left\langle \frac{df}{dD},  \dot{D}_0\right\rangle \\
=&  \left\langle \frac{df}{dD},  \Div\K_{D, \Pi} \nabla x\right\rangle \\
=&\left\langle \nabla\frac{df}{dD}, \K_{D, \Pi} \nabla x \right\rangle \\
=&  \left\langle \nabla\frac{df}{dD} ,  \nabla x \right\rangle_{D, \widehat{\Delta}} .
\end{align*}
where the second equality follows from Theorem \ref{thm:gradient} and $x\in \cM$.
This implies that $\displaystyle  \grad_{g, D}f= \nabla\frac{df}{dD}$.
\end{proof}

Let $f(\rho)=H(\rho\| D_{\Delta})=\tau(\rho\log \rho-\rho\log D_{\Delta})$ be the relative entropy functional on $\cM_{+}(D)$. 
A differentiable path $\{D_s\}_s$ is a gradient flow (See \cite{Maa11} for more details) for $f$ if $\displaystyle \dot{D}_s= -\Div \K_{D_s,  \Pi} \left(\nabla\frac{df}{dD_s} \right)$. 
Then we have the following theorem:

\begin{theorem}\label{thm:kmsgradient}
Suppose that $\{\Phi_t\}_{t\geq 0}$ is a bimodule quantum Markov semigroup which is relatively ergodic and bimodule KMS symmetric with respect to $\widehat{\Delta}$.
Suppose $\{D_s\}_{s\in (-\epsilon, \epsilon)}$ is a differentiable path of density operators in $\cM_{+}(D)$ satisfying
\begin{align*}
    D_0 = D,\quad \dot{D}_s=-\cL^*(D_s).
\end{align*}
Then it is the gradient flow for the relative entropy functional $f(\rho) = H(\rho\| D_{\Delta})$.
\end{theorem}
\begin{proof}
For any $x\in \cM$ with $\bE_{\cN}(x)=0$, we have 
\begin{align*}
&     \left.\frac{d}{ds}f(D+sx)\right|_{s=0} \\
    =& \tau(x(\log D-\log D_\Delta))+\int_0^\infty \tau\left(\frac{D}{(D+r)^2}x\right)dr \\
    =& \tau(x(\log D-\log D_\Delta)).
\end{align*}
Hence $\displaystyle \frac{df}{dD}=\log D-\log D_{\Delta}$.
On the other hand, by Theorem \ref{thm:sqrt}, we have that 
\begin{align*}
 \dot{D}_s=& - \cL^*(D_s)\\
 = & -  \Div\K_{D_s, \Pi} \left( \nabla \log D_s  - \nabla \log D_\Delta \right) \\
 =& -  \Div\K_{D_s, \Pi}\left( \nabla \frac{df}{dD_s} \right).
\end{align*}
This completes the proof of the theorem.
\end{proof}

\begin{corollary}\label{cor:diffgrad}
Suppose that $\{D_s\}_{s\in (-\epsilon, \epsilon)}$ is a differentiable path of density operators such that
\begin{align*}
    D_0 = D,\quad \dot{D}_s=- \cL^*(D_s). 
\end{align*}
Then 
\begin{align*}
   \left \|\frac{d}{ds} D_s\right\|_{g_{\cL}} = \|\grad_{g, D_s} f\|_{D_s, \widehat{\Delta}},
\end{align*}
where $f(D)=H(D\| D_{\Delta})$.
\end{corollary}
\begin{proof}
By Theorem \ref{thm:kmsgradient} and Lemma \ref{lem:vecterfieldform}, we have that 
\begin{align*}
  \left \|\frac{d}{ds} D_s\right\|_{g_{\cL}}
  = & \min\left\{\|X\|_{D_s, \widehat{\Delta}}: \frac{d}{ds}D_s= \Div \K_{D_s,  \Pi} X.\right\}\\
  =& \left\| \nabla \frac{df}{dD_s} \right\|_{D_s, \widehat{\Delta}}\\
  =&\|\grad_{g, D} f\|_{D_s, \widehat{\Delta}}.
\end{align*}
This completes the proof of the corollary.
\end{proof}

In the following, we extend the intertwining property introduced by Carlen and Maas \cite{CarMaa17} for our setting.
\begin{definition}[Intertwining Property]
Suppose that $\{\Phi_t\}_{t\geq 0}$ is a bimodule quantum Markov semigroup bimodule KMS symmetric with respect to $\widehat{\Delta}$. 
We say that $\{\Phi_t\}_{t\geq 0}$ satisfies the intertwining property if a trace-preserving completely positive bimodule map $\widetilde{\Phi_t}^*: \cM_1 \to \cM_1 $ such that $\widetilde{\Phi_t}^*|_{\cM}=\Phi_t^*$, and
\begin{align*}
   \Phi_t^* \Div (\cdot)=\Div \bigoplus_{\ell=1}^m  (\widetilde{\Phi_t}^*\otimes I_\ell) (e^{-T^{(\ell)} t}  \cdot)  , \quad  t\geq 0,
\end{align*}    
where $I_\ell$ is the identity map on $M_{d_\ell}(\bC)$ and 
\begin{align*}
T^{(\ell)}=\left( T_{j,k}^{(\ell)} F_{j,k}^{(\ell)}\right)_{j,k=1}^{d_\ell}
\end{align*}
with $\left(T_{j,k}^{(\ell)}\right)$ a positive definite matrix in $M_{d_\ell}(\bC)$.
Here we identify $T^{(\ell)}$ with the linear map $X\mapsto \lambda^{-1}(\bE_{\cM_1}\otimes I)(T^{(\ell)}X)$.

If there is no confusion, we shall denote  $\displaystyle \bigoplus_{\ell=1}^m (\widetilde{\Phi_t}^*\otimes I_\ell)$ by $\widetilde{\Phi_t}^*$ for simplicity.
Let $\displaystyle T=\bigoplus_{\ell=1}^m T^{(\ell)}$ and $\beta$ the minimal eigenvalue of $T$.
\end{definition}


\begin{proposition}\label{prop:lowerbd}
Suppose that 
\begin{align*}
\K_{D, \Pi}^{1/2} T \K_{D, \Pi}^{-1/2} + \K_{D, \Pi}^{-1/2} T \K_{D, \Pi}^{1/2} \geq 2\beta >0.
\end{align*}
Then 
\begin{align*}
\left\|  \K_{D, \Pi}^{1/2}  e^{-Ts }\K_{D, \Pi}^{-1} e^{-T s} \K_{D, \Pi}^{1/2}  \right\| \leq e^{-2\beta s}.
\end{align*}
\end{proposition}
\begin{proof}
Firstly, we shows that $\K_{D, \Pi}^{1/2}  e^{-T s }\K_{D, \Pi}^{-1/2} =\exp (-\K_{D, \Pi}^{1/2}  T \K_{D, \Pi}^{-1/2} s)$.
For any $X\in \oM$, we have that 
\begin{align*}
& \K_{D, \Pi}^{1/2} T \K_{D, \Pi}^{-1/2}\K_{D, \Pi}^{1/2} (T(\K_{D, \Pi}^{-1/2} X)) \\
 =& \K_{D, \Pi}^{1/2} T (T(\K_{D, \Pi}^{-1/2} X)) \\
 =&  \lambda^{-1} \K_{D, \Pi}^{1/2} T( \bE_{\cM_1}\otimes I )(T \K_{D, \Pi}^{-1/2} X e_2) \\
  =&  \lambda^{-2} \K_{D, \Pi}^{1/2} (\bE_{\cM_1}\otimes I) (T( \bE_{\cM_1}\otimes I )(T \K_{D, \Pi}^{-1/2} X e_2)e_2) \\
  =&  \lambda^{-1} \K_{D, \Pi}^{1/2} (\bE_{\cM_1}\otimes I) (T^2 \K_{D, \Pi}^{-1/2} X e_2).
\end{align*}
This shows that $\K_{D, \Pi}^{1/2}  e^{-T s }\K_{D, \Pi}^{-1/2} =\exp (-\K_{D, \Pi}^{1/2}  T \K_{D, \Pi}^{-1/2} s)$.
Let $Y =\K_{D, \Pi}^{1/2}  T \K_{D, \Pi}^{-1/2}.$
By the assumption, we have that $Y+Y^*\geq \beta$.
By Lie-Trotter formula, for any self-adjoint operators $H,\,K$  we have that 
\begin{equation*}
    \|e^{H+iK}\|=\left\|  \lim_{n \rightarrow \infty}(e^{H/n})^n (e^{iK/n})^n \right\| \le \| e^H\|,
\end{equation*}
which implies that 
\begin{align*}
\|e^{-Ys} e^{-Y^* s}\| \leq \|e^{-(Y+Y^*)s}\|\leq e^{-2\beta s}.
\end{align*}
This completes the proof of the proposition.
\end{proof}

\begin{remark}
Suppose that the semigroup is bimodule GNS symmetric.
The matrix $T$ commutes with $\K_{D, \Pi}$ and $\beta$ is the minimal eigenvalue of $T$.
\end{remark}

\begin{remark}
By taking adjoint, we see that the intertwining property is equivalent to 
\begin{align*}
   \nabla \Phi_t(\cdot) =  e^{-T t} \widetilde{\Phi_t}\nabla(\cdot) ,
\end{align*}
Suppose that $\widetilde{\Phi}_t$ is a bimodule quantum Markov semigroup and $\cJ$ is its generator.
Then 
\begin{align*}
-\nabla \cL =-\cJ \nabla-T \nabla .
\end{align*}
By taking the corresponding derivatives, we obtain that 
\begin{align}\label{eq:derivations}
\partial_k^{(\ell)} \cL- \cJ \partial_k^{(\ell)} =\lambda^{-1}\sum_{j=1}^{d_\ell} T_{kj}^{(\ell)} \bE_{\cM_1} \left(F_{kj}^{(\ell)}[\cdot, e_1]e_2\right)
\end{align}
for all $k=1, \ldots, d_\ell$ and $\ell=1, \ldots, m$.
When $e^{-T t} \widetilde{\Phi}_t$ is a semigroup, the equality \eqref{eq:derivations} is equivalent to the intertwining property.
\end{remark}

\begin{proposition}
Suppose that $\Phi: \cM_1\to \cM_1$ is trace-preserving completely positive bimodule map.
Then for strictly positive $D\in \cM$, $\ell=1, \ldots, m$ and $\mathbf{X}^{(\ell)}\in \cM_1\otimes M_{d_\ell}(\bC)$, we have that 
\begin{align*}
\langle \K_{\Phi(D), \Pi}^{(\ell)-1} (\Phi\otimes I)(\mathbf{X}^{(\ell)}), (\Phi\otimes I )(\mathbf{X}^{(\ell)})\rangle \leq \langle \K_{D, \Pi}^{(\ell)-1}\mathbf{X}^{(\ell)}, \mathbf{X}^{(\ell)}\rangle. 
\end{align*}
\end{proposition}
\begin{proof}
By Lieb's concavity theorem, we obtain that the map 
\begin{align*}
D \mapsto &\int_0^1  (\tau_1\otimes \Tr_{d_\ell})\left( F_\ell^{1-2s} \mathbf{D}^{(\ell) s} \Pi^{(\ell)} \left( F_\ell^{1-2s}  \mathbf{X}^{(\ell)} \mathbf{D}^{(\ell) (1-s)} \right) \mathbf{X}^{(\ell)*} \right)  ds\\
=&\int_0^1 \lambda^{-1}(\tau_2\otimes \Tr_{d_\ell})\left( F_\ell^{1-2s} \mathbf{D}^{(\ell) s} \Pi^{(\ell)} F_\ell^{1-2s}  \mathbf{X}^{(\ell)} e_2 \mathbf{D}^{(\ell) (1-s)}  e_2\mathbf{X}^{(\ell)*} \right) ds (:= f(\mathbf{X}^{(\ell)} , \mathbf{D}^{(\ell)}) )
\end{align*}
is concave.
The variational representation (Legendre-Fenchel transform) for $\K_{D,\Pi}^{(\ell)-1}$:
\begin{equation*}
     \langle \K_{D,\Pi}^{(\ell) -1} X, X \rangle = \sup_{Y} \left( 2\Re\langle X, Y \rangle - \langle \K_{D,\Pi}^{(\ell)} Y, Y \rangle \right)
\end{equation*}
implies the map 
\begin{align*}
( \mathbf{X}, D) \mapsto  \langle \K_{D, \Pi}^{(\ell)-1}\mathbf{X}, \mathbf{X}\rangle
\end{align*}
is jointly convex.

The proof of the following follows a similar argument to that of the data processing inequality\cite{Lin75}.
By the Stinespring dilation theorem for finite von Neumann algebras, any trace-preserving completely positive map $\Phi$ can be decomposed as $\Phi(X) = \mathbb{E}(V X V^*)$, where $V: \mathcal{H}_S \to \mathcal{H}_{SE}$ is an isometry ($V^*V = I$), and $\mathbb{E}$ is a trace-preserving conditional expectation.
Applying the joint convexity of $f$ and Jensen's inequality:
\begin{align*}
    f(\Phi(X), \Phi(D)) &= f(\mathbb{E}(VXV^*), \mathbb{E}(VDV^*)) \\
    &\le \lim_{n\to\infty} \sum_{j} \alpha_j f(U_j VXV^* U_j^*, U_j VDV^* U_j^*)\\
    &=f(X,D).
\end{align*}
since the conditional expectation $\mathbb{E}(\cdot)$ can be represented as a limit of convex combinations of unitary twirlings (Dixmier approximation).

\end{proof}

Thanks to the intertwining property, we are able to give the following inequalities.
We shall investigate the concrete generalized logarithmic Sobolev inequality in a forthcoming paper.

\begin{theorem}[Bimodule Logarithmic Sobolev Inequality]\label{thm:entropydecay}
Suppose that $\{\Phi_t\}_{t\geq 0}$ is a bimodule quantum Markov semigroup which is relatively ergodic and bimodule KMS symmetric with respect to $\widehat{\Delta}$ and satisfies the intertwining property and
\begin{align*}
\K_{D, \Pi}^{1/2} T \K_{D, \Pi}^{-1/2} + \K_{D, \Pi}^{-1/2} T \K_{D, \Pi}^{1/2} \geq \beta.
\end{align*}    
Then for any strictly positive density operator $D\in \cM$, we have that 
\begin{equation}\label{eq:entropy1}
\begin{aligned}
 & H(\Phi_t^*(D)\|D_\Delta) -H(\bE_{\cN}(D)\overline{F_\Phi} \|D_{\Delta}) \\
 \leq &  e^{-2\beta t}\left( H(D\| D_{\Delta})  - H(\bE_{\cN}(D)\overline{F_\Phi} \|D_{\Delta})\right). 
\end{aligned}
\end{equation}
Furthermore, we have that 
\begin{align}\label{eq:logsob}
H(D\|D_{\Delta})- H(\bE_{\cN}(D)\overline{F_\Phi} \|D_{\Delta}) \leq \frac{1}{2\beta} \tau(\cL^*(D)( \log D-\log D_{\Delta})).    
\end{align}
This is called the generalized logarithmic Sobolev inequality.
\end{theorem}
\begin{proof}
Suppose that $X(s)$ is the solution of $\displaystyle \frac{d}{ds} D_s=\Div \K_{D, \Pi} X(s)$ minimizing the norm $\|\cdot\|_{D, \widehat{\Delta}}$, where $\{D_s\}_s$ is a differentiable path in $\cM_+(D)$ such that $D_0=D$. 
We have that 
\begin{align*}
    \frac{d}{ds} \Phi_t^*(D_s) 
    = &   \Phi_t^*\Div \K_{D, \Pi} X(s) \\
    = &  \Div \widetilde{\Phi_t}^*( e^{-T t} \K_{D,\Pi} X(s) ) \\
    = & \Div \K_{\Phi_t^*(D), \Pi} \left( \K_{\Phi_t^*(D), \Pi}^{-1} \widetilde{\Phi_t}^*(e^{-T t} \K_{D ,\Pi} X(s))  \right) .
\end{align*}
By Proposition \ref{prop:lowerbd}, we have that 
\begin{align*}
 \left\|\frac{d}{ds}{\Phi_t}^*(D_s)\right\|_{g_\cL} ^2
\leq &  \left\| \K_{{\Phi_t}^*(D), \Pi}^{-1} \widetilde{\Phi_t}^* (e^{-T t}\K_{D ,\Pi} X(s))  \right\|_{{\Phi_t}^*(D), \widehat{\Delta}}^2 \\
\leq &   \left\langle  \K_{{\Phi_t}^*(D), \Pi}^{-1}\widetilde{\Phi_t}^* ( e^{-T t} \K_{D, \Pi} X(s)  ), \widetilde{\Phi_t}^*( e^{-T t} \K_{D, \Pi} X(s) )\right\rangle \\
\leq &   \left\langle   \K_{D, \Pi}^{-1} e^{-T t}  \K_{D, \Pi} X(s) ,  e^{-T t}  \K_{D,  \Pi} X(s)   \right\rangle\\
\leq &   e^{-2\beta t} \left\langle  X(s)  , \K_{D,  \Pi} X(s) \right\rangle   \\
= & e^{-2\beta t} \left\|\frac{d}{ds} D_s\right\|_{g_\cL}^2.
\end{align*}
This implies that $\displaystyle \lim_{t\to \infty}  \left\|\frac{d}{ds}{\Phi_t}^*(D_s)\right\|_{g_\cL}=0$ and 
\begin{align*}
   \frac{d}{dt}  \left\|\frac{d}{ds}{\Phi_t}^*(D_s)\right\|_{g_{\cL}} ^2  
 \leq -2\beta \left\|\frac{d}{ds}\Phi_t^*( D_s)\right\|_{g_{\cL}}^2.
\end{align*}
By Corollary \ref{cor:diffgrad}, we have that 
\begin{align*}
\|\grad_{g, \Phi_t^*(D_s)} f\|_{\Phi_t^*(D_s), \widehat{\Delta}}^2 
\leq e^{-2\beta t} \|\grad_{g, D_s} f\|_{D_s, \widehat{\Delta}}^2. 
\end{align*}
Hence 
\begin{align*}
  \left. \frac{d}{dt}  \|\grad_{g, \Phi_t^*(D_s)} f\|_{D, \widehat{\Delta}}^2\right|_{t=0^+} \leq -2\beta \|\grad_{g, D_s} f\|_{D_s, \widehat{\Delta}}^2.
\end{align*}
Note that 
\begin{align*}
    \left.\frac{d}{dt} f(\Phi_t^*(D))\right|_{t=0}
    =& -\tau\left(\frac{df}{d D}\cL^*(D)\right)\\
    =& -  \tau\left(\frac{df}{d D} \Div\K_{D, \Pi} \left( \nabla \frac{df}{dD}\right)\right) \\
    =& - \left\langle \K_{D,\Pi} \left( \nabla \frac{df}{dD}\right), \left(\nabla\frac{df}{d D}\right) \right\rangle \\
    =& -\|\grad_{g, D} f\|_{D, \widehat{\Delta}}^2.
\end{align*}
We see that 
\begin{align}\label{eq:energe}
 \frac{d}{dt} H(\Phi_t^*(D_s)\|D_{\Delta}) =-\left\|\grad_{g, \Phi_t^*(D_s)} f\right\|_{\Phi_t^*(D_s), \widehat{\Delta}}^2.
\end{align}
Hence
\begin{align*}
    \frac{d}{dt} H(\Phi_t^*(D_s)\|D_{\Delta})
    = & -\left\|\frac{d}{ds}\Phi_t^*(D_s)\right\|_{g_\cL} ^2\\
    =& \int_t^\infty \frac{d}{dr}\left\|\frac{d}{ds}{\Phi_r}^*(D_s)\right\|_{g_\cL} ^2 dr\\
    \leq & -2\beta \int_t^\infty \|\grad_{g, {\Phi_r}^*(D_s)} f\|_{{\Phi_r}^*(D_s), \widehat{\Delta}}^2 dr\\
    = & 2\beta H(\bE_{\Phi}^*(D_s)\|D_{\Delta}) -2\beta  H(\Phi_t^*(D_s)\|D_{\Delta}),
\end{align*}
i.e.
\begin{align}\label{eq:entropy0}
    \frac{d}{dt} H(\Phi_t^*(D_s)\|D_{\Delta}) \leq 2\beta H(\bE_{\Phi}^*(D_s)\|D_{\Delta}) -2\beta H(\Phi_t^*(D_s)\| D_{\Delta}).
\end{align}
By integrating Equation \eqref{eq:entropy0} with respect to $t$, we obtain that
\begin{align*}
 H(\Phi_t^*(D_s)\|D_\Delta)-H(\bE_{\Phi}^*(D_s)\|D_{\Delta}) \leq e^{-2\beta t} (H(D_s\| D_{\Delta}) - H(\bE_{\Phi}^*(D_s)\|D_{\Delta})).  
\end{align*}
By Equation \eqref{eq:energe}, we obtain that 
\begin{align*}
 H(D_s\| D_{\Delta}) - H(\bE_{\Phi}^*(D_s)\|D_{\Delta})\leq \frac{1}{2\beta} \|\grad_{g, D} f\|_{D, \widehat{\Delta}}^2. 
\end{align*}
Note that 
\begin{align*}
   H(\bE_{\Phi}^*(D_s)\|D_{\Delta}) =H(\bE_{\cN}(D_s)\overline{F_\Phi} \|D_{\Delta}).
\end{align*}
We see the first equality of the theorem is true.
By Lemma \ref{lem:vecterfieldform}, we have that 
\begin{align*}
    \left\|\grad_{g, D} f \right\|_{D, \widehat{\Delta}}^2 = \tau(\cL^*(D)(\log D-\log D_{\Delta})).
\end{align*}
Therefore, the second inequality of the theorem is true.

\end{proof}

By using the generalized logarithmic Sobolev inequality, we have the following bimodule Talagrand inequality.
\begin{theorem}[Bimodule Talagrand Inequality]\label{thm:talagrand}
Suppose that $\{\Phi_t\}_{t\geq 0}$ is a bimodule quantum Markov semigroup bimodule KMS symmetry with respect to $\widehat{\Delta}$ and relatively ergodic. 
Suppose that Equation \eqref{eq:logsob} holds. 
We have that 
\begin{align*}
    d(D, D_\Delta) \leq  \sqrt{ \frac{2 H(D\| D_\Delta)- 2 H(\bE_{\Phi}^*(D) \|D_{\Delta})}{\beta}}.
\end{align*}
\end{theorem}
\begin{proof}
Suppose that $D\in \cM$ is a density operator.
Let $\Delta_D=\bE_{\cN}(D)\overline{F_\Phi}$.
Then the distance between $D$ and $\Delta_D$ is described by the following equation:
\begin{align*}
d(D, \Delta_D)= & \int_0^\infty \left\|\frac{d}{dt} {\Phi_t}^*(D)\right\|_{g_\cL} dt .
\end{align*}
Note that $\displaystyle \left\|\frac{d}{dt} {\Phi_t}^*(D)\right\|_{g_\cL}^2 =- \frac{d}{dt} H({\Phi_t}^*(D)\| D_{\Delta})$.
By the Cauchy-Schwarz inequality, for any $0< t_1<  t_2< \infty$, we have that 
\begin{align*}
    \int_{t_1}^{t_2} \left\|\frac{d}{dt} {\Phi_t}^*(D)\right\|_{g_\cL} dt
    \leq  \sqrt{t_2-t_1} \sqrt{H({\Phi}_{t_1}^*(D)\| D_{\Delta})-H({\Phi}_{t_2}^*(D)\|D_\Delta)}
\end{align*}
Fix $\epsilon>0$, for each $k\in \bN$, we take $t_k\in \bR$ such that 
\begin{align*}
   &  H({\Phi}_{t_k}^*(D)\| D_{\Delta})-  H(\bE_{\cN}(D)\overline{F_\Phi} \|D_{\Delta})\\
    =&  e^{-k \epsilon} (H(D\| D_{\Delta})-  H(\bE_{\cN}(D)\overline{F_\Phi} \|D_{\Delta})).
\end{align*}
By Equation \eqref{eq:entropy1}, we have that $\displaystyle t_k-t_{k-1}\leq \frac{\epsilon}{2\beta}$.
By Equation \eqref{eq:entropy1} again, we have that
\begin{equation}\label{eq:taga}
\begin{aligned}
   &  \int_{t_{k-1}}^{t_k} \left\|\frac{d}{dt} \Phi_t^*(D)\right\|_{g_{\cL}}dt \\
     \leq &   \sqrt{t_k-t_{k-1}} \sqrt{H(\Phi_{t_{k-1}}^*(D)\| D_{\Delta})-H(\Phi_{t_k}^*(D)\|D_\Delta)} \\
    \leq & \sqrt{\frac{\epsilon}{2\beta} (e^{-(k-1)\epsilon}-e^{-k\epsilon}) (H(D\| D_{\Delta})-  H(\bE_{\cN}(D)\overline{F_\Phi} \|D_{\Delta}))}\\
    =& e^{-k\epsilon /2}\sqrt{\epsilon (e^\epsilon-1)} \sqrt{\frac{H(D\| D_{\Delta})-  H(\bE_{\cN}(D)\overline{F_\Phi} \|D_{\Delta})}{2\beta}} .
\end{aligned}
\end{equation}
Note that 
\begin{align*}
 \lim_{\epsilon\to 0} \sum_{k=1}^\infty e^{-k\epsilon /2}\sqrt{\epsilon (e^\epsilon-1)} =2.
\end{align*}
By taking summation of Equation \eqref{eq:taga} with respect to $k$, we see that the theorem is true.

\end{proof}

\section{Clifford Algebras}
In this section, we shall check the intertwining property for a KMS symmetric semigroup in Clifford algebra which might be known to the experts and obtain the generalized logarithmic Sobolev inequality.

Let $\mathfrak{C}$ be the Clifford algebra on two unitary generators $Q, P$, i.e.
\begin{align*}
    Q^2=P^2=1, \quad PQ+QP=0.
\end{align*}
The index $\lambda^{-1}=4$.
Let
\begin{align*}
z=\frac{1}{\sqrt{2}}(Q+iP), \quad \displaystyle z^*=\frac{1}{\sqrt{2}}(Q-iP), \quad w=iQP, \quad v=wz.
\end{align*}
We have that $v^*v+vv^*=z^*z+zz^*=2$.
Suppose that 
\begin{align*}
H=a
\vcenter{\hbox{\begin{tikzpicture}[scale=0.65]
    \begin{scope}[shift={(0,1.5)}]
    \draw [blue] (-0.5, 0.8)--(-0.5, 0) .. controls +(0, -0.6) and +(0,-0.6).. (0.5, 0)--(0.5, 0.8);    
\begin{scope}[shift={(0.5, 0.3)}]
\draw [fill=white] (-0.3, -0.3) rectangle (0.3, 0.3);
\node at (0, 0) {\tiny $P$};
\end{scope}
    \end{scope}
\draw [blue] (-0.5, -0.8)--(-0.5, 0) .. controls +(0, 0.6) and +(0,0.6).. (0.5, 0)--(0.5, -0.8);
\begin{scope}[shift={(0.5, -0.3)}]
\draw [fill=white] (-0.3, -0.3) rectangle (0.3, 0.3);
\node at (0, 0) {\tiny $P$};
\end{scope}
\end{tikzpicture}}}
+
b
\vcenter{\hbox{\begin{tikzpicture}[scale=0.65]
    \begin{scope}[shift={(0,1.5)}]
    \draw [blue] (-0.5, 0.8)--(-0.5, 0) .. controls +(0, -0.6) and +(0,-0.6).. (0.5, 0)--(0.5, 0.8);    
\begin{scope}[shift={(0.5, 0.3)}]
\draw [fill=white] (-0.3, -0.3) rectangle (0.3, 0.3);
\node at (0, 0) {\tiny $Q$};
\end{scope}
    \end{scope}
\draw [blue] (-0.5, -0.8)--(-0.5, 0) .. controls +(0, 0.6) and +(0,0.6).. (0.5, 0)--(0.5, -0.8);
\begin{scope}[shift={(0.5, -0.3)}]
\draw [fill=white] (-0.3, -0.3) rectangle (0.3, 0.3);
\node at (0, 0) {\tiny $Q$};
\end{scope}
\end{tikzpicture}}},
\end{align*}
where $a, b>0$ and $a+b=1$, 
and the bimodule modular operator $\widehat{\Delta}$ is given by 
\begin{align*}
    \widehat{\Delta}=e_2+\mu^{-1} \lambda^{1/2}
\vcenter{\hbox{\begin{tikzpicture}[scale=0.65]
    \begin{scope}[shift={(0,1.5)}]
    \draw [blue] (-0.5, 0.8)--(-0.5, 0) .. controls +(0, -0.6) and +(0,-0.6).. (0.5, 0)--(0.5, 0.8);    
\begin{scope}[shift={(0.5, 0.3)}]
\draw [fill=white] (-0.3, -0.3) rectangle (0.3, 0.3);
\node at (0, 0) {\tiny $v$};
\end{scope}
    \end{scope}
\draw [blue] (-0.5, -0.8)--(-0.5, 0) .. controls +(0, 0.6) and +(0,0.6).. (0.5, 0)--(0.5, -0.8);
\begin{scope}[shift={(0.5, -0.3)}]
\draw [fill=white] (-0.3, -0.3) rectangle (0.3, 0.3);
\node at (0, 0) {\tiny $v^*$};
\end{scope}
\end{tikzpicture}}}
+\mu\lambda^{1/2}
\vcenter{\hbox{\begin{tikzpicture}[scale=0.65]
    \begin{scope}[shift={(0,1.5)}]
    \draw [blue] (-0.5, 0.8)--(-0.5, 0) .. controls +(0, -0.6) and +(0,-0.6).. (0.5, 0)--(0.5, 0.8);    
\begin{scope}[shift={(0.5, 0.3)}]
\draw [fill=white] (-0.3, -0.3) rectangle (0.3, 0.3);
\node at (0, 0) {\tiny $v^*$};
\end{scope}
    \end{scope}
\draw [blue] (-0.5, -0.8)--(-0.5, 0) .. controls +(0, 0.6) and +(0,0.6).. (0.5, 0)--(0.5, -0.8);
\begin{scope}[shift={(0.5, -0.3)}]
\draw [fill=white] (-0.3, -0.3) rectangle (0.3, 0.3);
\node at (0, 0) {\tiny $v$};
\end{scope}
\end{tikzpicture}}}
+\lambda^{1/2}
\vcenter{\hbox{\begin{tikzpicture}[scale=0.65]
    \begin{scope}[shift={(0,1.5)}]
    \draw [blue] (-0.5, 0.8)--(-0.5, 0) .. controls +(0, -0.6) and +(0,-0.6).. (0.5, 0)--(0.5, 0.8);    
\begin{scope}[shift={(0.5, 0.3)}]
\draw [fill=white] (-0.4, -0.3) rectangle (0.4, 0.3);
\node at (0, 0) {\tiny $iPQ$};
\end{scope}
    \end{scope}
\draw [blue] (-0.5, -0.8)--(-0.5, 0) .. controls +(0, 0.6) and +(0,0.6).. (0.5, 0)--(0.5, -0.8);
\begin{scope}[shift={(0.5, -0.3)}]
\draw [fill=white] (-0.4, -0.3) rectangle (0.4, 0.3);
\node at (0, 0) {\tiny $iPQ$};
\end{scope}
\end{tikzpicture}}}.
\end{align*}
We have that $\overline{H}=H$, $\overline{\widehat{\Delta}}=\widehat{\Delta}^{-1}$,
\begin{align*}
    \widehat{\cL}_0=\widehat{\Delta}^{1/2} H \widehat{\Delta}^{1/2}
    =\frac{1}{2}\mu^{-1}
    \vcenter{\hbox{\begin{tikzpicture}[scale=0.65]
    \begin{scope}[shift={(0,1.5)}]
    \draw [blue] (-0.5, 0.8)--(-0.5, 0) .. controls +(0, -0.6) and +(0,-0.6).. (0.5, 0)--(0.5, 0.8);    
\begin{scope}[shift={(0.5, 0.3)}]
\draw [fill=white] (-0.3, -0.3) rectangle (0.3, 0.3);
\node at (0, 0) {\tiny $v$};
\end{scope}
    \end{scope}
\draw [blue] (-0.5, -0.8)--(-0.5, 0) .. controls +(0, 0.6) and +(0,0.6).. (0.5, 0)--(0.5, -0.8);
\begin{scope}[shift={(0.5, -0.3)}]
\draw [fill=white] (-0.3, -0.3) rectangle (0.3, 0.3);
\node at (0, 0) {\tiny $v^*$};
\end{scope}
\end{tikzpicture}}}
+\frac{1}{2}\mu\vcenter{\hbox{\begin{tikzpicture}[scale=0.65]
    \begin{scope}[shift={(0,1.5)}]
    \draw [blue] (-0.5, 0.8)--(-0.5, 0) .. controls +(0, -0.6) and +(0,-0.6).. (0.5, 0)--(0.5, 0.8);    
\begin{scope}[shift={(0.5, 0.3)}]
\draw [fill=white] (-0.3, -0.3) rectangle (0.3, 0.3);
\node at (0, 0) {\tiny $v^*$};
\end{scope}
    \end{scope}
\draw [blue] (-0.5, -0.8)--(-0.5, 0) .. controls +(0, 0.6) and +(0,0.6).. (0.5, 0)--(0.5, -0.8);
\begin{scope}[shift={(0.5, -0.3)}]
\draw [fill=white] (-0.3, -0.3) rectangle (0.3, 0.3);
\node at (0, 0) {\tiny $v$};
\end{scope}
\end{tikzpicture}}}
+\frac{1}{2}(b-a) \vcenter{\hbox{\begin{tikzpicture}[scale=0.65]
    \begin{scope}[shift={(0,1.5)}]
    \draw [blue] (-0.5, 0.8)--(-0.5, 0) .. controls +(0, -0.6) and +(0,-0.6).. (0.5, 0)--(0.5, 0.8);    
\begin{scope}[shift={(0.5, 0.3)}]
\draw [fill=white] (-0.3, -0.3) rectangle (0.3, 0.3);
\node at (0, 0) {\tiny $v$};
\end{scope}
    \end{scope}
\draw [blue] (-0.5, -0.8)--(-0.5, 0) .. controls +(0, 0.6) and +(0,0.6).. (0.5, 0)--(0.5, -0.8);
\begin{scope}[shift={(0.5, -0.3)}]
\draw [fill=white] (-0.3, -0.3) rectangle (0.3, 0.3);
\node at (0, 0) {\tiny $v$};
\end{scope}
\end{tikzpicture}}}
+
\frac{1}{2}(b-a) \vcenter{\hbox{\begin{tikzpicture}[scale=0.65]
    \begin{scope}[shift={(0,1.5)}]
    \draw [blue] (-0.5, 0.8)--(-0.5, 0) .. controls +(0, -0.6) and +(0,-0.6).. (0.5, 0)--(0.5, 0.8);    
\begin{scope}[shift={(0.5, 0.3)}]
\draw [fill=white] (-0.3, -0.3) rectangle (0.3, 0.3);
\node at (0, 0) {\tiny $v^*$};
\end{scope}
    \end{scope}
\draw [blue] (-0.5, -0.8)--(-0.5, 0) .. controls +(0, 0.6) and +(0,0.6).. (0.5, 0)--(0.5, -0.8);
\begin{scope}[shift={(0.5, -0.3)}]
\draw [fill=white] (-0.3, -0.3) rectangle (0.3, 0.3);
\node at (0, 0) {\tiny $v^*$};
\end{scope}
\end{tikzpicture}}}
\end{align*}
and $\cR(H)=\cR(\widehat{\cL}_0)$.
Now we see that 
\begin{align*}
    1*(\widehat{\Delta}^{1/2} H \widehat{\Delta}^{1/2}) =\frac{1}{2}\mu^{-1} v^*v+\frac{1}{2}\mu vv^*(=\mathbf{y}).
\end{align*}
The semigroup $\{\Phi_t\}_{t\geq 0}$ arising from $\cL$ is defined by
\begin{align*}
    \cL(x)=\frac{1}{2}\mathbf{y} x +\frac{1}{2} x\mathbf{y} -x*\widehat{\cL}_0.
\end{align*}
\begin{lemma}
The semigroup described above is KMS symmetric with $\cL_w=0$.
Let $\displaystyle \rho=\frac{2 \bfy}{\mu+\mu^{-1}}$.
We have that $\widehat{\Delta}=\widehat{\Delta}_{\rho, 1/2}$.
\end{lemma}
\begin{proof}
By Proposition \ref{prop:kms1}, it suffices to check the following condition
\begin{align}\label{eqn: KMS symmetry condition}
 \vcenter{\hbox{\begin{tikzpicture}[scale=0.65]
\draw [blue] (-0.5, -1.5)--(-0.5, 0) .. controls +(0, 0.6) and +(0,0.6).. (0.5, 0)--(0.5, -1.5);
\begin{scope}[shift={(0.5, -0.3)}]
\draw [fill=white] (-0.3, -0.3) rectangle (0.3, 0.3);
\node at (0, 0) {\tiny $\mathbf{y}$};
\end{scope}
\begin{scope}[shift={(0, -1)}]
\draw [fill=white] (-0.7, -0.3) rectangle (0.7, 0.3);
\node at (0, 0) {\tiny $\overline{\widehat{\Delta}}$};    
\end{scope}
\end{tikzpicture}}}
= \vcenter{\hbox{\begin{tikzpicture}[scale=0.65]
\draw [blue] (-0.5, -0.8)--(-0.5, 0) .. controls +(0, 0.6) and +(0,0.6).. (0.5, 0)--(0.5, -0.8);
\begin{scope}[shift={(0.5, -0.3)}]
\draw [fill=white] (-0.3, -0.3) rectangle (0.3, 0.3);
\node at (0, 0) {\tiny $\mathbf{y}$};
\end{scope}
\end{tikzpicture}}}.  
\end{align}
Note that
\begin{align*}
    \tau(\bfy)=\frac{1}{2}(\mu^{-1}+\mu), \quad \tau(v^*\bfy)=\tau(v\bfy)=0,\quad
    \tau(iPQ \bfy)=\frac{1}{2}(\mu-\mu^{-1}).
\end{align*}
Then we have that 
\begin{align*}
 \vcenter{\hbox{\begin{tikzpicture}[scale=0.65]
\draw [blue] (-0.5, -1.5)--(-0.5, 0) .. controls +(0, 0.6) and +(0,0.6).. (0.5, 0)--(0.5, -1.5);
\begin{scope}[shift={(0.5, -0.3)}]
\draw [fill=white] (-0.3, -0.3) rectangle (0.3, 0.3);
\node at (0, 0) {\tiny $\mathbf{y}$};
\end{scope}
\begin{scope}[shift={(0, -1)}]
\draw [fill=white] (-0.7, -0.3) rectangle (0.7, 0.3);
\node at (0, 0) {\tiny $\overline{\widehat{\Delta}}$};    
\end{scope}
\end{tikzpicture}}}&=\frac{1}{2}\mu^{-1}
    \left(
\vcenter{\hbox{\begin{tikzpicture}[scale=0.65]
\draw [blue] (-0.5, -0.8)--(-0.5, 0) .. controls +(0, 0.6) and +(0,0.6).. (0.5, 0)--(0.5, -0.8);
\end{tikzpicture}}}
+\vcenter{\hbox{\begin{tikzpicture}[scale=0.65]
\draw [blue] (-0.5, -0.8)--(-0.5, 0) .. controls +(0, 0.6) and +(0,0.6).. (0.5, 0)--(0.5, -0.8);
\begin{scope}[shift={(0.5, -0.3)}]
\draw [fill=white] (-0.4, -0.3) rectangle (0.4, 0.3);
\node at (0, 0) {\tiny $iPQ$};
\end{scope}
\end{tikzpicture}}} \right)
+\frac{1}{2}\mu
\left(
\vcenter{\hbox{\begin{tikzpicture}[scale=0.65]
\draw [blue] (-0.5, -0.8)--(-0.5, 0) .. controls +(0, 0.6) and +(0,0.6).. (0.5, 0)--(0.5, -0.8);
\end{tikzpicture}}}
-\vcenter{\hbox{\begin{tikzpicture}[scale=0.65]
\draw [blue] (-0.5, -0.8)--(-0.5, 0) .. controls +(0, 0.6) and +(0,0.6).. (0.5, 0)--(0.5, -0.8);
\begin{scope}[shift={(0.5, -0.3)}]
\draw [fill=white] (-0.4, -0.3) rectangle (0.4, 0.3);
\node at (0, 0) {\tiny $iPQ$};
\end{scope}
\end{tikzpicture}}} \right)\\
&=\frac{1}{2}\mu^{-1} \,
\vcenter{\hbox{\begin{tikzpicture}[scale=0.65]
\draw [blue] (-0.5, -0.8)--(-0.5, 0) .. controls +(0, 0.6) and +(0,0.6).. (0.5, 0)--(0.5, -0.8);
\begin{scope}[shift={(0.5, -0.3)}]
\draw [fill=white] (-0.4, -0.3) rectangle (0.4, 0.3);
\node at (0, 0) {\tiny $v^*v$};
\end{scope}\end{tikzpicture}}} +\frac{1}{2} \mu\,
\vcenter{\hbox{\begin{tikzpicture}[scale=0.65]
\draw [blue] (-0.5, -0.8)--(-0.5, 0) .. controls +(0, 0.6) and +(0,0.6).. (0.5, 0)--(0.5, -0.8);
\begin{scope}[shift={(0.5, -0.3)}]
\draw [fill=white] (-0.4, -0.3) rectangle (0.4, 0.3);
\node at (0, 0) {\tiny $vv^*$};
\end{scope}\end{tikzpicture}}}
=\vcenter{\hbox{\begin{tikzpicture}[scale=0.65]
\draw [blue] (-0.5, -0.8)--(-0.5, 0) .. controls +(0, 0.6) and +(0,0.6).. (0.5, 0)--(0.5, -0.8);
\begin{scope}[shift={(0.5, -0.3)}]
\draw [fill=white] (-0.3, -0.3) rectangle (0.3, 0.3);
\node at (0, 0) {\tiny $\mathbf{y}$};
\end{scope}
\end{tikzpicture}}}.
\end{align*}
Thus $\{\Phi_t\}_{t\geq 0}$ is KMS symmetric.

Note that 
\begin{align*}
\frac{1}{4}\vcenter{\hbox{\begin{tikzpicture}[scale=0.8]
    \draw [blue] (-1, 1.2)--(-1, -1.2) (-0.3, 1.2)--(-0.3, -1.2);
\begin{scope}[shift={(-1, 0)}]
\draw [fill=white] (-0.3, -0.3) rectangle (0.3, 0.3);
\node at (0, 0) {\tiny $\overline{v^*v}$};
\end{scope}
    \begin{scope}[shift={(-0.3, 0)}]
\draw [fill=white] (-0.3, -0.3) rectangle (0.3, 0.3);
\node at (0, 0) {\tiny $v^*v$};
\end{scope}
\end{tikzpicture}}} 
=\frac{1}{2}\lambda^{1/2}
    \vcenter{\hbox{\begin{tikzpicture}[scale=0.65]
    \begin{scope}[shift={(0,1.5)}]
    \draw [blue] (-0.5, 0.8)--(-0.5, 0) .. controls +(0, -0.6) and +(0,-0.6).. (0.5, 0)--(0.5, 0.8);    
    \end{scope}
\draw [blue] (-0.5, -0.8)--(-0.5, 0) .. controls +(0, 0.6) and +(0,0.6).. (0.5, 0)--(0.5, -0.8);
\end{tikzpicture}}}
+\frac{1}{2}\lambda^{1/2} \vcenter{\hbox{\begin{tikzpicture}[scale=0.65]
    \begin{scope}[shift={(0,1.5)}]
    \draw [blue] (-0.5, 0.8)--(-0.5, 0) .. controls +(0, -0.6) and +(0,-0.6).. (0.5, 0)--(0.5, 0.8);    
\begin{scope}[shift={(0.5, 0.3)}]
\draw [fill=white] (-0.3, -0.3) rectangle (0.3, 0.3);
\node at (0, 0) {\tiny $iPQ$};
\end{scope}
    \end{scope}
\draw [blue] (-0.5, -0.8)--(-0.5, 0) .. controls +(0, 0.6) and +(0,0.6).. (0.5, 0)--(0.5, -0.8);
\end{tikzpicture}}}
+\frac{1}{2}\lambda^{1/2} \vcenter{\hbox{\begin{tikzpicture}[scale=0.65]
    \begin{scope}[shift={(0,1.5)}]
    \draw [blue] (-0.5, 0.8)--(-0.5, 0) .. controls +(0, -0.6) and +(0,-0.6).. (0.5, 0)--(0.5, 0.8);    
    \end{scope}
\draw [blue] (-0.5, -0.8)--(-0.5, 0) .. controls +(0, 0.6) and +(0,0.6).. (0.5, 0)--(0.5, -0.8);
\begin{scope}[shift={(0.5, -0.3)}]
\draw [fill=white] (-0.3, -0.3) rectangle (0.3, 0.3);
\node at (0, 0) {\tiny $iPQ$};
\end{scope}
\end{tikzpicture}}} 
+
\frac{1}{2}\lambda^{1/2}\vcenter{\hbox{\begin{tikzpicture}[scale=0.65]
    \begin{scope}[shift={(0,1.5)}]
    \draw [blue] (-0.5, 0.8)--(-0.5, 0) .. controls +(0, -0.6) and +(0,-0.6).. (0.5, 0)--(0.5, 0.8);    
\begin{scope}[shift={(0.5, 0.3)}]
\draw [fill=white] (-0.4, -0.3) rectangle (0.4, 0.3);
\node at (0, 0) {\tiny $iPQ$};
\end{scope}
    \end{scope}
\draw [blue] (-0.5, -0.8)--(-0.5, 0) .. controls +(0, 0.6) and +(0,0.6).. (0.5, 0)--(0.5, -0.8);
\begin{scope}[shift={(0.5, -0.3)}]
\draw [fill=white] (-0.4, -0.3) rectangle (0.4, 0.3);
\node at (0, 0) {\tiny $iPQ$};
\end{scope}
\end{tikzpicture}}}, \\
\frac{1}{4}\vcenter{\hbox{\begin{tikzpicture}[scale=0.8]
    \draw [blue] (-1, 1.2)--(-1, -1.2) (-0.3, 1.2)--(-0.3, -1.2);
\begin{scope}[shift={(-1, 0)}]
\draw [fill=white] (-0.3, -0.3) rectangle (0.3, 0.3);
\node at (0, 0) {\tiny $\overline{vv^*}$};
\end{scope}
    \begin{scope}[shift={(-0.3, 0)}]
\draw [fill=white] (-0.3, -0.3) rectangle (0.3, 0.3);
\node at (0, 0) {\tiny $v^*v$};
\end{scope}
\end{tikzpicture}}} 
=
\lambda^{1/2} \vcenter{\hbox{\begin{tikzpicture}[scale=0.65]
    \begin{scope}[shift={(0,1.5)}]
    \draw [blue] (-0.5, 0.8)--(-0.5, 0) .. controls +(0, -0.6) and +(0,-0.6).. (0.5, 0)--(0.5, 0.8);    
\begin{scope}[shift={(0.5, 0.3)}]
\draw [fill=white] (-0.3, -0.3) rectangle (0.3, 0.3);
\node at (0, 0) {\tiny $v$};
\end{scope}
    \end{scope}
\draw [blue] (-0.5, -0.8)--(-0.5, 0) .. controls +(0, 0.6) and +(0,0.6).. (0.5, 0)--(0.5, -0.8);
\begin{scope}[shift={(0.5, -0.3)}]
\draw [fill=white] (-0.3, -0.3) rectangle (0.3, 0.3);
\node at (0, 0) {\tiny $v^*$};
\end{scope}
\end{tikzpicture}}},
\quad 
\frac{1}{4}\vcenter{\hbox{\begin{tikzpicture}[scale=0.8]
    \draw [blue] (-1, 1.2)--(-1, -1.2) (-0.3, 1.2)--(-0.3, -1.2);
\begin{scope}[shift={(-1, 0)}]
\draw [fill=white] (-0.3, -0.3) rectangle (0.3, 0.3);
\node at (0, 0) {\tiny $\overline{v^*v}$};
\end{scope}
    \begin{scope}[shift={(-0.3, 0)}]
\draw [fill=white] (-0.3, -0.3) rectangle (0.3, 0.3);
\node at (0, 0) {\tiny $vv^*$};
\end{scope}
\end{tikzpicture}}} 
=
\lambda^{1/2} \vcenter{\hbox{\begin{tikzpicture}[scale=0.65]
    \begin{scope}[shift={(0,1.5)}]
    \draw [blue] (-0.5, 0.8)--(-0.5, 0) .. controls +(0, -0.6) and +(0,-0.6).. (0.5, 0)--(0.5, 0.8);    
\begin{scope}[shift={(0.5, 0.3)}]
\draw [fill=white] (-0.3, -0.3) rectangle (0.3, 0.3);
\node at (0, 0) {\tiny $v^*$};
\end{scope}
    \end{scope}
\draw [blue] (-0.5, -0.8)--(-0.5, 0) .. controls +(0, 0.6) and +(0,0.6).. (0.5, 0)--(0.5, -0.8);
\begin{scope}[shift={(0.5, -0.3)}]
\draw [fill=white] (-0.3, -0.3) rectangle (0.3, 0.3);
\node at (0, 0) {\tiny $v$};
\end{scope}
\end{tikzpicture}}}.
\end{align*}
We obtain that $\widehat{\Delta}_{\rho, 1/2}=\widehat{\Delta}$ and $\widehat{\Delta}$ is implemented by the faithful state $\rho$.

\end{proof}

In the following, we shall check the intertwining property.
The directional derivations $\partial_1, \partial_2$ is given by
\begin{align*}
    \vcenter{\hbox{\begin{tikzpicture}[scale=0.65]
    \begin{scope}[shift={(0,1.5)}]
    \draw [blue] (-0.5, 0.8)--(-0.5, 0) .. controls +(0, -0.6) and +(0,-0.6).. (0.5, 0)--(0.5, 0.8);    
\begin{scope}[shift={(0.5, 0.3)}]
\draw [fill=white] (-0.3, -0.3) rectangle (0.3, 0.3);
\node at (0, 0) {\tiny $v$};
\end{scope}
    \end{scope}
\draw [blue] (-0.5, -0.8)--(-0.5, 0) .. controls +(0, 0.6) and +(0,0.6).. (0.5, 0)--(0.5, -0.8);
\begin{scope}[shift={(0.5, -0.3)}]
\draw [fill=white] (-0.3, -0.3) rectangle (0.3, 0.3);
\node at (0, 0) {\tiny $v^*$};
\end{scope}
\end{tikzpicture}}}, 
\quad 
\vcenter{\hbox{\begin{tikzpicture}[scale=0.65]
    \begin{scope}[shift={(0,1.5)}]
    \draw [blue] (-0.5, 0.8)--(-0.5, 0) .. controls +(0, -0.6) and +(0,-0.6).. (0.5, 0)--(0.5, 0.8);    
\begin{scope}[shift={(0.5, 0.3)}]
\draw [fill=white] (-0.3, -0.3) rectangle (0.3, 0.3);
\node at (0, 0) {\tiny $v^*$};
\end{scope}
    \end{scope}
\draw [blue] (-0.5, -0.8)--(-0.5, 0) .. controls +(0, 0.6) and +(0,0.6).. (0.5, 0)--(0.5, -0.8);
\begin{scope}[shift={(0.5, -0.3)}]
\draw [fill=white] (-0.3, -0.3) rectangle (0.3, 0.3);
\node at (0, 0) {\tiny $v$};
\end{scope}
\end{tikzpicture}}}.
\end{align*}
We denote by $\check{\partial}_1$, $\check{\partial}_2$ the associated twisted derivations.

To check the intertwining property, we have to evaluate $\check{\partial}_j\cL-\cL\check{\partial}_j$ for $j=1, 2$.
The following two items
\begin{align*}
\frac{1}{2}\vcenter{\hbox{\begin{tikzpicture}[scale=0.65]
     \draw [blue] (1.4, 2.3)--(1.4, -0.8);
    \begin{scope}[shift={(0,1.5)}]
    \draw [blue] (-0.5, 0.8)--(-0.5, 0) .. controls +(0, -0.6) and +(0,-0.6).. (0.5, 0)--(0.5, 0.8); 
\begin{scope}[shift={(0.5, 0.3)}]
\draw [fill=white] (-0.5, -0.3) rectangle (0.5, 0.3);
\node at (0, 0) {\tiny $\mathbf{y} v$};
\end{scope}
    \end{scope}
\draw [blue] (-0.5, -0.8)--(-0.5, 0) .. controls +(0, 0.6) and +(0,0.6).. (0.5, 0)--(0.5, -0.8);
\begin{scope}[shift={(0.5, -0.3)}]
\draw [fill=white] (-0.3, -0.3) rectangle (0.3, 0.3);
\node at (0, 0) {\tiny $w$};
\end{scope}
\begin{scope}[shift={(1.4, 0.8)}]
\draw [fill=white] (-0.3, -0.3) rectangle (0.3, 0.3);
\node at (0, 0) {\tiny $\overline{v^*}$};  
\end{scope}
\end{tikzpicture}}}
+
\frac{1}{2}\vcenter{\hbox{\begin{tikzpicture}[scale=0.65]
     \draw [blue] (1.4, 2.3)--(1.4, -0.8);
    \begin{scope}[shift={(0,1.5)}]
    \draw [blue] (-0.5, 0.8)--(-0.5, 0) .. controls +(0, -0.6) and +(0,-0.6).. (0.5, 0)--(0.5, 0.8); 
\begin{scope}[shift={(0.5, 0.3)}]
\draw [fill=white] (-0.4, -0.3) rectangle (0.4, 0.3);
\node at (0, 0) {\tiny $v$};
\end{scope}
    \end{scope}
\draw [blue] (-0.5, -0.8)--(-0.5, 0) .. controls +(0, 0.6) and +(0,0.6).. (0.5, 0)--(0.5, -0.8);
\begin{scope}[shift={(0.5, -0.3)}]
\draw [fill=white] (-0.3, -0.3) rectangle (0.3, 0.3);
\node at (0, 0) {\tiny $w\mathbf{y}$};
\end{scope}
\begin{scope}[shift={(1.4, 0.8)}]
\draw [fill=white] (-0.3, -0.3) rectangle (0.3, 0.3);
\node at (0, 0) {\tiny $\overline{v^*}$};  
\end{scope}
\end{tikzpicture}}}
-\frac{1}{2} \mu\vcenter{\hbox{\begin{tikzpicture}[scale=0.65]
     \draw [blue] (1.4, 2.3)--(1.4, -0.8);
    \begin{scope}[shift={(0,1.5)}]
    \draw [blue] (-0.5, 0.8)--(-0.5, 0) .. controls +(0, -0.6) and +(0,-0.6).. (0.5, 0)--(0.5, 0.8); 
\begin{scope}[shift={(0.5, 0.3)}]
\draw [fill=white] (-0.55, -0.3) rectangle (0.55, 0.3);
\node at (0, 0) {\tiny $v^2$};
\end{scope}
    \end{scope}
\draw [blue] (-0.5, -0.8)--(-0.5, 0) .. controls +(0, 0.6) and +(0,0.6).. (0.5, 0)--(0.5, -0.8);
\begin{scope}[shift={(0.5, -0.3)}]
\draw [fill=white] (-0.3, -0.3) rectangle (0.3, 0.3);
\node at (0, 0) {\tiny $wv^*$};
\end{scope}
\begin{scope}[shift={(1.4, 0.8)}]
\draw [fill=white] (-0.3, -0.3) rectangle (0.3, 0.3);
\node at (0, 0) {\tiny $\overline{v^*}$};  
\end{scope}
\end{tikzpicture}}}
-\frac{1}{2} \mu^{-1}\vcenter{\hbox{\begin{tikzpicture}[scale=0.65]
     \draw [blue] (1.4, 2.3)--(1.4, -0.8);
    \begin{scope}[shift={(0,1.5)}]
    \draw [blue] (-0.5, 0.8)--(-0.5, 0) .. controls +(0, -0.6) and +(0,-0.6).. (0.5, 0)--(0.5, 0.8); 
\begin{scope}[shift={(0.5, 0.3)}]
\draw [fill=white] (-0.55, -0.3) rectangle (0.55, 0.3);
\node at (0, 0) {\tiny $v^*v$};
\end{scope}
    \end{scope}
\draw [blue] (-0.5, -0.8)--(-0.5, 0) .. controls +(0, 0.6) and +(0,0.6).. (0.5, 0)--(0.5, -0.8);
\begin{scope}[shift={(0.5, -0.3)}]
\draw [fill=white] (-0.3, -0.3) rectangle (0.3, 0.3);
\node at (0, 0) {\tiny $wv$};
\end{scope}
\begin{scope}[shift={(1.4, 0.8)}]
\draw [fill=white] (-0.3, -0.3) rectangle (0.3, 0.3);
\node at (0, 0) {\tiny $\overline{v^*}$};  
\end{scope}
\end{tikzpicture}}}
-\frac{1}{2} (b-a)\vcenter{\hbox{\begin{tikzpicture}[scale=0.65]
     \draw [blue] (1.4, 2.3)--(1.4, -0.8);
    \begin{scope}[shift={(0,1.5)}]
    \draw [blue] (-0.5, 0.8)--(-0.5, 0) .. controls +(0, -0.6) and +(0,-0.6).. (0.5, 0)--(0.5, 0.8); 
\begin{scope}[shift={(0.5, 0.3)}]
\draw [fill=white] (-0.55, -0.3) rectangle (0.55, 0.3);
\node at (0, 0) {\tiny $v^2$};
\end{scope}
    \end{scope}
\draw [blue] (-0.5, -0.8)--(-0.5, 0) .. controls +(0, 0.6) and +(0,0.6).. (0.5, 0)--(0.5, -0.8);
\begin{scope}[shift={(0.5, -0.3)}]
\draw [fill=white] (-0.3, -0.3) rectangle (0.3, 0.3);
\node at (0, 0) {\tiny $wv$};
\end{scope}
\begin{scope}[shift={(1.4, 0.8)}]
\draw [fill=white] (-0.3, -0.3) rectangle (0.3, 0.3);
\node at (0, 0) {\tiny $\overline{v^*}$};  
\end{scope}
\end{tikzpicture}}}
-\frac{1}{2} (b-a)\vcenter{\hbox{\begin{tikzpicture}[scale=0.65]
     \draw [blue] (1.4, 2.3)--(1.4, -0.8);
    \begin{scope}[shift={(0,1.5)}]
    \draw [blue] (-0.5, 0.8)--(-0.5, 0) .. controls +(0, -0.6) and +(0,-0.6).. (0.5, 0)--(0.5, 0.8); 
\begin{scope}[shift={(0.5, 0.3)}]
\draw [fill=white] (-0.55, -0.3) rectangle (0.55, 0.3);
\node at (0, 0) {\tiny $vv^*$};
\end{scope}
    \end{scope}
\draw [blue] (-0.5, -0.8)--(-0.5, 0) .. controls +(0, 0.6) and +(0,0.6).. (0.5, 0)--(0.5, -0.8);
\begin{scope}[shift={(0.5, -0.3)}]
\draw [fill=white] (-0.3, -0.3) rectangle (0.3, 0.3);
\node at (0, 0) {\tiny $wv^*$};
\end{scope}
\begin{scope}[shift={(1.4, 0.8)}]
\draw [fill=white] (-0.3, -0.3) rectangle (0.3, 0.3);
\node at (0, 0) {\tiny $\overline{v^*}$};  
\end{scope}
\end{tikzpicture}}},
\end{align*}
and
\begin{align*}
\frac{1}{2}\vcenter{\hbox{\begin{tikzpicture}[scale=0.65]
     \draw [blue] (1.4, 2.3)--(1.4, -0.8);
    \begin{scope}[shift={(0,1.5)}]
    \draw [blue] (-0.5, 0.8)--(-0.5, 0) .. controls +(0, -0.6) and +(0,-0.6).. (0.5, 0)--(0.5, 0.8); 
\begin{scope}[shift={(0.5, 0.3)}]
\draw [fill=white] (-0.5, -0.3) rectangle (0.5, 0.3);
\node at (0, 0) {\tiny $\mathbf{y}$};
\end{scope}
    \end{scope}
\draw [blue] (-0.5, -0.8)--(-0.5, 0) .. controls +(0, 0.6) and +(0,0.6).. (0.5, 0)--(0.5, -0.8);
\begin{scope}[shift={(0.5, -0.3)}]
\draw [fill=white] (-0.4, -0.3) rectangle (0.4, 0.3);
\node at (0, 0) {\tiny $w v $};
\end{scope}
\begin{scope}[shift={(1.4, 0.8)}]
\draw [fill=white] (-0.3, -0.3) rectangle (0.3, 0.3);
\node at (0, 0) {\tiny $\overline{v^*}$};  
\end{scope}
\end{tikzpicture}}}
+
\frac{1}{2}\vcenter{\hbox{\begin{tikzpicture}[scale=0.65]
     \draw [blue] (1.4, 2.3)--(1.4, -0.8);
    \begin{scope}[shift={(0,1.5)}]
    \draw [blue] (-0.5, 0.8)--(-0.5, 0) .. controls +(0, -0.6) and +(0,-0.6).. (0.5, 0)--(0.5, 0.8); 
    \end{scope}
\draw [blue] (-0.5, -0.8)--(-0.5, 0) .. controls +(0, 0.6) and +(0,0.6).. (0.5, 0)--(0.5, -0.8);
\begin{scope}[shift={(0.5, -0.3)}]
\draw [fill=white] (-0.55, -0.3) rectangle (0.55, 0.3);
\node at (0, 0) {\tiny $wv\mathbf{y}$};
\end{scope}
\begin{scope}[shift={(1.4, 0.8)}]
\draw [fill=white] (-0.3, -0.3) rectangle (0.3, 0.3);
\node at (0, 0) {\tiny $\overline{v^*}$};  
\end{scope}
\end{tikzpicture}}}
-\frac{1}{2}\mu \vcenter{\hbox{\begin{tikzpicture}[scale=0.65]
     \draw [blue] (1.4, 2.3)--(1.4, -0.8);
    \begin{scope}[shift={(0,1.5)}]
    \draw [blue] (-0.5, 0.8)--(-0.5, 0) .. controls +(0, -0.6) and +(0,-0.6).. (0.5, 0)--(0.5, 0.8); 
\begin{scope}[shift={(0.5, 0.3)}]
\draw [fill=white] (-0.3, -0.3) rectangle (0.3, 0.3);
\node at (0, 0) {\tiny $v$};
\end{scope}
    \end{scope}
\draw [blue] (-0.5, -0.8)--(-0.5, 0) .. controls +(0, 0.6) and +(0,0.6).. (0.5, 0)--(0.5, -0.8);
\begin{scope}[shift={(0.5, -0.3)}]
\draw [fill=white] (-0.55, -0.3) rectangle (0.55, 0.3);
\node at (0, 0) {\tiny $wvv^*$};
\end{scope}
\begin{scope}[shift={(1.4, 0.8)}]
\draw [fill=white] (-0.3, -0.3) rectangle (0.3, 0.3);
\node at (0, 0) {\tiny $\overline{v^*}$};  
\end{scope}
\end{tikzpicture}}}
-\frac{1}{2}\mu^{-1}\vcenter{\hbox{\begin{tikzpicture}[scale=0.65]
     \draw [blue] (1.4, 2.3)--(1.4, -0.8);
    \begin{scope}[shift={(0,1.5)}]
    \draw [blue] (-0.5, 0.8)--(-0.5, 0) .. controls +(0, -0.6) and +(0,-0.6).. (0.5, 0)--(0.5, 0.8); 
\begin{scope}[shift={(0.5, 0.3)}]
\draw [fill=white] (-0.3, -0.3) rectangle (0.3, 0.3);
\node at (0, 0) {\tiny $v^*$};
\end{scope}
    \end{scope}
\draw [blue] (-0.5, -0.8)--(-0.5, 0) .. controls +(0, 0.6) and +(0,0.6).. (0.5, 0)--(0.5, -0.8);
\begin{scope}[shift={(0.5, -0.3)}]
\draw [fill=white] (-0.55, -0.3) rectangle (0.55, 0.3);
\node at (0, 0) {\tiny $wv^2$};
\end{scope}
\begin{scope}[shift={(1.4, 0.8)}]
\draw [fill=white] (-0.3, -0.3) rectangle (0.3, 0.3);
\node at (0, 0) {\tiny $\overline{v^*}$};  
\end{scope}
\end{tikzpicture}}}
-\frac{1}{2}(b-a)\vcenter{\hbox{\begin{tikzpicture}[scale=0.65]
     \draw [blue] (1.4, 2.3)--(1.4, -0.8);
    \begin{scope}[shift={(0,1.5)}]
    \draw [blue] (-0.5, 0.8)--(-0.5, 0) .. controls +(0, -0.6) and +(0,-0.6).. (0.5, 0)--(0.5, 0.8); 
\begin{scope}[shift={(0.5, 0.3)}]
\draw [fill=white] (-0.3, -0.3) rectangle (0.3, 0.3);
\node at (0, 0) {\tiny $v$};
\end{scope}
    \end{scope}
\draw [blue] (-0.5, -0.8)--(-0.5, 0) .. controls +(0, 0.6) and +(0,0.6).. (0.5, 0)--(0.5, -0.8);
\begin{scope}[shift={(0.5, -0.3)}]
\draw [fill=white] (-0.55, -0.3) rectangle (0.55, 0.3);
\node at (0, 0) {\tiny $wv^2$};
\end{scope}
\begin{scope}[shift={(1.4, 0.8)}]
\draw [fill=white] (-0.3, -0.3) rectangle (0.3, 0.3);
\node at (0, 0) {\tiny $\overline{v^*}$};  
\end{scope}
\end{tikzpicture}}}
-\frac{1}{2}(b-a)\vcenter{\hbox{\begin{tikzpicture}[scale=0.65]
     \draw [blue] (1.4, 2.3)--(1.4, -0.8);
    \begin{scope}[shift={(0,1.5)}]
    \draw [blue] (-0.5, 0.8)--(-0.5, 0) .. controls +(0, -0.6) and +(0,-0.6).. (0.5, 0)--(0.5, 0.8); 
\begin{scope}[shift={(0.5, 0.3)}]
\draw [fill=white] (-0.3, -0.3) rectangle (0.3, 0.3);
\node at (0, 0) {\tiny $v^*$};
\end{scope}
    \end{scope}
\draw [blue] (-0.5, -0.8)--(-0.5, 0) .. controls +(0, 0.6) and +(0,0.6).. (0.5, 0)--(0.5, -0.8);
\begin{scope}[shift={(0.5, -0.3)}]
\draw [fill=white] (-0.55, -0.3) rectangle (0.55, 0.3);
\node at (0, 0) {\tiny $wvv^*$};
\end{scope}
\begin{scope}[shift={(1.4, 0.8)}]
\draw [fill=white] (-0.3, -0.3) rectangle (0.3, 0.3);
\node at (0, 0) {\tiny $\overline{v^*}$};  
\end{scope}
\end{tikzpicture}}}
\end{align*}
is the Fourier multiplier of $\check{\partial}_1\cL$.
The following two items
\begin{align*}
 \frac{1}{2}\vcenter{\hbox{\begin{tikzpicture}[scale=0.65]
     \draw [blue] (1.4, 2.3)--(1.4, -0.8);
    \begin{scope}[shift={(0,1.5)}]
    \draw [blue] (-0.5, 0.8)--(-0.5, 0) .. controls +(0, -0.6) and +(0,-0.6).. (0.5, 0)--(0.5, 0.8); 
\begin{scope}[shift={(0.5, 0.3)}]
\draw [fill=white] (-0.5, -0.3) rectangle (0.5, 0.3);
\node at (0, 0) {\tiny $v \mathbf{y} $};
\end{scope}
    \end{scope}
\draw [blue] (-0.5, -0.8)--(-0.5, 0) .. controls +(0, 0.6) and +(0,0.6).. (0.5, 0)--(0.5, -0.8);
\begin{scope}[shift={(0.5, -0.3)}]
\draw [fill=white] (-0.3, -0.3) rectangle (0.3, 0.3);
\node at (0, 0) {\tiny $w$};
\end{scope}
\begin{scope}[shift={(1.4, 0.8)}]
\draw [fill=white] (-0.3, -0.3) rectangle (0.3, 0.3);
\node at (0, 0) {\tiny $\overline{v^*}$};  
\end{scope}
\end{tikzpicture}}}
+
\frac{1}{2}\vcenter{\hbox{\begin{tikzpicture}[scale=0.65]
     \draw [blue] (1.4, 2.3)--(1.4, -0.8);
    \begin{scope}[shift={(0,1.5)}]
    \draw [blue] (-0.5, 0.8)--(-0.5, 0) .. controls +(0, -0.6) and +(0,-0.6).. (0.5, 0)--(0.5, 0.8); 
\begin{scope}[shift={(0.5, 0.3)}]
\draw [fill=white] (-0.4, -0.3) rectangle (0.4, 0.3);
\node at (0, 0) {\tiny $v$};
\end{scope}
    \end{scope}
\draw [blue] (-0.5, -0.8)--(-0.5, 0) .. controls +(0, 0.6) and +(0,0.6).. (0.5, 0)--(0.5, -0.8);
\begin{scope}[shift={(0.5, -0.3)}]
\draw [fill=white] (-0.3, -0.3) rectangle (0.3, 0.3);
\node at (0, 0) {\tiny $\mathbf{y}w$};
\end{scope}
\begin{scope}[shift={(1.4, 0.8)}]
\draw [fill=white] (-0.3, -0.3) rectangle (0.3, 0.3);
\node at (0, 0) {\tiny $\overline{v^*}$};  
\end{scope}
\end{tikzpicture}}}
-\frac{1}{2} \mu\vcenter{\hbox{\begin{tikzpicture}[scale=0.65]
     \draw [blue] (1.4, 2.3)--(1.4, -0.8);
    \begin{scope}[shift={(0,1.5)}]
    \draw [blue] (-0.5, 0.8)--(-0.5, 0) .. controls +(0, -0.6) and +(0,-0.6).. (0.5, 0)--(0.5, 0.8); 
\begin{scope}[shift={(0.5, 0.3)}]
\draw [fill=white] (-0.55, -0.3) rectangle (0.55, 0.3);
\node at (0, 0) {\tiny $v^2$};
\end{scope}
    \end{scope}
\draw [blue] (-0.5, -0.8)--(-0.5, 0) .. controls +(0, 0.6) and +(0,0.6).. (0.5, 0)--(0.5, -0.8);
\begin{scope}[shift={(0.5, -0.3)}]
\draw [fill=white] (-0.4, -0.3) rectangle (0.4, 0.3);
\node at (0, 0) {\tiny $v^*w$};
\end{scope}
\begin{scope}[shift={(1.4, 0.8)}]
\draw [fill=white] (-0.3, -0.3) rectangle (0.3, 0.3);
\node at (0, 0) {\tiny $\overline{v^*}$};  
\end{scope}
\end{tikzpicture}}}
-\frac{1}{2}\mu^{-1}\vcenter{\hbox{\begin{tikzpicture}[scale=0.65]
     \draw [blue] (1.4, 2.3)--(1.4, -0.8);
    \begin{scope}[shift={(0,1.5)}]
    \draw [blue] (-0.5, 0.8)--(-0.5, 0) .. controls +(0, -0.6) and +(0,-0.6).. (0.5, 0)--(0.5, 0.8); 
\begin{scope}[shift={(0.5, 0.3)}]
\draw [fill=white] (-0.55, -0.3) rectangle (0.55, 0.3);
\node at (0, 0) {\tiny $vv^*$};
\end{scope}
    \end{scope}
\draw [blue] (-0.5, -0.8)--(-0.5, 0) .. controls +(0, 0.6) and +(0,0.6).. (0.5, 0)--(0.5, -0.8);
\begin{scope}[shift={(0.5, -0.3)}]
\draw [fill=white] (-0.3, -0.3) rectangle (0.3, 0.3);
\node at (0, 0) {\tiny $vw$};
\end{scope}
\begin{scope}[shift={(1.4, 0.8)}]
\draw [fill=white] (-0.3, -0.3) rectangle (0.3, 0.3);
\node at (0, 0) {\tiny $\overline{v^*}$};  
\end{scope}
\end{tikzpicture}}}  
-\frac{1}{2} (b-a)\vcenter{\hbox{\begin{tikzpicture}[scale=0.65]
     \draw [blue] (1.4, 2.3)--(1.4, -0.8);
    \begin{scope}[shift={(0,1.5)}]
    \draw [blue] (-0.5, 0.8)--(-0.5, 0) .. controls +(0, -0.6) and +(0,-0.6).. (0.5, 0)--(0.5, 0.8); 
\begin{scope}[shift={(0.5, 0.3)}]
\draw [fill=white] (-0.55, -0.3) rectangle (0.55, 0.3);
\node at (0, 0) {\tiny $v^2$};
\end{scope}
    \end{scope}
\draw [blue] (-0.5, -0.8)--(-0.5, 0) .. controls +(0, 0.6) and +(0,0.6).. (0.5, 0)--(0.5, -0.8);
\begin{scope}[shift={(0.5, -0.3)}]
\draw [fill=white] (-0.3, -0.3) rectangle (0.3, 0.3);
\node at (0, 0) {\tiny $vw$};
\end{scope}
\begin{scope}[shift={(1.4, 0.8)}]
\draw [fill=white] (-0.3, -0.3) rectangle (0.3, 0.3);
\node at (0, 0) {\tiny $\overline{v^*}$};  
\end{scope}
\end{tikzpicture}}}
-\frac{1}{2} (b-a)\vcenter{\hbox{\begin{tikzpicture}[scale=0.65]
     \draw [blue] (1.4, 2.3)--(1.4, -0.8);
    \begin{scope}[shift={(0,1.5)}]
    \draw [blue] (-0.5, 0.8)--(-0.5, 0) .. controls +(0, -0.6) and +(0,-0.6).. (0.5, 0)--(0.5, 0.8); 
\begin{scope}[shift={(0.5, 0.3)}]
\draw [fill=white] (-0.55, -0.3) rectangle (0.55, 0.3);
\node at (0, 0) {\tiny $vv^*$};
\end{scope}
    \end{scope}
\draw [blue] (-0.5, -0.8)--(-0.5, 0) .. controls +(0, 0.6) and +(0,0.6).. (0.5, 0)--(0.5, -0.8);
\begin{scope}[shift={(0.5, -0.3)}]
\draw [fill=white] (-0.3, -0.3) rectangle (0.3, 0.3);
\node at (0, 0) {\tiny $v^*w$};
\end{scope}
\begin{scope}[shift={(1.4, 0.8)}]
\draw [fill=white] (-0.3, -0.3) rectangle (0.3, 0.3);
\node at (0, 0) {\tiny $\overline{v^*}$};  
\end{scope}
\end{tikzpicture}}}
\end{align*}
and
\begin{align*}
\frac{1}{2}\vcenter{\hbox{\begin{tikzpicture}[scale=0.65]
     \draw [blue] (1.4, 2.3)--(1.4, -0.8);
    \begin{scope}[shift={(0,1.5)}]
    \draw [blue] (-0.5, 0.8)--(-0.5, 0) .. controls +(0, -0.6) and +(0,-0.6).. (0.5, 0)--(0.5, 0.8); 
\begin{scope}[shift={(0.5, 0.3)}]
\draw [fill=white] (-0.5, -0.3) rectangle (0.5, 0.3);
\node at (0, 0) {\tiny $\mathbf{y}$};
\end{scope}
    \end{scope}
\draw [blue] (-0.5, -0.8)--(-0.5, 0) .. controls +(0, 0.6) and +(0,0.6).. (0.5, 0)--(0.5, -0.8);
\begin{scope}[shift={(0.5, -0.3)}]
\draw [fill=white] (-0.4, -0.3) rectangle (0.4, 0.3);
\node at (0, 0) {\tiny $wv $};
\end{scope}
\begin{scope}[shift={(1.4, 0.8)}]
\draw [fill=white] (-0.3, -0.3) rectangle (0.3, 0.3);
\node at (0, 0) {\tiny $\overline{v^*}$};  
\end{scope}
\end{tikzpicture}}}
+
\frac{1}{2}\vcenter{\hbox{\begin{tikzpicture}[scale=0.65]
     \draw [blue] (1.4, 2.3)--(1.4, -0.8);
    \begin{scope}[shift={(0,1.5)}]
    \draw [blue] (-0.5, 0.8)--(-0.5, 0) .. controls +(0, -0.6) and +(0,-0.6).. (0.5, 0)--(0.5, 0.8); 
    \end{scope}
\draw [blue] (-0.5, -0.8)--(-0.5, 0) .. controls +(0, 0.6) and +(0,0.6).. (0.5, 0)--(0.5, -0.8);
\begin{scope}[shift={(0.5, -0.3)}]
\draw [fill=white] (-0.55, -0.3) rectangle (0.55, 0.3);
\node at (0, 0) {\tiny $\mathbf{y}wv_k$};
\end{scope}
\begin{scope}[shift={(1.4, 0.8)}]
\draw [fill=white] (-0.3, -0.3) rectangle (0.3, 0.3);
\node at (0, 0) {\tiny $\overline{v^*}$};  
\end{scope}
\end{tikzpicture}}}
-\frac{1}{2}\mu \vcenter{\hbox{\begin{tikzpicture}[scale=0.65]
     \draw [blue] (1.4, 2.3)--(1.4, -0.8);
    \begin{scope}[shift={(0,1.5)}]
    \draw [blue] (-0.5, 0.8)--(-0.5, 0) .. controls +(0, -0.6) and +(0,-0.6).. (0.5, 0)--(0.5, 0.8); 
\begin{scope}[shift={(0.5, 0.3)}]
\draw [fill=white] (-0.3, -0.3) rectangle (0.3, 0.3);
\node at (0, 0) {\tiny $v$};
\end{scope}
    \end{scope}
\draw [blue] (-0.5, -0.8)--(-0.5, 0) .. controls +(0, 0.6) and +(0,0.6).. (0.5, 0)--(0.5, -0.8);
\begin{scope}[shift={(0.5, -0.3)}]
\draw [fill=white] (-0.55, -0.3) rectangle (0.55, 0.3);
\node at (0, 0) {\tiny $v^*wv$};
\end{scope}
\begin{scope}[shift={(1.4, 0.8)}]
\draw [fill=white] (-0.3, -0.3) rectangle (0.3, 0.3);
\node at (0, 0) {\tiny $\overline{v^*}$};  
\end{scope}
\end{tikzpicture}}}
-\frac{1}{2}\mu^{-1}\vcenter{\hbox{\begin{tikzpicture}[scale=0.65]
     \draw [blue] (1.4, 2.3)--(1.4, -0.8);
    \begin{scope}[shift={(0,1.5)}]
    \draw [blue] (-0.5, 0.8)--(-0.5, 0) .. controls +(0, -0.6) and +(0,-0.6).. (0.5, 0)--(0.5, 0.8); 
\begin{scope}[shift={(0.5, 0.3)}]
\draw [fill=white] (-0.3, -0.3) rectangle (0.3, 0.3);
\node at (0, 0) {\tiny $v^*$};
\end{scope}
    \end{scope}
\draw [blue] (-0.5, -0.8)--(-0.5, 0) .. controls +(0, 0.6) and +(0,0.6).. (0.5, 0)--(0.5, -0.8);
\begin{scope}[shift={(0.5, -0.3)}]
\draw [fill=white] (-0.55, -0.3) rectangle (0.55, 0.3);
\node at (0, 0) {\tiny $ vwv$};
\end{scope}
\begin{scope}[shift={(1.4, 0.8)}]
\draw [fill=white] (-0.3, -0.3) rectangle (0.3, 0.3);
\node at (0, 0) {\tiny $\overline{v^*}$};  
\end{scope}
\end{tikzpicture}}}
-\frac{1}{2}(b-a)\vcenter{\hbox{\begin{tikzpicture}[scale=0.65]
     \draw [blue] (1.4, 2.3)--(1.4, -0.8);
    \begin{scope}[shift={(0,1.5)}]
    \draw [blue] (-0.5, 0.8)--(-0.5, 0) .. controls +(0, -0.6) and +(0,-0.6).. (0.5, 0)--(0.5, 0.8); 
\begin{scope}[shift={(0.5, 0.3)}]
\draw [fill=white] (-0.3, -0.3) rectangle (0.3, 0.3);
\node at (0, 0) {\tiny $v$};
\end{scope}
    \end{scope}
\draw [blue] (-0.5, -0.8)--(-0.5, 0) .. controls +(0, 0.6) and +(0,0.6).. (0.5, 0)--(0.5, -0.8);
\begin{scope}[shift={(0.5, -0.3)}]
\draw [fill=white] (-0.55, -0.3) rectangle (0.55, 0.3);
\node at (0, 0) {\tiny $vwv$};
\end{scope}
\begin{scope}[shift={(1.4, 0.8)}]
\draw [fill=white] (-0.3, -0.3) rectangle (0.3, 0.3);
\node at (0, 0) {\tiny $\overline{v^*}$};  
\end{scope}
\end{tikzpicture}}}
-\frac{1}{2}(b-a)\vcenter{\hbox{\begin{tikzpicture}[scale=0.65]
     \draw [blue] (1.4, 2.3)--(1.4, -0.8);
    \begin{scope}[shift={(0,1.5)}]
    \draw [blue] (-0.5, 0.8)--(-0.5, 0) .. controls +(0, -0.6) and +(0,-0.6).. (0.5, 0)--(0.5, 0.8); 
\begin{scope}[shift={(0.5, 0.3)}]
\draw [fill=white] (-0.3, -0.3) rectangle (0.3, 0.3);
\node at (0, 0) {\tiny $v^*$};
\end{scope}
    \end{scope}
\draw [blue] (-0.5, -0.8)--(-0.5, 0) .. controls +(0, 0.6) and +(0,0.6).. (0.5, 0)--(0.5, -0.8);
\begin{scope}[shift={(0.5, -0.3)}]
\draw [fill=white] (-0.55, -0.3) rectangle (0.55, 0.3);
\node at (0, 0) {\tiny $vwv^*$};
\end{scope}
\begin{scope}[shift={(1.4, 0.8)}]
\draw [fill=white] (-0.3, -0.3) rectangle (0.3, 0.3);
\node at (0, 0) {\tiny $\overline{v^*}$};  
\end{scope}
\end{tikzpicture}}}
\end{align*}
are the Fourier multiplier of $\cL \check{\partial}_1$.
Now we see that $\check{\partial}_1\cL-\cL\check{\partial}_1$ is depicted as
\begin{align*}
\frac{1}{2} (\mu^{-1}+\mu)
 \vcenter{\hbox{\begin{tikzpicture}[scale=0.65]
     \draw [blue] (1.4, 2.3)--(1.4, -0.8);
    \begin{scope}[shift={(0,1.5)}]
    \draw [blue] (-0.5, 0.8)--(-0.5, 0) .. controls +(0, -0.6) and +(0,-0.6).. (0.5, 0)--(0.5, 0.8); 
    \end{scope}
\draw [blue] (-0.5, -0.8)--(-0.5, 0) .. controls +(0, 0.6) and +(0,0.6).. (0.5, 0)--(0.5, -0.8);
\begin{scope}[shift={(0.5, -0.3)}]
\draw [fill=white] (-0.3, -0.3) rectangle (0.3, 0.3);
\node at (0, 0) {\tiny $ wv$};
\end{scope}
\begin{scope}[shift={(1.4, 0.8)}]
\draw [fill=white] (-0.3, -0.3) rectangle (0.3, 0.3);
\node at (0, 0) {\tiny $\overline{v^*}$};  
\end{scope}
\end{tikzpicture}}}  
- \frac{1}{2} (\mu+\mu^{-1}) \vcenter{\hbox{\begin{tikzpicture}[scale=0.65]
     \draw [blue] (1.4, 2.3)--(1.4, -0.8);
    \begin{scope}[shift={(0,1.5)}]
    \draw [blue] (-0.5, 0.8)--(-0.5, 0) .. controls +(0, -0.6) and +(0,-0.6).. (0.5, 0)--(0.5, 0.8); 
\begin{scope}[shift={(0.5, 0.3)}]
\draw [fill=white] (-0.3, -0.3) rectangle (0.3, 0.3);
\node at (0, 0) {\tiny $v$};
\end{scope}
    \end{scope}
\draw [blue] (-0.5, -0.8)--(-0.5, 0) .. controls +(0, 0.6) and +(0,0.6).. (0.5, 0)--(0.5, -0.8);
\begin{scope}[shift={(0.5, -0.3)}]
\draw [fill=white] (-0.55, -0.3) rectangle (0.55, 0.3);
\node at (0, 0) {\tiny $w$};
\end{scope}
\begin{scope}[shift={(1.4, 0.8)}]
\draw [fill=white] (-0.3, -0.3) rectangle (0.3, 0.3);
\node at (0, 0) {\tiny $\overline{v^*}$};  
\end{scope}
\end{tikzpicture}}}
-(b-a)\vcenter{\hbox{\begin{tikzpicture}[scale=0.65]
     \draw [blue] (1.4, 2.3)--(1.4, -0.8);
    \begin{scope}[shift={(0,1.5)}]
    \draw [blue] (-0.5, 0.8)--(-0.5, 0) .. controls +(0, -0.6) and +(0,-0.6).. (0.5, 0)--(0.5, 0.8); 
    \end{scope}
\draw [blue] (-0.5, -0.8)--(-0.5, 0) .. controls +(0, 0.6) and +(0,0.6).. (0.5, 0)--(0.5, -0.8);
\begin{scope}[shift={(0.5, -0.3)}]
\draw [fill=white] (-0.3, -0.3) rectangle (0.3, 0.3);
\node at (0, 0) {\tiny $wv^*$};
\end{scope}
\begin{scope}[shift={(1.4, 0.8)}]
\draw [fill=white] (-0.3, -0.3) rectangle (0.3, 0.3);
\node at (0, 0) {\tiny $\overline{v^*}$};  
\end{scope}
\end{tikzpicture}}}
+(b-a)\vcenter{\hbox{\begin{tikzpicture}[scale=0.65]
     \draw [blue] (1.4, 2.3)--(1.4, -0.8);
    \begin{scope}[shift={(0,1.5)}]
    \draw [blue] (-0.5, 0.8)--(-0.5, 0) .. controls +(0, -0.6) and +(0,-0.6).. (0.5, 0)--(0.5, 0.8); 
\begin{scope}[shift={(0.5, 0.3)}]
\draw [fill=white] (-0.3, -0.3) rectangle (0.3, 0.3);
\node at (0, 0) {\tiny $v^*$};
\end{scope}
    \end{scope}
\draw [blue] (-0.5, -0.8)--(-0.5, 0) .. controls +(0, 0.6) and +(0,0.6).. (0.5, 0)--(0.5, -0.8);
\begin{scope}[shift={(0.5, -0.3)}]
\draw [fill=white] (-0.55, -0.3) rectangle (0.55, 0.3);
\node at (0, 0) {\tiny $w$};
\end{scope}
\begin{scope}[shift={(1.4, 0.8)}]
\draw [fill=white] (-0.3, -0.3) rectangle (0.3, 0.3);
\node at (0, 0) {\tiny $\overline{v^*}$};  
\end{scope}
\end{tikzpicture}}}.
\end{align*}
Similar result is true for $\check{\partial}_2$.
This implies that 
\begin{align*}
    \check{\partial}_1 \cL -\cL\check{\partial}_1 =\frac{1}{2}(\mu+\mu^{-1})\check{\partial}_1+(a-b)\check{\partial}_2, \\
    \check{\partial}_2 \cL -\cL\check{\partial}_2 =\frac{1}{2}(\mu+\mu^{-1})\check{\partial}_2+(a-b)\check{\partial}_1.
\end{align*}
Let 
\begin{align*}
F_{1,1}=\lambda^{1/2}\,\vcenter{\hbox{\begin{tikzpicture}[scale=0.65]
    \begin{scope}[shift={(0,1.5)}]
    \draw [blue] (-0.5, 0.8)--(-0.5, 0) .. controls +(0, -0.6) and +(0,-0.6).. (0.5, 0)--(0.5, 0.8); 
\begin{scope}[shift={(0.5, 0.3)}]
\draw [fill=white] (-0.3, -0.3) rectangle (0.3, 0.3);
\node at (0, 0) {\tiny $v$};
\end{scope}
    \end{scope}
\draw [blue] (-0.5, -0.8)--(-0.5, 0) .. controls +(0, 0.6) and +(0,0.6).. (0.5, 0)--(0.5, -0.8);
\begin{scope}[shift={(0.5, -0.3)}]
\draw [fill=white] (-0.3, -0.3) rectangle (0.3, 0.3);
\node at (0, 0) {\tiny $v^*$};
\end{scope}
\end{tikzpicture}}},\,
F_{1,2}=\lambda^{1/2}\,\vcenter{\hbox{\begin{tikzpicture}[scale=0.65]
    \begin{scope}[shift={(0,1.5)}]
    \draw [blue] (-0.5, 0.8)--(-0.5, 0) .. controls +(0, -0.6) and +(0,-0.6).. (0.5, 0)--(0.5, 0.8); 
\begin{scope}[shift={(0.5, 0.3)}]
\draw [fill=white] (-0.3, -0.3) rectangle (0.3, 0.3);
\node at (0, 0) {\tiny $v^*$};
\end{scope}
    \end{scope}
\draw [blue] (-0.5, -0.8)--(-0.5, 0) .. controls +(0, 0.6) and +(0,0.6).. (0.5, 0)--(0.5, -0.8);
\begin{scope}[shift={(0.5, -0.3)}]
\draw [fill=white] (-0.3, -0.3) rectangle (0.3, 0.3);
\node at (0, 0) {\tiny $v^*$};
\end{scope}
\end{tikzpicture}}},\,
F_{2,1}=\lambda^{1/2}\,\vcenter{\hbox{\begin{tikzpicture}[scale=0.65]
    \begin{scope}[shift={(0,1.5)}]
    \draw [blue] (-0.5, 0.8)--(-0.5, 0) .. controls +(0, -0.6) and +(0,-0.6).. (0.5, 0)--(0.5, 0.8); 
\begin{scope}[shift={(0.5, 0.3)}]
\draw [fill=white] (-0.3, -0.3) rectangle (0.3, 0.3);
\node at (0, 0) {\tiny $v$};
\end{scope}
    \end{scope}
\draw [blue] (-0.5, -0.8)--(-0.5, 0) .. controls +(0, 0.6) and +(0,0.6).. (0.5, 0)--(0.5, -0.8);
\begin{scope}[shift={(0.5, -0.3)}]
\draw [fill=white] (-0.3, -0.3) rectangle (0.3, 0.3);
\node at (0, 0) {\tiny $v$};
\end{scope}
\end{tikzpicture}}},\,
F_{2,2}=\lambda^{1/2}\,\vcenter{\hbox{\begin{tikzpicture}[scale=0.65]
    \begin{scope}[shift={(0,1.5)}]
    \draw [blue] (-0.5, 0.8)--(-0.5, 0) .. controls +(0, -0.6) and +(0,-0.6).. (0.5, 0)--(0.5, 0.8); 
\begin{scope}[shift={(0.5, 0.3)}]
\draw [fill=white] (-0.3, -0.3) rectangle (0.3, 0.3);
\node at (0, 0) {\tiny $v^*$};
\end{scope}
    \end{scope}
\draw [blue] (-0.5, -0.8)--(-0.5, 0) .. controls +(0, 0.6) and +(0,0.6).. (0.5, 0)--(0.5, -0.8);
\begin{scope}[shift={(0.5, -0.3)}]
\draw [fill=white] (-0.3, -0.3) rectangle (0.3, 0.3);
\node at (0, 0) {\tiny $v$};
\end{scope}
\end{tikzpicture}}}.
\end{align*}
We have that $\displaystyle T=\begin{pmatrix}
    \frac{1}{2}(\mu+\mu^{-1}) F_{1,1} &  (a-b)F_{1,2} \\ (a-b)F_{2,1} & \frac{1}{2}(\mu+\mu^{-1}) F_{2,2}
\end{pmatrix}$.
Suppose that $D=\begin{pmatrix}
        c & 0 \\
        0 & 1-c
    \end{pmatrix}$.
Note that the terms $\begin{pmatrix} D & \\ & D\end{pmatrix}$ commutes with the other terms.
It is reasonable to assume that $D$ is diagonal to estimate $\beta$.
The linear map $\K_{D, \Pi}$ can be written in matrix form as follows
\begin{align}
    K_{D, \Pi} = \int_0^1 A(s) \otimes B(s) ds,
\end{align}
where
\begin{align*}
A(s) &= \begin{pmatrix} \mu^{2s-1} & \gamma \\ \gamma & \mu^{1-2s} \end{pmatrix} \otimes D^s = M(s) \otimes D^s,\quad 
B(s)= I_2 \otimes D^{1-s}.
\end{align*}
Let $\gamma=b-a$ and 
\begin{align*}
h_{11}=& \int_0^1 \mu^{2s-1} c^s (1-c)^{1-s} ds, \\
h_{22}=& \int_0^1 \mu^{1-2s}) c^s (1-c)^{1-s} ds, \\
h_{12}=h_{21}=& \gamma \int_0^1 c^s (1-c)^{1-s} ds.
\end{align*}
We obtain that the minimial eigenvalues of $K_{D, \Pi}^{1/2} TK_{D, \Pi}^{-1/2} + K_{D, \Pi}^{-1/2} TK_{D, \Pi}^{1/2}$ is 
\begin{align*}
\lambda_{\min}(c) = (\mu + \mu^{-1}) - 2|\gamma| \sqrt{ 1 + \varphi(c) },
\end{align*}
where
$$
\varphi(c) = \frac{(h_{11} - h_{22})^2}{4(h_{11} h_{22} - h_{12}^2)}.
$$
\begin{proposition}\label{prop:fermbound1}
The global minimum eigenvalue $\lambda_{\min}(c)$ is monotonically increasing on $c \in (0, 1/2)$.
\end{proposition}
\begin{proof}
It suffices to show that the term $\varphi(c)$ is monotonically decreasing with respect to $c$. 
Let $\displaystyle \epsilon = \frac{1-c}{c}$. 
When $c \in (0, 1/2)$, $\epsilon \in (1, +\infty)$. 
As $c$ increases, $\epsilon$ monotonically decreases. 
Therefore, we have to show that $\varphi(\epsilon )$ is a monotonically increasing function with respect to $\epsilon $.
Let $\displaystyle s = \frac{1}{2} + t$.
We have that 
\begin{align*}
  h_{11}-h_{22}=&  \int_0^1 (\mu^{2s-1} - \mu^{1-2s}) c^s (1-c)^{1-s} ds \\
  =& 2\sqrt{c(1-c)} \epsilon ^{1/2} \int_{-1/2}^{1/2} \sinh(2t \ln \mu) \epsilon ^{-t} dt \\
  =& 2\sqrt{c(1-c)} \epsilon ^{1/2} \int_0^{1/2} \sinh(2t \ln \mu) \left( \epsilon ^{-t} - \epsilon ^t \right) dt,
\end{align*}
where the third equality hold since $\sinh(2t \ln \mu)$ is an odd function of $t$.
Noting that $\epsilon ^{-t} - \epsilon ^t = -2\sinh(t \ln \epsilon)$, we arrive at an  integral expression:
$$
(h_{11} - h_{22})^2 = 16 c(1-c) \epsilon  \left( \int_0^{1/2} \sinh(2t \ln \mu) \sinh(t \ln \epsilon ) dt \right)^2.
$$
Let 
\begin{align*}
\Upsilon_1(\epsilon)=& 2\int_0^{1/2} \cosh(2t \ln \mu) \cosh(t \ln \epsilon ) dt,  \\
\Upsilon_2(\epsilon)=& 2\int_0^{1/2} \sinh(2t \ln \mu) \sinh(t \ln \epsilon ) dt , \\
\Upsilon_3(\epsilon)=& 2\int_0^{1/2}  \cosh(t \ln \epsilon ) dt.
\end{align*}
Similarly, we can write the cross terms and volume terms in the denominator as similar integrals. Therefore, $\varphi(\epsilon )$ can essentially be expressed as:
$$
\varphi(\epsilon ) = \frac{\Upsilon_2^2}{\Upsilon_1^2-\Upsilon_2^2-\gamma^2 \Upsilon_3^2}.
$$
By taking differentiation, we see that $\displaystyle \frac{\Upsilon_1}{\Upsilon_2}$ is decreasing and $\displaystyle \frac{\Upsilon_3}{\Upsilon_2}$ is increasing with respect to $\epsilon$.
This shows that $\varphi(\epsilon)$ is increasing with respect to $\epsilon$.

When $c \to 1/2$ (i.e., $\epsilon  \to 1$), we have that $\ln \epsilon = 0$, leading to $\sinh(t \ln \epsilon) = 0$. 
The numerator integral is $0$, and thus $\varphi(1) = 0$. 
In this case, there is no non-commutative penalty, and $\lambda_{\min}$ reaches its maximum value $(\mu+\mu^{-1}) - 2|\gamma|$.

\end{proof}

\begin{proposition}\label{prop:fermbound2}
For any $\mu>0$ and $|\gamma|<1$, we have that 
$$
\lim_{c \to 0} \lambda_{\min} = (\mu + \mu^{-1}) - |\gamma| \sqrt{ \frac{(\mu + \mu^{-1})^2 - 4\gamma^2}{1 - \gamma^2} }.
$$
When $\displaystyle |\gamma| \leq  \frac{\sqrt{2}}{2}$, we have that 
$$
\inf_{c \in (0,1/2)} \lambda_{\min} = \lim_{c \to 0} \lambda_{\min} = (\mu + \mu^{-1}) - \sqrt{\mu^2+\mu^{-2}}>0.
$$
\end{proposition}

\begin{proof}
When $c \to 0$, $1-c \to 1$, the integral is essentially dominated by the $c^s$ term.
All logarithmic mean denominators contain $\ln c - \ln(1-c) \approx \ln c$. 
Since $\ln c \to -\infty$ as $c \to 0^+$, we extract the common infinitesimal dominant factor $\displaystyle \tilde{c} = \frac{-1}{\ln c}$.

We perform asymptotic expansion on the three core integral quantities:
\begin{align*}
    h_{12} = |\gamma| \frac{c - (1-c)}{\ln c - \ln(1-c)} \approx |\gamma| \frac{-1}{\ln c},\\
   h_{11} = \mu^{-1}(1-c) \frac{\frac{\mu^2 c}{1-c} - 1}{\ln\left(\frac{\mu^2 c}{1-c}\right)} \approx \mu^{-1} \frac{-1}{\ln c},\\
   h_{22} = \mu(1-c) \frac{\frac{\mu^{-2} c}{1-c} - 1}{\ln\left(\frac{\mu^{-2} c}{1-c}\right)} \approx \mu\frac{-1}{\ln c} .
\end{align*}
Hence
$$
\varphi\approx \frac{(\mu^{-1}\tilde{c} - \mu \tilde{c})^2}{4(\mu^{-1}\tilde{c} \cdot \mu \tilde{c} - |\gamma|^2 T^2)} = \frac{\tilde{c}^2 (\mu^{-1} - \mu)^2}{4 \tilde{c}^2 (1 - \gamma^2)}
= \frac{(\mu - \mu^{-1})^2}{4(1 - \gamma^2)}.
$$
Thus
\begin{align*}
    \lim_{c \to 0} \lambda_{\min} &= (\mu + \mu^{-1}) - 2|\gamma| \sqrt{ 1 + \frac{(\mu - \mu^{-1})^2}{4(1 - \gamma^2)} }\\
    &=(\mu + \mu^{-1}) - |\gamma| \sqrt{ \frac{(\mu + \mu^{-1})^2 - 4\gamma^2}{1 - \gamma^2} }.
\end{align*}
By taking $\displaystyle \gamma^2=\frac{1}{2}$, we see the second equality is true.
\end{proof}

\begin{theorem}[Generalized Logarithmic Sobolev Inequality]
Suppose that $\displaystyle |\gamma|\leq \frac{\sqrt{2}}{2}$ and $\displaystyle \beta=\frac{(\mu + \mu^{-1}) - \sqrt{\mu^2+\mu^{-2}}}{2}$.
We have that 
\begin{align*}
 H(\Phi_t^*(D)\|\rho)\leq   e^{-2\beta t} H(D\| \rho). 
\end{align*}
\end{theorem}
\begin{proof}
It follows from Proposition \ref{prop:fermbound2} and Theorem \ref{thm:entropydecay}.    
\end{proof}

\begin{remark}
The constant $\beta$ depends on the density matrix $D$.
Suppose that $\displaystyle c=\frac{1}{3}$.
We have that 
\begin{align*}
 \sup_{\mu>0, |\gamma|<1} \sqrt{1+\varphi(1/3)}=\frac{3\sqrt{2}}{4}\approx 1.0606.
\end{align*}
This implies that $\beta=\lambda_{\min}>0$ for any $|\gamma|<1$ when $\mu$ large enough and the generalized logarithmic Sobolev inequality holds.
\end{remark}

\begin{remark}
The method can help us construct KMS symmetric semigroups for any finite dimensional von Neumann algebra $\cA$ described in Section 3.
For instance, we can take
\begin{align*}
\widehat{\Delta}
=e_2 +  \sum_{j=1}^m e^{\beta a_j/4} \lambda^{1/2} \vcenter{\hbox{\begin{tikzpicture}[scale=0.65]
    \begin{scope}[shift={(0,1.5)}]
    \draw [blue] (-0.5, 0.8)--(-0.5, 0) .. controls +(0, -0.6) and +(0,-0.6).. (0.5, 0)--(0.5, 0.8);    
\begin{scope}[shift={(0.5, 0.3)}]
\draw [fill=white] (-0.3, -0.3) rectangle (0.3, 0.3);
\node at (0, 0) {\tiny $v_j$};
\end{scope}
    \end{scope}
\draw [blue] (-0.5, -0.8)--(-0.5, 0) .. controls +(0, 0.6) and +(0,0.6).. (0.5, 0)--(0.5, -0.8);
\begin{scope}[shift={(0.5, -0.3)}]
\draw [fill=white] (-0.3, -0.3) rectangle (0.3, 0.3);
\node at (0, 0) {\tiny $v_j^*$};
\end{scope}
\end{tikzpicture}}}
+ \sum_{j=1}^m e^{-\beta a_j/4} \lambda^{1/2} \vcenter{\hbox{\begin{tikzpicture}[scale=0.65]
    \begin{scope}[shift={(0,1.5)}]
    \draw [blue] (-0.5, 0.8)--(-0.5, 0) .. controls +(0, -0.6) and +(0,-0.6).. (0.5, 0)--(0.5, 0.8);    
\begin{scope}[shift={(0.5, 0.3)}]
\draw [fill=white] (-0.3, -0.3) rectangle (0.3, 0.3);
\node at (0, 0) {\tiny $v_j^*$};
\end{scope}
    \end{scope}
\draw [blue] (-0.5, -0.8)--(-0.5, 0) .. controls +(0, 0.6) and +(0,0.6).. (0.5, 0)--(0.5, -0.8);
\begin{scope}[shift={(0.5, -0.3)}]
\draw [fill=white] (-0.3, -0.3) rectangle (0.3, 0.3);
\node at (0, 0) {\tiny $v_j$};
\end{scope}
\end{tikzpicture}}}
+\left( 1- e_2- \sum_{j=1}^m \lambda^{1/2}\vcenter{\hbox{\begin{tikzpicture}[scale=0.65]
    \begin{scope}[shift={(0,1.5)}]
    \draw [blue] (-0.5, 0.8)--(-0.5, 0) .. controls +(0, -0.6) and +(0,-0.6).. (0.5, 0)--(0.5, 0.8);    
\begin{scope}[shift={(0.5, 0.3)}]
\draw [fill=white] (-0.3, -0.3) rectangle (0.3, 0.3);
\node at (0, 0) {\tiny $v_j$};
\end{scope}
    \end{scope}
\draw [blue] (-0.5, -0.8)--(-0.5, 0) .. controls +(0, 0.6) and +(0,0.6).. (0.5, 0)--(0.5, -0.8);
\begin{scope}[shift={(0.5, -0.3)}]
\draw [fill=white] (-0.3, -0.3) rectangle (0.3, 0.3);
\node at (0, 0) {\tiny $v_j^*$};
\end{scope}
\end{tikzpicture}}}
-\sum_{j=1}^m \lambda^{1/2}\vcenter{\hbox{\begin{tikzpicture}[scale=0.65]
    \begin{scope}[shift={(0,1.5)}]
    \draw [blue] (-0.5, 0.8)--(-0.5, 0) .. controls +(0, -0.6) and +(0,-0.6).. (0.5, 0)--(0.5, 0.8);    
\begin{scope}[shift={(0.5, 0.3)}]
\draw [fill=white] (-0.3, -0.3) rectangle (0.3, 0.3);
\node at (0, 0) {\tiny $v_j^*$};
\end{scope}
    \end{scope}
\draw [blue] (-0.5, -0.8)--(-0.5, 0) .. controls +(0, 0.6) and +(0,0.6).. (0.5, 0)--(0.5, -0.8);
\begin{scope}[shift={(0.5, -0.3)}]
\draw [fill=white] (-0.3, -0.3) rectangle (0.3, 0.3);
\node at (0, 0) {\tiny $v_j$};
\end{scope}
\end{tikzpicture}}}\right)
\end{align*}
be the associated bimodule modular operator, where $a_j$ are distinct for $j=1, \ldots, m$ and
\begin{align*}
H=\sum_{j,k=1}^m \gamma_{j,k}^+ \vcenter{\hbox{\begin{tikzpicture}[scale=0.65]
    \begin{scope}[shift={(0,1.5)}]
    \draw [blue] (-0.5, 0.8)--(-0.5, 0) .. controls +(0, -0.6) and +(0,-0.6).. (0.5, 0)--(0.5, 0.8);    
\begin{scope}[shift={(0.5, 0.3)}]
\draw [fill=white] (-0.5, -0.3) rectangle (0.5, 0.3);
\node at (0, 0) {\tiny $iwQ_j$};
\end{scope}
    \end{scope}
\draw [blue] (-0.5, -0.8)--(-0.5, 0) .. controls +(0, 0.6) and +(0,0.6).. (0.5, 0)--(0.5, -0.8);
\begin{scope}[shift={(0.5, -0.3)}]
\draw [fill=white] (-0.5, -0.3) rectangle (0.5, 0.3);
\node at (0, 0) {\tiny $iwQ_k$};
\end{scope}
\end{tikzpicture}}}
+ \gamma_{j,k}^- \vcenter{\hbox{\begin{tikzpicture}[scale=0.65]
    \begin{scope}[shift={(0,1.5)}]
    \draw [blue] (-0.5, 0.8)--(-0.5, 0) .. controls +(0, -0.6) and +(0,-0.6).. (0.5, 0)--(0.5, 0.8);    
\begin{scope}[shift={(0.5, 0.3)}]
\draw [fill=white] (-0.5, -0.3) rectangle (0.5, 0.3);
\node at (0, 0) {\tiny $iwP_j$};
\end{scope}
    \end{scope}
\draw [blue] (-0.5, -0.8)--(-0.5, 0) .. controls +(0, 0.6) and +(0,0.6).. (0.5, 0)--(0.5, -0.8);
\begin{scope}[shift={(0.5, -0.3)}]
\draw [fill=white] (-0.5, -0.3) rectangle (0.5, 0.3);
\node at (0, 0) {\tiny $iwP_k$};
\end{scope}
\end{tikzpicture}}},
\end{align*}
where $(\gamma_{j,k}^{\pm})_{j,k=1}^m$ are real positive definite matrices such that $\gamma_{j,j}^++\gamma_{j,j}^{-}=1$ and $\gamma_{j,j}^+ \neq \gamma_{j,j}^-$ for all $j=1, \ldots, m$ to obtain KMS symmetric quantum Markov semigroups.
By certain constrains on the matrices $(\gamma_{j,k}^{\pm})_{j,k=1}^m$, the algebra $\cA$ can be any finite dimensional von Neumann algebra.
\end{remark}

\bibliographystyle{abbrv}
\bibliography{kms}

\end{document}